\documentclass[11 pt, reqno]{amsart}

\usepackage{setspace} 
\usepackage{geometry} 
\usepackage[dvipsnames]{xcolor}
\usepackage{amsfonts}
\usepackage{amsmath}
\usepackage{graphicx}
\usepackage{mathrsfs}
\usepackage{amssymb}
\usepackage{amsthm}
\usepackage{tikz}
\usepackage{tikz-cd}
\usepackage{color}
\usepackage{enumitem}
\usepackage{imakeidx}
\usepackage[colorinlistoftodos]{todonotes}
\usepackage{upgreek}

\usepackage[hypertexnames=false]{hyperref}
\hypersetup{colorlinks=true, linkcolor=., citecolor=ForestGreen, urlcolor=blue}

\usepackage[capitalise, nameinlink]{cleveref}
\crefformat{equation}{(#2#1#3)} 

\usepackage{mdframed}
\usepackage{manfnt}
\usepackage[new]{old-arrows}
\usepackage{changepage}
\usepackage{stmaryrd}
\usepackage{appendix}
\usepackage[normalem]{ulem}
\usepackage{blkarray}
\usepackage{mathtools}
\usepackage{verbatim}
\usepackage{wasysym}

\usepackage{cancel}

\newcommand\linkcolor{BrickRed}
\let\crefinner=\cref
\newcommand\mycref[1]{{\color{\linkcolor}\crefinner{#1}}}
\def\cref{\mycref}
\let\refinner=\ref
\newcommand\myref[1]{{\color{\linkcolor}\refinner{#1}}}
\def\ref{\myref}

\DeclareMathAlphabet{\matheuler}{U}{eus}{m}{n}

\setlist[enumerate]{align=right, label=\bfseries{(\arabic*)}} 

\newenvironment{alphabetize}{\begin{enumerate}[label=\bfseries{(\alph*)}]}{\end{enumerate}}
\newenvironment{romanize}{\begin{enumerate}[label=\bfseries{(\roman*)}]}{\end{enumerate}}

\tikzcdset{
  cells={font=\everymath\expandafter{\the\everymath\displaystyle}},
}

\newcommand\define[1]{{\color{violet}{\texttt{#1}}}}
\newcommand\important[1]{\emph{#1}}

\newcommand{\MAT}[7]{\begin{pmatrix}{#1}&{#2}&{#3}&{#4}\\{#5}&{#6}&{#7}&\MATHelper}
\newcommand{\MATHelper}[9]{{#1}\\{#2}&{#3}&{#4}&{#5}\\{#6}&{#7}&{#8}&{#9}\end{pmatrix}}
\newcommand{\Matrix}[7]{\begin{pmatrix}{#1}&{#2}&{#3}&{#4}&{#5}\\{#6}&{#7}&\MatrixHelper}
\newcommand{\MatrixHelper}[9]{{#1}&{#2}&{#3}\\{#4}&{#5}&{#6}&{#7}&{#8}\\{#9}&\MatrixHelperii}
\newcommand{\MatrixHelperii}[9]{{#1}&{#2}&{#3}&{#4}\\{#5}&{#6}&{#7}&{#8}&{#9}\end{pmatrix}}

\newcommand{\units}[1]{#1^{\times}}

\newcommand{\inv}[2][1]{#2^{-#1}}
\newcommand{\pinv}[2][1]{\inv[#1]{\parens{#2}}}

\newcommand{\ur}[1]{#1^{\mrm{ur}}}

\DeclareMathOperator{\Hom}{Hom}

\DeclareMathOperator{\coker}{coker}

\renewcommand{\hom}{\operatorname H}

\DeclareMathOperator{\im}{im}

\newcommand{\zmod}[1]{\mathbb Z/#1\mathbb Z}

\newcommand\uhom{\hom_{\mrm{ur}}}
\DeclareMathOperator\opinv{inv}
\DeclareMathOperator\oploc{loc}
\DeclareMathOperator{\Disc}{Disc}
\DeclareMathOperator{\Nm}{Nm}

\DeclareMathOperator{\ord}{ord}

\DeclareMathOperator{\disc}{disc}
\DeclareMathOperator{\Gal}{Gal}
\DeclareMathOperator{\Br}{Br}
\DeclareMathOperator{\vBr}{Br_{vert}}

\DeclareMathOperator{\rank}{rank}

\DeclareMathOperator{\Sel}{Sel}
\DeclareFontFamily{U}{wncy}{}
\DeclareFontShape{U}{wncy}{m}{n}{<->wncyr10}{}
\DeclareSymbolFont{mcy}{U}{wncy}{m}{n}
\DeclareMathSymbol{\Sha}{\mathord}{mcy}{"58} 

\newcommand{\val}{\operatorname{val}}

\newcommand{\openset}{\overset{\text{open}}\subset}

\newcommand{\Jac}{\operatorname{Jac}}

\newcommand{\abs}[1]{\left|#1\right|}

\let\templim\lim
\renewcommand{\lim}{\templim\limits}

\newcommand{\dX}[0]{\d X}
\newcommand{\dx}[0]{\d x}

\renewcommand{\d}{\mathrm d}

\newcommand\homses[4]{
    \begin{tikzcd}[ampersand replacement=\&]
        0\ar[r]\&{\displaystyle #1}\ar[r, "{#2}"]\ar[r]\&{\displaystyle #3}\ar[r, "{#4}"]\&
        \homsesHELPER
}
\newcommand\homsesHELPER[9]{
        {\displaystyle #1}\ar[r]\& 0\\
        0\ar[r]\&{\displaystyle #5}\ar[r, "{#6}"]\&{\displaystyle #7}\ar[r, "{#8}"]\&{\displaystyle #9}\ar[r]\& 0
        \ar[from=lllu, to=lll, "{#2}"]\ar[from=llu, to=ll, "{#3}"]\ar[from=lu, to=l, "{#4}"]
    \end{tikzcd}
}

\newcommand\frontthreesimplex[1]{
    \begin{tikzcd}[every arrow/.append style={dash}, ampersand replacement=\&, row sep = large]
        \&{\displaystyle#1}
        \frontthreesimplexHELPER
}
\newcommand\frontthreesimplexHELPER[9]{
        \ar[ddl, "{#1}"']\ar[d, "{#2}" description]\ar[ddr, "{#3}"]\\
        \&{\displaystyle#4}\ar[dl, "{#5}"]\ar[dr, "{#6}"']\\
        {\displaystyle#7}\ar[rr, "{#8}"']\&\&{\displaystyle#9}
    \end{tikzcd}
}

\renewcommand\d{\mrm d}
\renewcommand{\phi}{\varphi}

\newcommand{\mbf}{\mathbf}

\newcommand{\mbfv}{\mathbf v}

\newcommand{\msO}{\mathscr O}

\newcommand{\mf}{\mathfrak}
\newcommand{\mfc}{\mf c}

\newcommand{\msP}{\mathscr P}

\newcommand{\mc}{\mathcal} 

\newcommand{\mcE}{\mc E}

\newcommand{\mcM}{\mc M}

\newcommand{\mcL}{\mc L}

\newcommand{\mrm}{\mathrm}
\newcommand{\opname}{\operatorname}

\newcommand{\F}{\mathbb F}
\newcommand{\Q}{\mathbb Q}
\newcommand{\Z}{\mathbb Z}
\newcommand{\R}{\mathbb R}
\newcommand{\C}{\mathbb C}

\newcommand{\N}{\mathbb N}

\newcommand{\A}{\mathbb A}

\newcommand{\eps}{\varepsilon}

\renewcommand{\tau}{\uptau}
\renewcommand{\P}{\mathbb P}

\renewcommand{\b}{\beta}

\newcommand{\parens}[1]{\!\left(#1\right)}
\newcommand{\p}[1]{\!\left(#1\right)} 
\newcommand{\pfrac}[2]{\parens{\frac{#1}{#2}}}
\newcommand{\brackets}[1]{\left\{#1\right\}}
\newcommand{\bracks}[1]{\brackets{#1}}
\renewcommand{\b}{\bracks} 
\newcommand{\sqbracks}[1]{\!\left[#1\right]} 
\newcommand{\sq}{\sqbracks}
\newcommand{\angles}[1]{\left\langle#1\right\rangle}

\newcommand{\dparens}[1]{\!\left(\!\left(#1\right)\!\right)}
\newcommand\dps\dparens 

\newcommand{\xto}{\xrightarrow}

\newcommand{\too}{\longrightarrow}
\newcommand{\xtoo}{\xlongrightarrow}
\newcommand{\iso}{\xto\sim}

\newcommand{\xiso}[1]{\xto[#1]\sim}

\newcommand{\into}{\hookrightarrow}
\newcommand{\xinto}[1]{\overset{#1}\into}

\newcommand{\bij}{\leftrightarrow}

\makeatletter
\newcommand*{\da@rightarrow}{\mathchar"0\hexnumber@\symAMSa 4B }
\newcommand*{\da@leftarrow}{\mathchar"0\hexnumber@\symAMSa 4C }
\newcommand*{\xdashrightarrow}[2][]{%
  \mathrel{%
    \mathpalette{\da@xarrow{#1}{#2}{}\da@rightarrow{\,}{}}{}%
  }%
}
\newcommand{\xdashleftarrow}[2][]{%
  \mathrel{%
    \mathpalette{\da@xarrow{#1}{#2}\da@leftarrow{}{}{\,}}{}%
  }%
}
\newcommand*{\da@xarrow}[7]{%
  \sbox0{$\ifx#7\scriptstyle\scriptscriptstyle\else\scriptstyle\fi#5#1#6\m@th$}%
  \sbox2{$\ifx#7\scriptstyle\scriptscriptstyle\else\scriptstyle\fi#5#2#6\m@th$}%
  \sbox4{$#7\dabar@\m@th$}%
  \dimen@=\wd0 %
  \ifdim\wd2 >\dimen@
    \dimen@=\wd2 %
  \fi
  \count@=2 %
  \def\da@bars{\dabar@\dabar@}%
  \@whiledim\count@\wd4<\dimen@\do{%
    \advance\count@\@ne
    \expandafter\def\expandafter\da@bars\expandafter{%
      \da@bars
      \dabar@ 
    }%
  }%
  \mathrel{#3}%
  \mathrel{%
    \mathop{\da@bars}\limits
    \ifx\\#1\\%
    \else
      _{\copy0}%
    \fi
    \ifx\\#2\\%
    \else
      ^{\copy2}%
    \fi
  }%
  \mathrel{#4}%
}
\makeatother
\newcommand\xdashto\xdashrightarrow
\newcommand\xdashfrom\xdashleftarrow

\makeatletter
\providecommand{\leftsquigarrow}{%
  \mathrel{\mathpalette\reflect@squig\relax}%
}
\newcommand{\reflect@squig}[2]{%
  \reflectbox{$\m@th#1\rightsquigarrow$}%
}
\makeatother

\newcommand\notsquare{\cancel\square}

\newcommand*{\rom}[1]{\textup{\uppercase\expandafter{\romannumeral#1}}}
\newcommand{\tbf}{\textbf}
\newcommand{\bp}[1]{\tbf{(#1)}} 
\newcommand{\sm}{\setminus}

\renewcommand{\l}{\ell}
\newcommand{\wt}{\widetilde}

\newcommand{\twocases}[4][.]{
	\begin{cases}
		\hfill\displaystyle{#2}\hfill&\text{if }\displaystyle{#3}\\
		\hfill\displaystyle{#4}\hfill&\text{otherwise#1}
	\end{cases}
}

\newcommand{\Twocases}[4]{
	\begin{cases}
		\hfill\displaystyle{#1}\hfill&\text{if }\displaystyle{#2}\\
		\hfill\displaystyle{#3}\hfill&\text{if }\displaystyle{#4}
	\end{cases}
}

\newcommand{\push}[1]{#1_*}

\newcommand{\pull}[1]{#1^*}

\newcommand{\by}{\times}
\newcommand{\mapdesc}[5]{
	\begin{matrix}
		\ifblank{#1}{}{\displaystyle#1:}&\displaystyle#2&\longrightarrow&\displaystyle#3\\
		&\displaystyle#4&\longmapsto&\displaystyle #5
	\end{matrix}
}

\DeclareMathOperator{\sign}{sign}
\renewcommand{\Re}{\operatorname{Re}}

\renewcommand{\bar}{\overline}

\newcommand\colonequals\coloneqq
\newcommand\equalscolon\eqqcolon

\renewcommand\t\text
\newcommand\ttt\texttt
\newcommand{\tand}{\,\text{ and }\,} 
\newcommand{\tor}{\,\text{ or }\,}
\newcommand{\twhen}{\,\text{ when }\,}
\newcommand{\twhere}{\,\text{ where }\,}
\newcommand{\twith}{\,\text{ with }\,}
\newcommand{\twhile}{\,\text{ while }\,}
\newcommand{\tif}{\,\text{ if }\,}

\newcommand{\tfor}{\,\text{ for }\,}

\newcommand{\tbut}{\,\text{ but }\,}
\newcommand{\tcomma}{\text{, }\,\,}

\newcommand{\tforall}{\,\text{ for all }\,}

\newcommand{\tforsome}{\,\text{ for some }\,}

\newcommand{\tandso}{\,\text{ and so }\,}

\newcommand{\tin}{\,\text{ in }\,}

\tikzset{%
	symbol/.style={%
		,draw=none
		,every to/.append style={%
			edge node={node [sloped, allow upside down, auto=false]{$#1$}}}
	}
}

\makeatletter
\def\renewtheorem#1{%
	\expandafter\let\csname#1\endcsname\relax
	\expandafter\let\csname c@#1\endcsname\relax
	\gdef\renewtheorem@envname{#1}
	\renewtheorem@secpar
}
\def\renewtheorem@secpar{\@ifnextchar[{\renewtheorem@numberedlike}{\renewtheorem@nonumberedlike}}
\def\renewtheorem@numberedlike[#1]#2{\newtheorem{\renewtheorem@envname}[#1]{#2}}
\def\renewtheorem@nonumberedlike#1{  
	\def\renewtheorem@caption{#1}
	\edef\renewtheorem@nowithin{\noexpand\newtheorem{\renewtheorem@envname}{\renewtheorem@caption}}
	\renewtheorem@thirdpar
}
\def\renewtheorem@thirdpar{\@ifnextchar[{\renewtheorem@within}{\renewtheorem@nowithin}}
\def\renewtheorem@within[#1]{\renewtheorem@nowithin[#1]}
\makeatother

\DeclareDocumentEnvironment{MyFrame}{O{1cm}O{0.4pt}O{0.8cm}O{black}O{3}O{2ex}}
{\par\hfill\rlap{%
		\bgroup\color{#4}%
		\hskip-\dimexpr#1-#3\relax\rule{#1}{#2}%
		\hskip-\dimexpr#1/#5\relax\rule[-\dimexpr#1-\dimexpr#1/#5\relax]{#2}{#1}%
		\egroup
	}%
	\vskip-\dimexpr#1/#5+\dimexpr#1/#5-#6\relax%
}
{\par\nobreak\offinterlineskip\vskip-\dimexpr#1/#5+\dimexpr#1/#5-#6\relax\noindent%
	\hskip-#3\bgroup\color{#4}%
	\rule{#1}{#2}\hskip-\dimexpr#1-\dimexpr#1/#5-#2\relax%
	\rule[-\dimexpr#1/#5-#2\relax]{#2}{#1}\egroup\par
}

\newtheorem{innercustomgeneric}{\customgenericname}
\providecommand{\customgenericname}{}
\newcommand{\newcustomtheorem}[2]{%
	\newenvironment{#1}[1]
	{%
		\renewcommand\customgenericname{#2}%
		\renewcommand\theinnercustomgeneric{##1}%
		\innercustomgeneric
	}
	{\endinnercustomgeneric}
}

\makeatletter
\@ifclassloaded{amsart}{
}{%
    \@ifundefined{subjclass}{
    \newcommand{\subjclass}[2][1991]{%
      \let\@oldtitle\@title%
      \gdef\@title{\@oldtitle\footnotetext{#1 \emph{Mathematics subject classification.} #2}}%
    }
    }{}
    \@ifundefined{keywords}{
    \newcommand{\keywords}[1]{%
      \let\@@oldtitle\@title%
      \gdef\@title{\@@oldtitle\footnotetext{\emph{Key words and phrases.} #1.}}%
    }
    }{}
}
\makeatother

\theoremstyle{plain}
\newtheorem{thm}{Theorem}
\newtheorem{thma}{Theorem}
\newtheorem{lemma}[thm]{Lemma}
\newtheorem{cor}[thm]{Corollary}

\renewtheorem{conj}[thm]{Conjecture}

\newtheorem{prop}[thm]{Proposition}

\newcustomtheorem{Prob}{Problem}
\newcustomtheorem{Exc}{Theoretical Exercise}

\theoremstyle{definition}
\newtheorem{defninner}[thm]{Definition}

\newtheorem{notn}[thm]{Notation}

\newtheorem{recinner}[thm]{Recall}
\newtheorem{setinner}[thm]{Setup}
\newtheorem{exinner}[thm]{Example}
\newtheorem{reminner}[thm]{Remark}

\newtheorem{constructinner}[thm]{Construction}

\newtheorem*{exampleinner}{Example}

\newtheorem*{nonexinner}{Non-example}
\newtheorem{warninner}[thm]{Warning}

\newtheorem*{ansinner}{Answer}

\newtheorem*{histinner}{History}

\newtheorem{conv}[thm]{Convention}

\theoremstyle{remark}

\crefname{thm}{Theorem}{Theorems}
\crefname{thma}{Theorem}{Theorems}
\crefname{prop}{Proposition}{Propositions}
\crefname{lemma}{Lemma}{Lemmas}
\crefname{cor}{Corollary}{Corollaries}
\crefname{prob}{Problem}{Problems}
\crefname{exinner}{Example}{Examples}
\crefname{reminner}{Remark}{Remarks}
\crefname{ansinner}{Answer}{Answers}
\crefname{warninner}{Warning}{Warnings}
\crefname{nonexinner}{Non-example}{Non-examples}
\crefname{recinner}{Recall}{Recall}
\crefname{defninner}{Definition}{Definitions}
\crefname{histinner}{History}{History}
\crefname{constructinner}{Construction}{Constructions}

\newcommand\exsymbol{$\triangle$}
\newcommand\remsymbol{$\circ$}
\newcommand\anssymbol{$\star$}
\newcommand\warnsymbol{$\bullet$}
\newcommand\proofsymbol{$\blacksquare$}
\newcommand\nonexsymbol{$\triangledown$}
\newcommand\recsymbol{$\odot$}
\newcommand\defnsymbol{$\diamond$}
\newcommand\histsymbol{$\ominus$}
\newcommand\constructsymbol{$\octagon$}
\newcommand\setupsymbol{$\circledcirc$}

\newenvironment{ex}[1][]{\begin{exinner}[#1]\ifdefined\exsymbol\pushQED{\qed}\renewcommand\qedsymbol\exsymbol\fi}{\ifdefined\exsymbol\popQED\fi\end{exinner}}
\newenvironment{example}[1][]{\begin{exampleinner}[#1]\ifdefined\exsymbol\pushQED{\qed}\renewcommand\qedsymbol\exsymbol\fi}{\ifdefined\exsymbol\popQED\fi\end{exampleinner}}
\newenvironment{rem}[1][]{\begin{reminner}[#1]\ifdefined\remsymbol\pushQED{\qed}\renewcommand\qedsymbol\remsymbol\fi}{\ifdefined\remsymbol\popQED\fi\end{reminner}}

\newenvironment{warn}[1][]{\begin{warninner}[#1]\ifdefined\warnsymbol\pushQED{\qed}\renewcommand\qedsymbol\warnsymbol\fi}{\ifdefined\warnsymbol\popQED\fi\end{warninner}}

\newenvironment{rec}[1][]{\begin{recinner}[#1]\ifdefined\recsymbol\pushQED{\qed}\renewcommand\qedsymbol\recsymbol\fi}{\ifdefined\recsymbol\popQED\fi\end{recinner}}
\newenvironment{defn}[1][]{\begin{defninner}[#1]\ifdefined\defnsymbol\pushQED{\qed}\renewcommand\qedsymbol\defnsymbol\fi}{\ifdefined\defnsymbol\popQED\fi\end{defninner}}

\newenvironment{construct}[1][]{\begin{constructinner}[#1]\ifdefined\constructsymbol\pushQED{\qed}\renewcommand\qedsymbol\constructsymbol\fi}{\ifdefined\constructsymbol\popQED\fi\end{constructinner}}
\newenvironment{set}[1][]{\begin{setinner}[#1]\ifdefined\setupsymbol\pushQED{\qed}\renewcommand\qedsymbol\setupsymbol\fi}{\ifdefined\setupsymbol\popQED\fi\end{setinner}}

\let\proofinner=\proof
\newcommand\myproof[1][Proof]{\proofinner[#1]\renewcommand\qedsymbol\proofsymbol}
\def\proof{\myproof}

\renewcommand{\qedsymbol}{\proofsymbol}

\usepackage[style=alphabetic,url=true,backref=true,doi=false,isbn=false,maxnames=99,maxalphanames=99]{biblatex}
\usepackage{autonum}

\usepackage{orcidlink}
\makeatletter
\newcommand{\blx@noerroretextools}{}
\makeatother

\numberwithin{thm}{section}
\numberwithin{equation}{section}

\author{Niven Achenjang \orcidlink{0000-0001-9551-5821}}
\address{Niven T. Achenjang \\
	Department of Mathematics \\
	Harvard University \\
	Cambridge, MA 02138 }
\email{achenjang@math.harvard.edu}

\author{Katharine Woo \orcidlink{0009-0005-2837-106X}}
\address{Katharine Woo \\
    ETH Institute for Theoretical Studies \\
    Scheuchzerstrasse 70 \\ 8092 Zürich
	 }
\email{k.woo@eth-its.ethz.ch}

\title[Arithmetic applications of additive combinatorics]{Prime values of linear forms and a quadratic form: \\ one proof and two arithmetic applications}

\date{\today}

\allowdisplaybreaks

\let\exsymbol\undefined
\let\remsymbol\undefined
\let\anssymbol\undefined
\let\warnsymbol\undefined
\let\nonexsymbol\undefined
\let\recsymbol\undefined
\let\defnsymbol\undefined
\let\histsymbol\undefined
\let\constructsymbol\undefined
\let\setupsymbol\undefined

\newcommand{\Cramer}{\textrm{Cram\'er}}
\newcommand{\scrP}{\mathscr{P}}
\newcommand{\calQ}{\mathcal{Q}}
\newcommand{\Inv}{\textrm{inv}}
\newcommand{\fS}{\mathfrak{S}}
\newcommand{\fp}{\mathfrak{p}}
\newcommand{\fq}{\mathfrak{q}}
\newcommand{\fa}{\mathfrak{a}}
\newcommand{\fb}{\mathfrak{b}}
\newcommand{\fd}{\mathfrak{d}}
\newcommand{\fm}{\mathfrak{m}}

\newcommand{\calI}{\mathcal{I}}

\newcommand{\bx}{\mathbf{x}}

\DeclareMathOperator*{\E}{\mathbb{E}}

\newcommand\kumtwist[2][{}]{\delta^{#2}(E^{#2}[2](\Q_{#1}))}

\theoremstyle{plain}
\newtheorem*{theorem*}{Theorem}

\begin{document}
\begin{abstract}
Let $L_1,\dots,L_k\in \Z[x,y]$ be a system of linear polynomials whose homogeneous parts are pairwise linearly independent. Let $K=\Q(\alpha)$ be an imaginary quadratic field. Utilizing the higher order Fourier analysis machinery of Green and Sawhney in their breakthrough work on the Gaussian prime conjecture, we show that if there are no local obstructions, then there are infinitely many $(x,y)\in \Z$ such that $L_1(x,y),\dots,L_k(x,y),$ and $N_{K/\Q}(x+y\alpha)$ are simultaneously prime. We further establish an asymptotic formula which matches the prediction by the multivariate Bateman-Horn conjecture. 

We present two arithmetic applications of this prime counting result. First, we show that the Brauer--Manin obstruction is the only obstruction to weak approximation on a new family of varieties, including diagonal conic bundles all of whose degenerate fibres are defined over $\Q$ except one fibre that may be defined over an imaginary quadratic field. Second, we construct elliptic curves over $\Q$ with infinitely many quadratic twists of rank exactly two.

\end{abstract}
\maketitle
\vspace{-2em}

\noindent\hrulefill
\tableofcontents
\vspace{-3em}
\noindent\hrulefill
\hypersetup{linkcolor=\linkcolor}

\section{Introduction}

Given an irreducible polynomial $f(x)\in \Z[x]$, is $f(n)$ prime infinitely often? In 1857, Bunyakovsky conjectured that under certain necessary conditions, the answer to this question was positive. In 1962, Bateman and Horn \cite{BatemanHorn} gave a quantitative estimate for how often $f(x)$ is expected to attain prime values; however, both Bunyakovsky's conjecture and the Bateman-Horn conjecture are only known when $f$ is linear.

One may extend these conjecture to systems of multivariate polynomials, where they sometimes become more amenable to study. See \cite[Conjecture A.3]{DS:multi_BH} for a heuristic explaining the following statement. 
\begin{conj}[The multivariate Bateman-Horn conjecture]\label{conj: multi BH}
    Let $\scrP$ denote the set of positive and negative primes. Let $F_1,\dots,F_k\in \Z[X_1,\dots,X_n]$ be distinct nonconstant irreducible polynomials. As $X\rightarrow\infty$, 
    \begin{equation}
        \#\{\bx\in \Z^n: \|\bx\|_\infty\leq X, \forall i\in \{1,\dots,k\}, F_i(\bx)\in \scrP\} \sim \fS_{F_1,\dots,F_k} \int_{\|\bx\|_\infty\leq X} \frac{d\bx}{\prod_{i=1}^k \log^+ F_i(\bx)}, 
    \end{equation}
    where we define
    \begin{equation}
        \log^+(n) = \begin{cases}
            \log(|n|), & |n|>2, \\
            \infty, & \text{otherwise},
        \end{cases}
    \end{equation}
    and the singular series $\fS_{F_1,\dots,F_k}$\footnote{It is shown in \cite[Appendix A]{DS:multi_BH} that the singular series converges; however, it is not guaranteed to converge absolutely.} as
    \begin{equation}\label{def: singular series}
        \fS_{F_1,\dots,F_k} = \prod_{p} \frac{p^{-n}\#\{\bx\in \F_p^n: F_1\dots F_k(\bx)\in \F_p^\times\}}{(1-p^{-1})^k}.
    \end{equation}
\end{conj}
In the setting when $n$ is exponentially large compared to the degrees of $F_1,\dots,F_k$, the circle method is used to establish instances of Conjecture \ref{conj: multi BH} (see \cite{Birch, DS:multi_BH}). When the polynomial $F(\bx)\in \Z[\bx]$ has certain algebraic structure (in particular, when it is related to the norm of a number field), algebraic number theory and sieve methods combine to give some striking results (see the breakthrough works of \cite{FriedlanderIwaniec24,FouvryIwaniec,HeathBrownx32y3} and the generalizations by \cite{HBMoroz,HBLi,LSX:quadratic-prime-var,MaynardNormForm,Merikoski,Xiao:prime-quadratic}). The determinant polynomials and other polynomials invariant under group actions were studied using methods from homogeneous dynamics (see \cite{KotsovolisWoo}). We make no claim that the list of known instances of the multivariate Bateman-Horn conjecture given above and below is complete, only a survey of the land.

In this paper, we focus on the story of using additive combinatorics to study Conjecture \ref{conj: multi BH}. From Dirichlet's theorem on arithmetic progressions, Conjecture \ref{conj: multi BH} is known for a single linear form. For systems of two and three linear forms, Conjecture \ref{conj: multi BH} can be established via the Hardy-Littlewood circle method; this is a consequence of the work of Vinogradov in \cite{Vinogradov}. In 2004, Green and Tao \cite{GT1} established that the primes contain arbitrarily long arithmetic progressions using methods from additive combinatorics. This breakthrough was soon followed by their work on Conjecture \ref{conj: multi BH} for systems of linear polynomials in \cite{GT2,GT3,GT4}; the asymptotic for Conjecture \ref{conj: multi BH} for systems of linear polynomials with finite complexity was established by Green, Tao, and Ziegler in \cite{GTZ:inverse_theorem}. In all of the aforementioned results for more than three linear forms, higher order Fourier analysis is a key ingredient. These methods were then used to consider other systems of polynomials in \cite{Tao:GaussianPrimesConstellations, TZ:polynomial-progressions}. We also note the result of Frantzikinakis and Host \cite{FH:HOFA-mult-func}, which used higher order Fourier analysis to study completely multiplicative and aperiodic functions along products of linear forms and one quadratic form.

Recently, Green and Sawhney \cite{Green_Sawhney} showed Conjecture \ref{conj: multi BH} for the system of polynomials $x,y,$ and $x^2+ny^2$ (for $n\geq 1$), resolving the Gaussian Primes Conjecture of Friedlander and Iwaniec \cite{FI:coordinates_of_Gaussian_primes}.  To do so, Green and Sawhney utilized tools from sieve theory, as well as recent developments in higher order Fourier analysis on concatenation theorems by Kuca \cite{Kuca} and Kravitz--Kuca--Leng \cite{KravitzKucaLeng}, quantitative Gowers inverse theorems by Leng--Sah--Sawhney \cite{LengSahSawhney}, and quantitative bounds on Gowers norms of the von Mangoldt function by Leng \cite{Leng}. We utilize the innovative techniques of Green and Sawhney to establish the following generalization of their theorem. 

\begin{thma}\label{thm: prime count}
    Let $L_1,\dots,L_k$ be linear polynomials in $x,y$ such that the homogeneous parts of $L_1,\dots,L_k$ are pairwise linearly independent. Let $T_1,T_2 \in \Z[t]$ be linear polynomials with leading coefficients $t_1,t_2$. Let $K=\Q(\alpha)$ denote an imaginary quadratic field, with $\alpha$ an integral element, and define the quadratic polynomial: 
    \begin{equation}
        Q(x,y) = N_{K/\Q}(T_1(x)+T_2(y)\alpha).
    \end{equation}
    Let $\rho_{Q} = 2\pi [\msO_K:\Z[\alpha]]^{-1} \disc(K)^{-1/2} |t_{1}t_2|^{-1}.$
    Then as $X\rightarrow\infty$, 
    \begin{equation}
        \sum_{\substack{x,y\in \Z \\ Q(x,y)\leq X}} \prod_{i=1}^k \Lambda(L_i(x,y)) \mathbf{1}_{\scrP}(Q(x,y)) = \fS_{L_1,\dots,L_k,Q}' \cdot \frac{\rho_Q X}{\log X} + O_{L_1,\dots,L_k,T_1,T_2,\alpha}\left(\frac{X(\log\log X)^2}{(\log X)^{2}}\right),
    \end{equation}
    where $\Lambda(\cdot)$ denotes the von Mangoldt function
    and the truncated singular series is defined as 
    \begin{equation}\label{eq: truncated singular series}
        \fS'_{L_1,\dots,L_k,Q} = \prod_{p\leq \exp((\log X)^{1/10k})} \frac{p^{-n} \#\{\bx\in \F_p^2: L_1\dots L_k(\bx) Q(\bx)\in \F_p^\times\}}{(1-p^{-1})^{k+1}}. 
    \end{equation}
    Moreover, if $\fS_{L_1,\dots,L_k,Q}$ (as defined in \eqref{def: singular series}) is strictly positive, then there are infinitely many $(x,y)\in \Z$ such that $L_1(x,y),\dots,L_k(x,y)$, and $Q(x,y)$ are simultaneously prime.
\end{thma}
Observe that $\rho_Q \cdot X$ is the asymptotic area of the ellipse $Q(x,y)\leq X$ as $X\rightarrow\infty$. By the convergence of $\fS_{L_1,\dots,L_k,Q}$, Conjecture \ref{conj: multi BH} for $L_1,\dots,L_k,Q$ is a consequence of Theorem \ref{thm: prime count}.
We remark that we anticipate no major
difficulties to generalize the above theorem statement to all nondegenerate binary quadratic forms (including indefinite binary quadratic forms); Theorem \ref{thm: prime count} is, however, sufficient for our two arithmetic applications and reduces technical details in the proof.

\begin{rem}
    We remark that Theorem \ref{thm: prime count} is a joint version of the theorem of Green, Tao, and Ziegler in \cite{GTZ:inverse_theorem} and of Green and Sawhney in \cite{Green_Sawhney}. One immediate consequence in analytic number theory is towards Sarnak's saturation problem (see \cite{NS:saturation} for a definition of the saturation number and \cite{Sarnak:eqi-and-primes} for further discussion on the problem). Let $(x,y,z)$ be a Pythagorean triple. The work of Green and Tao \cite{GT2} showed that the (modified) area of the triangle (i.e. $xy/12$) has saturation number exactly 4. Theorem \ref{thm: prime count} implies that the (modified) product of the coordinates $xyz/60$ has saturation number exactly 5. Since we can write $x=2uv$, $y=u^2-v^2$, and $z = u^2+v^2$, we know that $xyz = 2uv(u+v)(u-v)(u^2+v^2)$ is a product of linear forms and one quadratic form. Previously, the best known bound for the saturation number for the (modified) product of the coordinates was $\leq 18$ by Marasingha \cite{Marasingha:saturation} (and if $v=1$ then Diamond and Halberstam \cite{DH} show that $r\leq 17$).
\end{rem}


Our arithmetic applications follow a long history of applying Conjecture \ref{conj: multi BH} towards the study of rational points on varieties, only a piece of which is summarized here. Colliot-Th\'el\`ene and Sansuc \cite{CT-Sansuc:schinzel-hasse} showed that Conjecture \ref{conj: multi BH} implies the Hasse principle and weak approximation for all conic bundles of the form $X\colon y^2-az^2=P(x)$ where $a\in\units\Q$ and $P(x)\in\Q[x]$ is irreducible. This idea was later expanded upon by \cite{Serre:Galois_cohomology, SD:pencils-of-conics, CT-SD:schinzel-hasse-severi-brauer, CTSSD:generalize-CTSD} where they established that Conjecture \ref{conj: multi BH} implies the sufficiency of the Brauer--Manin obstruction for certain, more general, class of varieties fibered over $\P^1$.

These results of Colliot-Th\'el\`ene, Sansuc, Skorobogatov, Swinnerton-Dyer, and Serre were conditional on the single variable case of Conjecture \ref{conj: multi BH}. Later seminal work of Green and Tao \cite{GT1} and Green, Tao, and Ziegler \cite{GT2,GTZ:inverse_theorem} gave new instances of Conjecture \ref{conj: multi BH}, which lead to unconditional consequences for the Hasse principle and weak approximation. Their work was used by Browning, Matthiesen, and Skorobogatov \cite{browning-matthiesen-skorobogatov} to prove unconditionally that the Brauer--Manin obstruction is the only obstruction to weak approximation on conic bundles \textit{all of whose degenerate fibers are defined over $\Q$}. Their work was then generalized by Harpaz, Skorobogatov, Wittenberg \cite{HSW:RationalPoints}, with further followup work by Browning and Matthiesen \cite{BM:norm-forms}, Harpaz and Wittenberg \cite{HW:fibration-method} and Harpaz, Wei, and Wittenberg \cite{HWW:fibration} applying to wider classes of varieties fibered over $\P^1$.

Even more recently, breakthrough work of Koymans and Pagano \cite{Koymans-Pagano:H10,koymans2025ellipticcurvesranknumber} has combined instances of Conjecture \ref{conj: multi BH} with descent in families of elliptic curves in order to better control the ranks of members of the family, leading to a resolution of Hilbert's 10th Problem over rings of integers of number fields (see also \cite{ABHS} for a different method that also applied additive combinatorics). In this paper, we give two arithmetic applications of our Theorem \ref{thm: prime count} along the above lines. \\

For our first application, in Part \ref{part:bmo}, we study the arithmetic of conic bundles (and other varieties fibered over $\P^1$), establishing new instances of weak approximations.
\begin{thma}[special case of Theorem \ref{thm:general-BMO}]\label{thma:BMO}
    Let $P(x)\in\Q[x]$ be a polynomial factoring as
    \[P(x)=L_1(x)\cdots L_k(x)\cdot Q(x)\twhere\deg L_i=1\tand\deg Q=2.\]
    If $K=\Q[x]/Q(x)$ is an imaginary quadratic field, then the Brauer--Manin obstruction on the conic bundle $X$ with affine equation
    \[X\colon y^2-az^2=P(x)\twith a\in\units\Q\]
    is the only obstruction to weak approximation on $X$. 
\end{thma}
Theorem \ref{thma:BMO} (and the more general Theorem \ref{thm:general-BMO}) generalizes the main result of \cite{HSW:RationalPoints}, and establishes new instances of Colliot-Th\'el\`ene's conjecture that the Brauer--Manin obstruction is the only obstruction to weak approximation on conic bundles (and rationally connected varieties more generally); see \cite[Section 14.1]{CT-Skorobogatov:BrauerGrothendieck} for a discussion of this conjecture.

\begin{rem}
    The hypotheses on the main result of \cite{HSW:RationalPoints} were weakened in later work of Harpaz and Wittenberg; see \cite[Theorem 1.6]{HW:fibration-method}. Our Theorem \ref{thm:general-BMO} is more analogous to \cite[Theorem 3.3]{HSW:RationalPoints} than to \cite[Theorem 1.6]{HW:fibration-method}. In principle, our Theorem \ref{thm:general-BMO} may be generalizable to give a version of \cite[Theorem 1.6]{HW:fibration-method} allowing for a single nonsplit fiber defined over a imaginary quadratic field, but our Theorem \ref{thm: prime count} does not prove sufficiently general instances of \cite[Conjecture 9.1]{HW:fibration-method}/\cite[Conjecture F]{HWW:fibration} to immediately give a direct analogue of \cite[Theorem 1.6]{HW:fibration-method}. However, it does give enough instances of this conjecture (and the modified \cite[Conjecture F$_+$]{HWW:fibration}) for \cite[Theorem 4.1]{HWW:fibration} to apply; see the proof of Theorem \ref{thm:general-BMO}.
\end{rem}

Our second application concerns the arithmetic of elliptic curves. For a fixed small number $r$, it has long been expected that, over every number field there should exist infinitely many elliptic curves $E$ of rank exactly $r$. This was confirmed for $r=0$ by Mazur and Rubin \cite{Mazur-Rubin:H10}, for $r=1$ independently by Zywina \cite{Zywina:Rank1_NumberFields} and Koymans and Pagano \cite{koymans2025ellipticcurvesranknumber}, for $r=2$ over $\Q$ by Zywina \cite{Zywina:Rank2} and over $\Q(i)$ by Savoie \cite{Savoie:rank2}, and most recently for $r\le4$ by Zywina \cite{Zywina:SmallRanks}. 
\begin{rem}
    This question for $r=1$ over $\Q$, i.e. the existence of infinitely many $E/\Q$ with $\rank E(\Q)=1$, has been resolved for much longer. Notably, Satg\'e \cite{Satge:calcul-de-Heegner} showed that, for odd primes $p$, the elliptic curves $E_{1,p}\colon X^3+Y^3=2pZ^3$ have rank 1 over $\Q$, when $p\equiv2\bmod9$, as do the elliptic curves $E_{2,p}\colon X^3+Y^3=2p^2Z^3$, when $p\equiv5\bmod9$. This was inspired by even earlier work of Heegner \cite{Heegner} who showed that the quartic curves $C_{3,p}\colon Y^2=-p(X^4-64)$ have infinitely many rational points, when $p\equiv7\bmod8$ (and still a prime). The curves $C_{3,p}$ have Jacobians $E_{3,p}=\Jac C_{3,p}\colon y^2=x^3+p^2x$, and it is not difficult to show, via an appropriate descent argument (analogous to but not the same as \cite[Proposition X.6.2]{silverman}), that $\rank E_{3,p}(\Q)\le1$, when $p\equiv7\bmod8$, so Heegner's curves also give infinitely many rank 1 elliptic curves over $\Q$ (since $C_{3,p}\cong E_{3,p}$ when $C_{3,p}(\Q)\neq\emptyset$).

    The authors are unaware of such rank 1 constructions over $\Q$ predating Heegner's paper, but cannot claim full knowledge of the history of this problem.
\end{rem}

The results of Koymans, Pagano, Savoie, and Zywina for $r\ge1$ all leverage a combination of descent with known instances of Conjecture \ref{conj: multi BH}. Savoie's work utilized Tao's constellation theorem in Gaussian primes \cite{Tao:GaussianPrimesConstellations}, Zywina's rank $2$ result over $\Q$ utilized Tao and Ziegler's work on primes in polynomial progressions \cite{TZ:polynomial-progressions}, and Koymans and Pagano's work -- as well as Zywina's paper constructing curves of rank up to $4$ -- utilized Kai's generalization of the Green-Tao-Ziegler theorem for number fields \cite{kai2025linearpatternsprimeelements}. Additionally, Savoie worked with quartic twists of $y^2=x^3+x$, Koymans and Pagano worked with quadratic twist families of elliptic curves, and Zywina worked with certain non-isotrivial families of elliptic curves. Thus, even over $\Q$, Zywina's work leaves open the question of whether one can find infinitely many rank $2$ elliptic curves in a quadratic twist family. In Part \ref{part:rank2}, we answer this in the affirmative.
\begin{thma}[special case of Theorem \ref{thm:rank2}; see Example \ref{ex:17}]\label{thma:EC}
    There exists infinitely many squarefree integers $d$ such that the elliptic curve
    \[E^d\colon dy^2=x(x-17)(x-184,009,700)\]
    has rank exactly $2$ over $\Q$.
\end{thma}

Related to the construction of elliptic curves of exact rank $r$, there is the separate question of constructing families of elliptic curves of rank $\ge r$. There have been many recent improvements in this direction; notably the most recent development by Elkies \cite{Elkies:rank1819} for rank 19 elliptic curves. Since our paper considers a rank $=r$ question as opposed to a rank $\ge r$ question, for a history of this problem and its foundational methods, we contend ourselves with referring to the papers \cite{Elkies:lectures,RS:ranks-of-ec}. 

Of course, the rank $\ge r$ and rank $=r$ questions are related. By a theorem of Silverman \cite[Theorem C]{Silverman:heights-specialization}, given an elliptic fibration $\mcE\to\P^1_K$, over a number field $K$, all but finitely many fibers have rank at least equal to the generic rank, i.e. $\rank\mcE(K(t))\le\rank\mcE_t(K)$ for all but finitely many $t\in\P^1(K)$. Beyond this, Silverman has conjectured \cite{Silverman:divisilibity-specialization} that in fact
\[\rank\mcE(K(t))\le\rank E_t(K)\le\rank\mcE(K(t))+1\quad\tfor100\%\t{ of }t\in\P^1(K),\]
with the proportion of $t$'s at either end expected to be accounted for by the distribution of root numbers of fibers in the family $\mcE$. Given this, the main difficulty in the rank $=r$ question over the rank $\ge r$ question is establishing an appropriate upper bound on the rank for infinitely many members of a given family $\mcE$. In Theorem \ref{thma:EC}, and most of the results cited above it, this difficulty is overcome by leveraging instances of Conjecture \ref{conj: multi BH} to search for fibers of $\mcE$ which are practically amenable to descent arguments.

In the setting of Theorem \ref{thma:BMO}, it is understood that new instances Conjecture \ref{conj: multi BH} lead to new classes of varieties for which one may study weak approximation and sufficiency of the Brauer--Manin obstruction. We hypothesize that this sort of philosophy holds also in the setting of Theorem \ref{thma:EC} as well, i.e., that new instances of Conjecture \ref{conj: multi BH} will give rise to new families of elliptic curves of generic rank $r$ for which one can show that there are infinitely many members of the family with rank exactly $r$.

\section*{Acknowledgments}

This material is based upon work supported by the National Science Foundation under Grant No. DMS-2502864 and DMS-2503296.
Any opinions, findings, and conclusions or recommendations expressed in this material are those of the author(s) and do not necessarily reflect the views of the National Science Foundation.

This research was worked on while the first author visited the second author at Stanford University and while the second author visited the first author at Harvard University. They are grateful to both universities for hosting the collaboration. The research was also conducted while both authors were junior research visitors at the CRM thematic program ``Universal Statistics in Number Theory'' in Montreal; we thank the organizers for a lovely program and working environment. 

The authors would also like to thank Mehtaab Sawhney and Peter Sarnak for conversations about the project. They are especially grateful to James Leng for explaining the ideas of concatenation and to Carlo Pagano for both pointing them towards the Rubin--Silverberg paper \cite{Rubin-Silverberg} (and so leading them to the family of elliptic curves in Proposition \ref{prop:RS-family}) and for the trick of using linear transformations in Section \ref{sect:AC-twist}. The authors also thank Bjorn Poonen and the participants of his ``informal number theory group meeting'' for suggestions simplifying the proof of Lemma \ref{lem:Lang-Weil-curve-trick}, as well as thank Olivier Wittenberg for pointing out the connections between Theorem \ref{thm:general-BMO} and \cite{HW:fibration-method, HWW:fibration} -- including the more streamlined proof presented in Part \ref{part:bmo} -- and thank Tim Browning for pointing out the application of Theorem \ref{thm: prime count} to Sarnak's saturation problem for the product of coordinates of Pythagorean triples. 
Finally, the authors thank Jean-Louis Colliot-Th\'el\`ene and Lilian Matthiesen for helpful comments on an earlier draft.

\section*{AI disclosure}
The background text (\S\ref{sec: HOFA background}) was written by the second author while experimenting with the capabilities of Gemini provided by the Stanford AI playground in March 2025 (Gemini 1.5 and Gemini 2.0). Since then, there has been no AI usage in the creation of this paper.  

\part{Prime values of linear forms and a quadratic form}

\section{Proof outline} \label{sec: HOFA outline}
Let us recall again the set up of our main Theorem \ref{thm: prime count}. Let $L_1,\dots,L_k$ be linear polynomials in $x,y$: 
\begin{equation}
    L_i(x,y) = L_i^0(x,y) + l_i,
\end{equation}
where $L_i^0(x,y)$ is a linear form. We assume that for each $i\neq j$ that $L_i^0$ and $L_j^0$ are linearly independent.

Let $T_x,T_y\in \Z[x]$ be linear polynomials:
\begin{equation}
    T_x(x) = t_{1,x} x + t_{0,x}, \quad T_y(y) = t_{1,y} y + t_{0,y},
\end{equation}
where $t_{1,x},t_{1,y}\in \Z_{\neq 0}$. Let $K=\Q(\alpha)$ denote an imaginary quadratic field and $\alpha$ is integral. Recall our quadratic polynomial is of the form: 
\begin{equation}
    Q(x,y) = N_{K/\Q}(T_x(x) + T_y(y)\alpha),
\end{equation}
and $\rho_Q = 2\pi |t_{1,x}t_{1,y}|^{-1} [\msO_K:\Z[\alpha]]^{-1}\sqrt{\disc(K)}^{-1}$.
Then our main theorem states that
\begin{equation}
    \sum_{\substack{x,y\in \Z \\ Q(x,y)\leq X}} \prod_{i=1}^k \Lambda(L_i(x,y)) \mathbf{1}_{\scrP}(Q(x,y)) = \fS' \cdot \frac{\rho_Q X}{\log X} + O_{L_1,\dots,L_k,T_x,T_y,K}\left(\frac{X(\log\log X)^2}{(\log X)^2}\right).
\end{equation}
(Note that our truncated singular series $\fS' := \fS_{L_1,\dots,L_k,Q}'$ (as defined in \eqref{eq: truncated singular series}) is allowed to be zero -- in this case we say that there is a local obstruction.) 

Before embarking on the proof of Theorem \ref{thm: prime count} (which largely follows the proof of \cite{Green_Sawhney}), we attempt to give an outline of the method for experts to read. In this outline, we will freely use notation from analytic number theory and additive combinatorics (and refer the reader to \S\ref{sec: HOFA background} for further background on the subject). We strongly encourage the reader to study \cite{Green_Sawhney} alongside our paper; it contains context and insight that will be missing in our paper, and because we have strived for brevity in the sections of our proof where the argument most closely resembles that in \cite{Green_Sawhney}. 

Define $w = \exp((\log X)^{1/10k})$. Define Cram\'er's model at $w$ as:
\begin{equation}
    \Lambda_{\Cramer} (n) = \prod_{p\leq w} \left(1-\frac{1}{p}\right)^{-1} \mathbf{1}_{\gcd(n,p)=1}.
\end{equation}
We can expand our prime counting function in terms of the main term:
\begin{equation}\label{eq: HOFA main term}
    M(X) := \sum_{\substack{x,y\in \Z \\ Q(x,y)\leq X}} \prod_{i=1}^k \Lambda_{\Cramer}(L_i(x,y)) \cdot \mathbf{1}_{\scrP}(Q(x,y))
\end{equation}
and error terms of the form:
\begin{equation}\label{eq: HOFA error terms}
    S_j(X):= \sum_{\substack{x,y\in \Z \\ Q(x,y)\leq X}} \prod_{i=1}^{j-1} \Lambda_{\Cramer}(L_i(x,y)) \cdot (\Lambda-\Lambda_{\Cramer})(L_j(x,y)) \cdot \prod_{i=j+1}^k \Lambda(L_i(x,y)) \cdot \mathbf{1}_{\scrP}(Q(x,y))
\end{equation}
for each $j\in \{1,\dots,k\}$.  
The main term $M(X)$ is treated classically and we will not expand further on it in this outline (see \S\ref{sec: primes MT} for a treatment of this term). We make the remark that since $\fS_{L_1,\dots,L_k,Q}$ (as defined in \eqref{def: singular series}) converges, the asymptotic given in Theorem \ref{thm: prime count} implies both Conjecture \ref{conj: multi BH} for this system of forms and the final statement of the theorem. We also make the observation that if $\fS_{L_1,\dots,L_k,Q}>0$ then at any value of $X$, $\fS_{L_1,\dots,L_k,Q}'>0$.

For the error terms $S_j(X),$ we follow the inventive ideas of Green and Sawhney in \cite{Green_Sawhney}. First, they (and hence we) apply the sieve  of Duke, Friedlander, and Iwaniec \cite{DFI} to bound $S_j(X)$ using Type I and Type II estimates; that is, sums of the form: 
\begin{equation}
    \sum_{\substack{\fd\subset \msO_K \\ N(\fd)\sim D}} \left|\sum_{\substack{N(\fa)\leq X \\ \fd\mid \fa\\ \fa = (x+y\alpha)}} \prod_{i=1}^j \Lambda_{\Cramer}(L_i(x,y)) (\Lambda-\Lambda_{\Cramer})(L_j(x,y)) \prod_{i=j+1}^k \Lambda(L_i(x,y))\right|,
\end{equation}
\begin{equation}
    \sum_{\substack{\fa,\fb\subset \msO_K \\ N(\fa\fb)\leq X \\ N(\fa)\sim D\\ \fa\fb = (x+y\alpha)}} \alpha_\fa \beta_\fb \prod_{i=1}^j \Lambda_{\Cramer}(L_i(x,y)) (\Lambda-\Lambda_{\Cramer})(L_j(x,y)) \prod_{i=j+1}^k \Lambda(L_i(x,y)),
\end{equation}
for appropriate ranges of $D$ (for Type I sums, $D\leq \delta^{O(1)}X^{1/2}$, and for Type II sums, $D\in [\delta^{O(1)},X^{3/8}]$). Here $\alpha_\fa,\beta_\fb$ are $1$-bounded complex sequences on ideals of $\msO_K$. We add the caveat that the expressions above are not fully factual; we've omitted the dependence on the linear polynomials $T_x$,$T_y$ for simplicity, but we assure the reader that they do not cause any additional complications. We also remark that, as in \cite{Green_Sawhney}, we need to separately bound the contribution from semprimes. This is done using fairly classical estimates on $\Lambda$ and $\Lambda_{\Cramer}$, and using the large sieve; this proof matches that in \cite{Green_Sawhney} and so we refrain from discussing it further in the outline (see \S\ref{sec: primes sieve}).

Now we apply the strategy of Green and Sawhney in \cite{Green_Sawhney} to bound both the Type I and Type II sums using Gowers norms. The recent development of Leng \cite{Leng} (building on the inverse theory of Leng, Sah, and Sawhney \cite{LengSahSawhney}) gives that for all $A>0$, 
\begin{equation}
    \|\Lambda-\Lambda_{\Cramer}\|_{U^s[N]} \ll_{s,A} N (\log N)^{-A}. 
\end{equation}
Green and Sawhney show that if either the Type I and Type II bound are $\gg (\log X)^{-O(1)}X$ above, then $\|\Lambda-\Lambda_{\Cramer}\|_{U^s[N]}$ must be too large for some $s\geq 1$, hence deriving a contradiction. We refer to Proposition \ref{prop: type I} and Proposition \ref{prop: type II} for explicit statements.  

While Green and Sawhney study when there are two forms -- $L_1(x,y) =x$, $L_2(x,y)=y$ -- their method generalizes to arbitrarily large number of linear forms after applying many more iterations of Cauchy-Schwarz (and hence studying $U^{s_k}[N]$ norms where $s_k$ grows with $k$). After Cauchy-Schwarz, Green and Sawhney (and hence we) apply the recent concatenation results of Kuca, Kravitz, and Leng \cite{KravitzKucaLeng}; our proof for the Type II estimates matches exactly their strategy. Where our proof minimally diverges from \cite{Green_Sawhney} is in the Type I estimates. Both our proof and \cite{Green_Sawhney} utilize that $\{n\in \Z: \exists x,y\in \Z, Q(x,y)=n\}$ is a set that while sparse is only sparse by a power of $\log X$; that is, 
\begin{equation}
    \#\{n\leq X: \exists x,y\in \Z, Q(x,y)=n\} \gg \frac{X}{(\log X)^c}
\end{equation}
for some $c>0$. Consequently, we can over replace $\{n\in \Z: \exists x,y\in \Z, Q(x,y)=n\}$ with $\Z$ itself up to an error of $(\log X)^{O(1)}$ -- we call this step ``linearization'' in our proof. In \cite{Green_Sawhney}, the final step for the Type I estimate is an explicit instance of concatenation; we instead use the more general statement of \cite{KravitzKucaLeng}. As such, we end up applying the ``linearization'' step  slightly earlier in the proof (in that we apply this ``linearization'' fact before studying and bounding the $L^2$-norm of the convolution of measures, rather than using the ``linearization'' inside of the proof of bound for the $L^2$ norm). 

The Type I estimates will be completed in \S\ref{sec: primes type I} and the Type II estimates will be completed in \S\ref{sec: primes type II}.

\section{Background and notation} \label{sec: HOFA background}
Let us begin our foray into the technical heart of the paper. We now set up some notation and introduce some concepts in analytic number theory and additive combinatorics. This section can be safely skipped by experts.

Recall that $K=\Q(\alpha)$ is an imaginary quadratic field; we denote the ring of integers as $\msO_K$. 
We will use the notation that Fraktured letters (i.e. $\fa,\fb,\dots$) denote ideals in $\msO_K$ and $\calI(\msO_K)$ will be the set of integral ideals. We will often use the shorthand that $N(\cdot) = N_{K/\Q}(\cdot)$. 

\begin{defn}
    The \define{von Mangoldt function} $\Lambda:\N\rightarrow \R$ is defined as 
    \begin{equation}
        \Lambda(n) := \begin{cases}
            \log(p), & n=p^k \text{ for some prime }p,\\
            0, & \text{otherwise.}
        \end{cases}
    \end{equation}
    We will not discriminate between positive and negative primes in this paper. Hence our definition will extend to $\Lambda:\Z\rightarrow \R$ with $\Lambda(n) := \Lambda(|n|)$ and $\Lambda(0)=0.$

    The von Mangoldt function for a number field $K$ is given by $\Lambda_K: \calI(\msO_K)\rightarrow \R$ where 
    \begin{equation}
        \Lambda_K(\fa) := \begin{cases}
            \log(N(\fp)), & \fa= \fp^k, \text{ for some prime ideal }\fp,\\
            0, & \text{otherwise.}
        \end{cases}
    \end{equation}
\end{defn}

We use Vinogradov notation. That is, given functions $f(x)$ and $g(x)$, we say that $f\ll g$, or equivalently $f= O(g)$, if there exists some constant $C$ such that $|f(x)|\leq C g(x)$ for all sufficiently large $x$. We say that $f\asymp g$, or equivalently $f= \Theta(g)$, if $f\ll g$ and $g\ll f$. We write that $f= o(g)$ if $f(x)/g(x)\rightarrow 0$ as $x\rightarrow\infty$. In our notation, all implicit constants will depend arbitrarily on $K$, $\alpha$, $T_x,T_y$, and $L_1,\dots,L_k$.

With this notation, we can now motivate the definition of the von Mangoldt function, which at first glance may seem unnatural. By partial summation, it is straightforward to show that the prime number theorem (resp. Landau's prime ideal theorem) is equivalent to the statement that 
\begin{equation}
    \sum_{n\leq X} \Lambda(n) = X + o(X), \quad (\text{resp. }\sum_{N(\fa)\leq X} \Lambda_K(\fa) = X + o(X)).
\end{equation}
Indeed, this is the pathway commonly used to prove the prime number theorem, since it is the Von Mangodlt function that directly ties to the Riemann zeta function. For $\Re(s)>1$, the following relationship holds: 
\begin{equation}
    \sum_{n\geq 1} \frac{\Lambda(n)}{n^s} = \frac{\zeta'(s)}{\zeta(s)}. 
\end{equation}
To evaluate our main term, we will need a prime ideal theorem along arithmetic progressions -- in other words, we need a version of Mitsui's prime ideal theorem with Hecke characters \cite{mitsui}. We will use the version stated by Green and Sawhney that records the explicit dependency on potential Siegel zeros. 

\begin{thm}[{\cite[Proposition 8.3]{Green_Sawhney}}]\label{thm: prime ideal theorem with characters}
Let $K$ be an imaginary quadratic field and $X>1$. There is at most one exceptional primitive Gr\"ossencharacter $\chi$ to a modulus $\fd\in \calI(\msO_K)$ and associated real number $\beta\in (1/2,1)$ such that the following is true for all Gr\"ossencharacters $\psi$ of $K$ with modulus $\fm$ and frequency $\ell$ satisfying $N(\fm), |\ell|\leq e^{\sqrt{\log X}}$: for a constant $c>0$,
\begin{equation}
    \sum_{N(\fa)\leq X} \Lambda_K(\fa)\psi(\fa) = X\cdot \mathbf{1}_{\psi=\psi_0,\ell=0} - \frac{X^\beta}{\beta} \mathbf{1}_{\psi\sim \chi,\ell=0} + O(Xe^{-c\sqrt{\log X}}). 
\end{equation}
Here $\psi_0$ denotes the principal character and $\psi\sim \chi$ if $\psi$ is induced from $\chi$. Finally, $\beta$ satsifes that 
\begin{equation}
    1-\beta \gg_\epsilon N(\fd)^{-\epsilon}.
\end{equation}
\end{thm}

\begin{rem}
    The number $\beta\in (1/2,1)$ above is referred to as a Siegel zero. Very recently, OpenAI \cite{OAI:Uniform-exclusion-of-Landau-Siegel-zeros-October-1-2026} has posted a proof of the non-existence of Landau-Siegel zeros that is widely believed to be true; in addition OpenAI \cite{OAI:The-Quasi-Riemann-Hypothesis-October-5-2026} posted works on a quasi-Riemann Hypothesis statement for $\zeta(s)$ -- that is, a zero-free vertical strip. Since then, the AI-assisted work of Cureton and Thorner \cite{hexagon:2610.00200v2}, building on the proof of OpenAI, contains a quasi-Riemann Hypothesis for Hecke $L$-functions that would exclude the existence of Siegel zeros in Theorem \ref{thm: prime ideal theorem with characters} and establish a power-saving error term in the asymptotic. Since the human verification of these recent developments has not been fully resolved, we refrain from making our result dependent on any of them. In our case, the error term from Theorem \ref{thm: prime ideal theorem with characters} as stated above will not be our final constraining error term, so improvements on this estimate do not immediately benefit our result.
\end{rem}

Next, we introduce the $W$-trick, which is a technique that commonly appears in applications of higher order Fourier analysis to the study of prime numbers (see \cite{Green_Tao:Linear_equations}). There are a few benefits to $W$-tricking arithmetic functions: first, when properly normalized, it estimates the behavior of arithmetic functiosn by their behavior at ``small'' moduli (which becomes a local problem); second, the $W$-tricked version of the Von Mangodlt function behaves ``pseudorandomly'' which was a requirement in earlier iterations of Gowers Inverse Theorems. This ``pseudorandomness'' condition was recently removed by the more powerful inverse theorem developed by Leng, Sah, and Sawhney \cite{LengSahSawhney}. 

For a given $w\in \N$, we define the integer $W$ as
\begin{equation}
    W = P(w) := \prod_{p \leq w} p.
\end{equation}
More generally, for a finite set $S\subset \scrP$ of primes, we take $P(S):= \prod_{p\in S} p$. We will choose
\begin{equation}
    w := \exp\left((\log X)^{1/10k}\right).
\end{equation}
The reason for this choice will be elucidated in \S\ref{sec: primes MT}. (In \cite{Green_Sawhney}, $w$ is chosen to be $\exp((\log X^{1/2})^{1/10})$ which is consistent with the choice taken in \cite{TaoTeravainen,Leng}; however, we need it to decrease with $k$. Fortunately, it is only important for \cite{TaoTeravainen,Leng} that the shape of the function of $w$ remains unchanged, so this modification brings no extra complications.)

\begin{defn}

    Define \define{Cram\'er's model at level Q} as the function $\Lambda_{\Cramer,Q}: \mathbb{Z} \rightarrow \mathbb{R}$,
    \begin{equation}
        \Lambda_{\Cramer,Q}(n) := \prod_{p \leq Q} \left(1 - \frac{1}{p}\right)^{-1} \mathbf{1}_{\gcd(n, p) = 1}
        =\twocases{\prod_{p\le Q}\pinv{1-\frac1p}}{\t{the smallest prime factor of }n\t{ is}>Q}0
        .
    \end{equation}
   
\end{defn}
We will often use the notation that $\Lambda_{\Cramer} = \Lambda_{\Cramer,w}$; however there will be instances where we have to consider other values of the level $Q$. \\ 

Unsurprisingly, we need to introduce the notion of Gowers norms for our proof. For a function $f: \mathbb{Z} \rightarrow \mathbb{C}$ and a shift $h \in \mathbb{Z}$, we define the difference operator 
\begin{equation}
    \Delta_h f(x) := f(x+h)\overline{f(x)}.
\end{equation}
Introduced by Gowers in \cite{Gowers:Gowers_norms}, the \define{$U^s[N]$-norm} of a function $f$ supported on $[-N, N]$ is defined as
\begin{equation}
    \|f\|_{U^s[N]}^{2^s} := \frac{1}{N} \sum_{x \in \mathbb{Z}} \mathbb{E}_{h_1, \dots, h_s \in [-N, N]} \Delta_{h_1} \dots \Delta_{h_s} f(x).
\end{equation}
For a more comprehensive treatment of the properties of these norms, including the Gowers--Cauchy--Schwarz inequality, we direct the reader to  \cite{Tao:HOFA_book}. 

These $U^s$-norms are strongly tied to the theory of classical Fourier analysis. Indeed, Plancherel's theorem tells us that if $f$ is instead viewed as supported on the finite group $\Z/(2N+1)\Z$, then 
\begin{equation}
    \|f\|_{U^2[\Z/(2N+1)\Z]} = \|\hat{f}\|_{\ell^4(\Z/(2N+1)\Z)}.
\end{equation}
It is because of this connection that one often applies higher order Fourier analysis when classical Fourier analytic techniques, such as the circle method, hit their limitations; for instance, the work of Green and Tao \cite{Green_Tao:Linear_equations} on primes in $k$-term arithmetic progressions when reduced to $k=3$ would be similar to Vinogradov's original proof of primes in $3$-term arithmetic progressions via the circle method.

In our analysis, we will require a variant of these norms. In her seminal work establishing quantitative bounds for sets lacking polynomial progressions, Peluse \cite{Peluse:sets_wo_poly_prog} developed a generalized Gowers norm adapted to shifts sampled from arbitrary probability measures rather than uniform intervals (\cite[Definition 2.1]{Peluse:sets_wo_poly_prog}). We adopt the precise formulation utilized by Green and Sawhney \cite[Definition 5.1]{Green_Sawhney}.

\begin{defn}[Gowers--Peluse Norm]\label{defn: gowers peluse norm}
    Let $f:\mathbb{Z}\rightarrow\mathbb{C}$ be supported on $[-N,N]$. Define the coupled difference operator: 
    \begin{equation}
        \Delta_{(h,h')} f(x) := f(x+h) \overline{f(x+h')} = \Delta_{h'-h} f(x+h).
    \end{equation}
    Given a set $\Omega = \{\mu_1,\dots,\mu_k\}$ of probability measures on $\mathbb{Q}$, we define the Gowers--Peluse norm as
    \begin{equation}
        \|f\|_{U_{\text{GP}}[N;\Omega]}^{2^k} = \|f\|_{U_{\text{GP}}[N;\mu_1,\dots,\mu_k]}^{2^k} := \frac{1}{N} \sum_{x\in \mathbb{Z}} \mathbb{E}_{h_i,h_i'\sim \mu_i} \Delta_{(h_1,h_1')}\dots \Delta_{(h_k,h_k')} f(x).
    \end{equation}

\end{defn}

When the measures $\mu_1,\dots,\mu_k$ are all uniform on $[-N,N]$, this norm is comparable to the standard Gowers norm $\|f\|_{U^k[N]}$. Crucially, if the underlying measures are sufficiently close to the uniform distribution on $[-N,N]$, information regarding the Gowers--Peluse norm can be cleanly transferred to information about the standard Gowers norm. 

\begin{lemma}[{\cite[Lemma 5.2]{Green_Sawhney}}]\label{lem: Gowers Peluse to Gowers}
    Let $\delta\in (0,1)$. Let $N,T\geq 1$ and let $f:\mathbb{Z}\rightarrow \mathbb{C}$ be a $1$-bounded function supported on $[\pm N].$ Let $\mu$ be a probability measure supported on $[\pm N]$ which satisfies that $\|\mu\|_2^2 := \sum_{[\pm N]}|x|^2 \mu(x)\leq T/N$.
    Suppose that 
    \begin{equation}
        \|f\|_{U_{\text{GP}}[N,\mu]}^2 \geq \delta.
    \end{equation}
    Then the following is true:
    \begin{equation}
        \|f\|_{U^2[N]} \gg \left(\frac{\delta^2}{T}\right)^{1/4}.
    \end{equation}
\end{lemma}
This lemma can consequently be applied iteratively for a multiset of probability measures. 

\begin{cor}[{\cite[Corollary 5.3]{Green_Sawhney}}]\label{cor: Gowers Peluse to Gowers}
Let $\delta\in (0,1)$, Let $N,T\geq 1$ and let $f:\Z\rightarrow \C$ be a 1-bounded function supported on $[\pm N]$. Let $\mu_1,\dots,\mu_k$ be probability measures supported on $[\pm N]$ such that $\|\mu_i\|_2^2\leq T/N$ for each $i.$ Suppose that 
\begin{equation}
    \|f\|_{U_{\text{GP}[N;\mu_1,\dots,\mu_k]}} \geq \delta. 
\end{equation}
Then the following is true: 
\begin{equation}
    \|f\|_{U^{2k}[N]} \gg \left(\frac{\delta^{2^k}}{T^k}\right)^{2^{-2k}}.
\end{equation}
    
\end{cor}

\section{Decomposing the count}
Recall that we are aiming to prove an asymptotic for the sum
\begin{equation}
    \sum_{\substack{x,y\in \Z \\ Q(x,y)\leq X}} \prod_{i=1}^k \Lambda(L_i(x,y)) \mathbf{1}_{\scrP}(Q(x,y)),
\end{equation}
and that we defined $M(X)$ and $S_j(X)$ for $j\in \{1,\dots,k\}$ as such:
\begin{equation}\label{eq: HOFA main term}
    M(X) := \sum_{\substack{x,y\in \Z \\ Q(x,y)\leq X}} \prod_{i=1}^k \Lambda_{\Cramer}(L_i(x,y)) \cdot \mathbf{1}_{\scrP}(Q(x,y)),
\end{equation}
\begin{equation}\label{eq: HOFA error terms}
    S_j(X):= \sum_{\substack{x,y\in \Z \\ Q(x,y)\leq X}} \prod_{i=1}^{j-1} \Lambda_{\Cramer}(L_i(x,y)) \cdot (\Lambda-\Lambda_{\Cramer})(L_j(x,y)) \cdot \prod_{i=j+1}^k \Lambda(L_i(x,y)) \cdot \mathbf{1}_{\scrP}(Q(x,y)).
\end{equation}
Then, Theorem \ref{thm: prime count} follows from the propositions below. 
\begin{prop}\label{prop: prime main term}
    Let $X\geq 1$. Then, 
    \begin{equation}
        M(X) = \mathfrak{S}(w)\cdot \frac{\rho_Q X}{\log X} + O\left(\frac{X}{(\log X)^{2}}\right),
    \end{equation}
    where the singular series is defined by 
    \begin{equation}
        \fS(w) = \prod_{p\leq w} \frac{p^{-2}\cdot \#\{(x,y)\in \F_p^2: L_1(x,y),\dots,L_k(x,y),Q(x,y)\in \F_p^\times \}}{(1-p^{-1})^{k+1}},
    \end{equation}
    and $\rho_Q$ is given by
    \begin{equation}
        \rho_Q = \frac{2\pi}{|t_{1,x}t_{1,y}| [\msO_K:\Z[\alpha]] \sqrt{\disc(K)}}.
    \end{equation}
\end{prop}

\begin{prop}\label{prop: prime error term}
    Let $X\geq 1$. For each $j\in \{1,\dots,k\}$, the following holds: 
    \begin{equation}
        S_j(X) \ll \frac{X (\log \log X)^2}{(\log X)^2}.
    \end{equation}
\end{prop}

\section{Analyzing the main term} \label{sec: primes MT}
Recall that we chose $w = \exp((\log X)^{1/10k})$ and $P(S) = \prod_{p\in S} p$. Set $t = 2^{k+1}\cdot 20\log\log X$; we do not claim that these are the optimal choices of $w$ and $t$, only that they are sufficient. Define the set of primes: $\scrP_w = \{p\in \scrP:p\leq w\}$. We now further decompose $\Lambda_{\Cramer}$:

\begin{equation}
   \Lambda_{\Cramer}(x)\coloneqq\prod_{p\le Q}\pinv{1-\frac1p}\mbf1_{\gcd(x,p)=1}=\prod_{p\le Q}\pinv{1-\frac1p}\cdot\sum_{S\subset\msP_Q}(-1)^{\abs S}\mbf1_{P(S)\mid x}=\Lambda^\#(x)+\Lambda^\flat(x),
\end{equation}
where we define
\begin{equation}
    \Lambda^\#(x) = \prod_{p\leq Q} \left(1-\frac{1}{p}\right)^{-1} \sum_{\substack{S\subset \scrP_Q \\ |S|\leq t}} (-1)^{|S|} \mathbf{1}_{P(S)\mid x},
\end{equation}
\begin{equation}
    \Lambda^{\flat}(x) = \prod_{p\leq Q} \left(1-\frac{1}{p}\right)^{-1} \sum_{\substack{S\subset \scrP_Q \\ |S|>t}} (-1)^{|S|} \mathbf{1}_{P(S)\mid x}.
\end{equation}
Consequently, we can write 
\begin{equation}
    M(X) = \sum_{\substack{x,y\in \Z\\ Q(x,y)\leq X}} \prod_{i=1}^k \Lambda_{\Cramer}(L_i(x,y)) \mathbf{1}_{\scrP}(Q(x,y)) = \sum_{\Vec{v}\in \{\#,\flat\}^k} \sum_{\substack{x,y\in \Z \\ Q(x,y)\leq X}} \prod_{i=1}^k \Lambda^{v_i}(L_i(x,y)) \mathbf{1}_{\scrP}(Q(x,y)).
\end{equation}

We claim that the contribution to $M(X)$ from any term where $v_i=\flat$ for some $i$ is negligible. 
\begin{lemma}
For $\Vec{v}\in \{\#,\flat\}^k$ such that $v_i=\flat$ for some $i$, as $X\rightarrow \infty$
\begin{equation}
    \sum_{\substack{x,y\in \Z \\ Q(x,y)\leq X}} \prod_{i=1}^k \Lambda^{v_i}(L_i(x,y)) \mathbf{1}_{\scrP}(Q(x,y))\ll X (\log X)^{-2}.
\end{equation}
\end{lemma}
\begin{proof}
    Assume that $v_j = \flat$ for some index $j$. Since $w = \exp((\log X)^{1/10k})$, we have that 
    \begin{equation}
        \prod_{p\leq w} \left(1-\frac{1}{p}\right)^{-1} \ll \log w \ll (\log X)^{1/10k}.
    \end{equation}
    Hence, we have that 
    \begin{equation}
        |\Lambda^{\flat}(x)| \ll (\log X)^{1/10k} \sum_{\substack{S\subset \scrP_w \\ |S|>t}} \mathbf{1}_{P(S)|x} \ll (\log X)^{1/10k} \tau(x), \quad \text{where }\tau(x):= \sum_{d\mid x}1.
    \end{equation}
    Since $\Lambda^\flat$ is supported on those integers with at least $t$ prime factors, we know that when $\Lambda^{\flat}(x)$ is nonzero, we must have that $2^t\leq \tau(x)$. Thus, 
    \begin{equation}
        \tau(x)\leq 2^{-t} \tau^2(x) \leq 2^{-2^{k+1}\cdot 20\log\log X} \tau^2(x).
    \end{equation}

    Next, we also note that we have the bound $|\Lambda^\#(x)|\ll (\log X)^{1/10k} \tau(x).$ So, 
    \begin{multline}
        \sum_{\substack{x,y\in \Z \\ Q(x,y)\leq X}} \prod_{i=1}^k \Lambda^{v_i}(L_i(x,y))\mathbf{1}_{\scrP}(Q(x,y)) \ll (\log X)^{1/10} \sum_{\substack{x,y\in \Z \\ Q(x,y)\leq X}}\prod_{i\neq j} \tau(L_i(x,y)) \cdot 2^{-2^{k+1}\cdot 20\log\log X} \tau^2(L_j(x,y)))\\
        \ll (\log X)^{1/10-2^{k+1}20\ln(2)} \sum_{\substack{|x|,|y|\leq c_QX^{1/2}}}\prod_{i\neq j} \tau(L_i(x,y)) \tau^2(L_j(x,y)).
    \end{multline}
    Here we have used that there is a constant $c_Q>0$ such that $\{Q(x,y)\leq X\}\subset \{|x|,|y|\leq c_Q X^{1/2}\}$. Applying \cite{Matthiesen:mult-functions} with the multiplicative functions $\tau$ and $\tau^2$ gives us the upper bound:
    \begin{equation}
        \ll (\log X)^{1/10-2^{k+1} 20\ln(2)} X (\log X)^{2^{k+1}-1} \ll X(\log X)^{-2}.
    \end{equation}

\end{proof}

It therefore suffices to find an asymptotic for the main term:
\begin{equation}
    M^\#(X) = \sum_{\substack{x,y\in \Z \\ Q(x,y)\leq X}} \prod_{i=1}^k \Lambda^{\#}(L_i(x,y)) \mathbf{1}_{\scrP}(Q(x,y)).
\end{equation}
It will be useful to convert this to a question about prime ideals; in other words, we want to replace $\mathbf{1}_{\scrP}$ with $\Lambda_K$. Define the function $\Lambda_i^\#:\calI(\msO_K)\rightarrow\R$ as:
\begin{equation}
    \Lambda_i^\#(\fa) = \begin{cases}
        \sum_{\fa = (T_x(x)+T_y(y)\sqrt{-n})} \Lambda^\#(L_i(x,y)), & \fa = (T_x(x)+T_y(y)\alpha) \text{ for some $x,y\in \Z$},\\
        0 & \text{ otherwise}.
    \end{cases}
\end{equation}
We know define the comparable sum: 
\begin{equation}
    M_K^\#(X) = \sum_{N(\fa)\leq X} \prod_{i=1}^k \Lambda_i^\#(\fa)\Lambda_K(\fa).
\end{equation}
Observe that the partial summation proof of \cite[p.41, proof of Proposition 8.5 assuming Proposition 8.6]{Green_Sawhney} directly translates over to our setting. So, it suffices to show that the following asymptotic holds for $M_K^\#(X)$. 

\begin{lemma}
    As $X\rightarrow\infty$, we have that 
    \begin{equation}
        M_K^\#(X) = \sum_{N(\fa)\leq X} \prod_{i=1}^k \Lambda_i^\#(\fa)\Lambda_K(\fa) = \fS(w) \rho_Q X + O(X(\log X)^{-1}),
    \end{equation}
    where the truncated singular series is defined by 
    \begin{equation}
        \fS(w) = \prod_{p} \frac{p^{-2}\cdot \#\{(x,y)\in \F_p^2: L_1(x,y),\dots,L_k(x,y),Q(x,y)\in \F_p^\times \}}{(1-p^{-1})^{k+1}}.
    \end{equation}
    and $\rho_{Q} = 2\pi [\msO_K:\Z[\alpha]]^{-1} \disc(K)^{-1/2} |t_{1}t_2|^{-1}.$
\end{lemma}
\begin{proof}
    Expanding the definition of $\Lambda^\#$ gives us that 
    \begin{equation}
    M_K^\#(X)=
    \sum_{N(\fa)\leq X} \prod_{i=1}^k \Lambda^\#_i(\fa) \Lambda_K(\fa) = \prod_{p\leq w}\left(1-\frac{1}{p}\right)^{-k} \sum_{i=1}^k \sum_{\substack{S_i\subset \scrP_w\\ |S_i|\leq t}} E(S_1,\dots,S_k), 
\end{equation}
where we define 
\begin{equation}
    E(S_1,\dots,S_k) = \sum_{\substack{\fa = (T_x(x)+T_y(y)\alpha)\\ N(\fa)\leq X}} \Lambda_K(\fa) \prod_{i=1}^k \mathbf{1}_{P(S_i)\mid L_i(x,y)}.
\end{equation}

For now, fix the sets $S_1,\dots,S_k$. Recall that 
\begin{equation}
    T_x(x) = t_{1,x}x+ t_{0,x}, \quad T_y(y) = t_{1,y} y + t_{0,y}.
\end{equation}
Define the function: 
\begin{equation}
    T:\Z[\alpha] \mapsto \Z[\alpha]: \quad T(u+v\alpha) = T_x(u) + T_y(v)\alpha.
\end{equation}
We can now define the following lattice in $\msO_K:$
\begin{equation}
    \Gamma(S_1,\dots,S_k) := \{u+v\alpha\in T(\Z[\alpha]): \forall i = 1,\dots,k, P(S_i) \mid L_i(t_{1,x}^{-1}(u-t_{0,x}), t_{1,y}^{-1}(v-t_{0,y}))\}.
\end{equation}
Observe that these functions above are still linear in $u$ and $v$ (and that we treat all of the coefficients of $L_1,\dots,L_k,T_x,T_y$ as fixed constants). Then the lattice $\Gamma(S_1,\dots,S_k)$ will be a finite union of congruence classes mod $\fm_{\Vec{S}}$, where $\fm_{\Vec{S}}$ is the principal ideal:
\begin{equation}
    \fm_{\Vec{S}} = \left(\disc(K)\prod_{i=1}^k P(S_i) |t_{1,x}t_{1,y}|\right)\subset \msO_K.
\end{equation}
Therefore, we can express the indicator function of $\Gamma(S_1,\dots,S_k)$ in terms of Dirichlet characters $\chi:\msO_K/\fm_{\Vec{S}} \rightarrow \C^\times$:
\begin{equation}
    \mathbf{1}_{\Gamma(S_1,\dots,S_k)}(u+v\alpha) = \sum_{\chi\bmod \fm_{\Vec{S}}} c_\chi^{\Vec{S}} \chi(u+v\alpha).
\end{equation}
We can rewrite $E(S_1,\dots,S_k)$ as:
\begin{align}
    E(S_1,\dots,S_k) &= \sum_{u\in \msO_K^*} \sum_{\substack{\gamma\in \msO_K/\msO_K^* \\ N(\gamma)\leq X}} \sum_{\chi \bmod \mathfrak{m}} c_\chi^{\Vec{S}} \chi(u\gamma) \Lambda_K((u\gamma))\\
    &= |\msO_K^*| \sum_{\chi\bmod \mathfrak{m}} c_\chi^{\Vec{S}} \sum_{N(\fa)\leq X} \chi(\fa) \Lambda_K(\fa) \cdot \frac{1}{|C|} \sum_{\psi\in \hat{C}} \psi(\fa),
\end{align}
where $C$ denotes the class group of $K.$ Note that by \cite[(8.23)]{Green_Sawhney}, $|c_\chi^{\Vec{S}}|\leq |c_{\chi_0}^{\Vec{S}}|\leq 1$ for all characters $\chi\bmod \fm_{\Vec{S}}$. 

Unsurprisingly, our main term comes from the principal character, i.e. when $\psi=1$ and $\chi = \chi_0$. By Theorem \ref{thm: prime ideal theorem with characters}, we have that 
\begin{multline}
    M^\#_K(X) = X \prod_{p\leq w} \left(1-\frac{1}{p}\right)^{-k} \sum_{\substack{S_1,\dots,S_k\subset \scrP_w\\ |S_i|\leq t}} (-1)^{S_1+\dots +S_k}  \frac{|\msO_K^*|}{|C|} \frac{\#\{\gamma\in \Gamma(S_1,\dots,S_k)/\mathfrak{m}_{\Vec{S}}: \gcd(\gamma,\mathfrak{m}_{\Vec{S}})=1\}}{|(\msO_K/\mathfrak{m}_{\Vec{S}})^*|} \\
    + O(E_{Siegel} + E_{error}). 
\end{multline}
Here, $E_{Siegel}$ tracks the contribution from the potential Siegel zeros (note that our moduli range up to $w^k = \exp((\log X)^{1/10})$) and $E_{error}$ is given by all of the other error terms in Theorem \ref{thm: prime ideal theorem with characters}.

Let $\chi_{D_*}$ denote the potential exceptional character with Siegel zero $\beta\in (1/2,1)$ from Theorem \ref{thm: prime ideal theorem with characters}. Then the contribution from the Siegel zeros is bounded by:
\begin{align}
    E_{Siegel} &= (\log w)^k X^{\beta} \sum_{\substack{S_1,\dots,S_k\subset \scrP_w\\ |S_i|\leq t}} \sum_{\substack{\psi^* \bmod \mathfrak{m}^{S_1,\dots,S_k}\\ \psi^* \sim \chi_{D_*}}} |c_{\psi^*}^{\Vec{S}}|\\
    &\ll (\log w)^{k+1} X^{\beta} \sum_{\substack{S_1,\dots,S_k\subset \scrP_w\\ D_*\mid \prod P(S_i)}} \frac{1}{P(S_1)\dots P(S_k)}
\end{align}
Here we have used the upper bound 
\begin{equation}
    |c_\chi^{\Vec{S}}| \ll \frac{\log w}{P(S_1)\dots P(S_k)},
\end{equation}
which follows from \cite[(8.23),(8.24),(8.30)]{Green_Sawhney}. Define $\tau_k(D_*)$ to be the number of tuples $(d_1,\dots,d_k)$ such that $d_1\dots d_k\mid D_*$. 
Then we have that 
\begin{equation}
    E_{Siegel}\ll (\log w)^{k+1} X^{\beta} \tau_k(D_*) \sup_{d_1\dots d_k = D_*} \sum_{\substack{S_1,\dots,S_k\subset \scrP_w \\ d_i\mid P(S_i)}} \frac{1}{P(S_1)\dots P(S_k)}.
\end{equation}
Next, recall that for any $\epsilon>0$, we have that $\beta\ll 1-c_\epsilon D_*^\epsilon$ and that for any $\eta>0$, $\tau_k(D_*)\ll D_*^\eta$. Further, since $\log w = (\log X)^{1/10k}$, we have that 
\begin{equation}
    E_{Siegel} \ll (\log w)^{k+1} X^{1-c_\epsilon D_*^\epsilon} D_*^\eta (\log w)^k D_*^{-1} \ll X(\log X)^{1/4} X^{-c_\epsilon D_*^{\eta}}D_*^{\eta-1}.
\end{equation}
Taking the maximum over the possible values of $D_*$ (which can range up to $e^{\sqrt{\log X}}$), we have that 
\begin{equation}
    E_{Siegel} \ll X (\log X)^{-1}.
\end{equation}

Next, we consdier the contribution from $E_{error}$. 
\begin{equation}
    E_{error} = X(\log w)^k \binom{|\scrP_w|}{t}^{k}e^{-c\sqrt{\log X}}.
\end{equation}
We estimate this binomial as:
\begin{equation}
    \binom{|\scrP_w|}{t} \ll \exp(t(\log X)^{\frac{1}{10k}}) \ll \exp(2^{k+1}20 \log\log X (\log X)^{\frac{1}{10k}})\ll \exp((\log X)^{\frac{1}{9k}}).
\end{equation}
Hence, we know that for any $A>0$,
\begin{equation}
    E_{error} \ll X \exp(-c\sqrt{\log X} + (\log X)^{\frac{1}{9k}} + \frac{1}{10}\log\log X) \ll_A \frac{X}{(\log X)^A}.
\end{equation}

Finally, we turn our attention to the main term. Let us first set up some notation: we say that $\gcd(\gamma,\fm) =1$ if and only if for all prime ideals $\fp\mid \fm$, we have that $\gamma\not\in \fp$. 
Consider further the expression:
\begin{equation}
    X \prod_{p\leq w} \left(1-\frac{1}{p}\right)^{-k} \sum_{\substack{S_1,\dots,S_k\subset \scrP_w\\ |S_i|\leq t}} (-1)^{|S_1|+\dots +|S_k|} \cdot \frac{|\msO_K^*|}{|C|} \frac{\#\{\gamma\in \Gamma(S_1,\dots,S_k)/\mathfrak{m}_{\Vec{S}}: \gcd(\gamma,\mathfrak{m}_{\Vec{S}})=1\}}{|(\msO_K/\mathfrak{m}_{\Vec{S}})^*|} 
\end{equation}
Similarly to in \cite[\S8.8]{Green_Sawhney}, we remove the constraint that $|S_1|,\dots,|S_k|\leq t$ with a Brun-sieve type computation; we refer to \cite[\S8.8]{Green_Sawhney} for a more detailed argument. Without loss of generality, let us bound the contribution when $|S_1|>t$. This accrues an error of the form:
\begin{align}
    X(\log w)^{2k-1} \sum_{\substack{S\subset \scrP_w\\ |S|>t}}\prod_{p\in S} \frac{1}{p-1} \ll X (\log w)^{2k-1} (\log X)^{-2^{k+1} 20 \ln(2)} \ll X(\log X)^{-1}.
\end{align}
This error is sufficient for the remaining analysis in \cite[\S8.8]{Green_Sawhney}, so we can instead analyze the expression:
\begin{equation}
    X\cdot \frac{|\msO_K^*|}{|C|}\prod_{p\leq w} \left(1-\frac{1}{p}\right)^{-k} \sum_{\substack{S_1,\dots,S_k\subset \scrP_w}} (-1)^{\sum |S_i|}  \times \frac{\#\{\gamma\in \Gamma(S_1,\dots,S_k)/\mathfrak{m}_{\Vec{S}}: \gcd(\gamma,\mathfrak{m}_{\Vec{S}})=1\}}{|(\msO_K/\mathfrak{m}_{\Vec{S}})^*|}.
\end{equation}

Let $\mathfrak{m} = (\disc(K)P(\scrP_w)|t_{1,x}t_{1,y}|)$. By the Chinese remainder theorem and M\"obius inversion, we can rewrite this sum as:
\begin{equation}
    \frac{\#\left\{u+v\alpha \in T(\Z[\alpha])/\mathfrak{m}: \begin{array}{c}
         \gcd(L_i(t_{1,x}^{-1}(u-t_{0,x}), t_{1,y}^{-1}(v-t_{0,y})),P(\scrP_w)|t_{1,x}t_{1,y}|)=1,\\
         u+v\alpha \in (\msO_K/\mathfrak{m})^*
    \end{array}\right\}}{|(\msO_K/\mathfrak{m})^*|}.
\end{equation}
Rearranging this numerator (noting that $\fm$ is a principal ideal in $\msO_K$), recalling the formula for $|(\msO_K/\fm)^*|$ in \cite[(8.27)]{Green_Sawhney}, and applying the class number formula for $|\msO_K^*|/|C|$ tells us that the main term is given by
\begin{equation}
    \frac{2\pi X}{|t_{1,x}t_{1,y}| [\msO_K:\Z[\alpha]]\sqrt{\disc(K)}}\prod_{p\leq w} \frac{p^{-2}\#\{x,y\bmod p: p\nmid L_i(x,y), p\nmid Q(x,y)\}}{(1-p^{-1})^{k+1}} \prod_{p>w} \left(1-\frac{\chi(p)}{p}\right),
\end{equation}
where $\chi$ is the real quadratic Dirichlet character such that $\zeta_K(s) = \zeta(s)L(s,\chi)$. Since $\sum_{p>w} \chi(p)/p = O(\exp(-c\sqrt{\log w}))$, we get the desired asymptotic below:


\begin{multline}
    M_K^\#(X) = \sum_{N(\fa)\leq X} \prod_{i=1}^k \Lambda^\#_i(\fa)\Lambda_K(\fa) \\  = \frac{2\pi X}{|t_{1,x}t_{1,y}|[\msO_K:\Z[\alpha]]\sqrt{\disc(K)}} \prod_{p\leq w} \frac{p^{-2}\#\{x,y\bmod p: p\nmid L_i(x,y), p\nmid Q(x,y)\}}{(1-p^{-1})^{k+1}} \\ + O(X/\log X). 
\end{multline}

\end{proof}

\section{Sieves and semiprimes}\label{sec: primes sieve}

We now embark on the proof of Proposition \ref{prop: prime error term}, that is, it is time to bound the error term $S_j(X)$. The first step of our proof is a sieve based on the sieve of Duke, Friedlander, and Iwaniec \cite{DFI}; we will use the version formulated by Green and Sawhney \cite{Green_Sawhney} that is adapted to the number field setting. Before doing so, let us first define the notion of Type I and Type II estimates. 

\begin{defn}
Let $w:\calI(\msO_K)\rightarrow\C$. Then we say that $w$ satisfies a \define{Type I bound (with savings of $(\log X)^{-B}$) at level $D$} if for any $D_0\leq D$, 
    \begin{equation}\label{eq: type I def}
        \sum_{\substack{\frak{d}\subset \mathcal{O}_K \\ N(\frak{d})\sim D_0}} \left|\sum_{\substack{N(\fa)\leq X \\ \frak{d} \mid \frak{a}}} w(\frak{a}) \right| \ll_B X(\log X)^{-B} \|w\|_\infty.
    \end{equation}
     
\end{defn}
This is also known as having level of distribution $D$. Bounds of these forms measure how well we understand the function $w:\calI(\msO_K)\rightarrow \C$ along arithmetic progressions of growing moduli. 

\begin{defn}
    We say that $w$ satisfies a \define{Type II bound (with savings $(\log X)^{-B}$) between $[D_1,D_2]$} if for any $D\in [D_1,D_2]$ and any $1$-bounded sequences $\alpha,\beta:\calI(\msO_K)\rightarrow \C$, 
    \begin{equation}\label{eq: type II def}
        \sum_{\substack{\frak{a},\frak{b}\subset \mathcal{O}_K \\ N(\fa\fb)\leq X\\  N(\frak{a})\sim D}} \alpha_{\frak{a}}\beta_{\frak{b}} w(\frak{a}\frak{b}) \ll_B X(\log X)^{-B} \|w\|_\infty.
    \end{equation}
\end{defn}
Type II information is also considered bilinear information (in comparison to Type I information, which is linear information). With these two definitions in hand, we can now record the sieve statement. 

\begin{thm}[{\cite[Lemma 3.2]{Green_Sawhney}}]\label{thm: DFI sieve}
    Let $A,C>1$ be parameters; let $B=2A+4$ and suppose $C\geq B$. Suppose that $w:\calI(\msO_K)\rightarrow \C$ satisfies Type I bounds (with savings $(\log X)^{-B}$) at $X^{1/2}(\log X)^{-C}$ and Type II bounds (with savings $(\log X)^{-B}$) between $[(\log X)^C, X^{3/8}]$. Then, 
    \begin{equation}
        \sum_{N(\fp)\leq X} w(\fp) + \sum_{\substack{\fp,\fq\subset\msO_K\\ N(\fp\fq)\leq X \\ N(\fp),N(\fq) \geq X^{1/2} \exp(-(\log \log X)^2)}} w(\fp\fq) \ll_{A,C} X(\log X)^{-A} \|w\|_\infty.
    \end{equation}
\end{thm}

To see how we will apply the sieve above to $S_j(X)$, we need to further define some functions. Fix an index $j\in \{1,\dots,k\}$ and define $w_j:\calI(\msO_K)\rightarrow\C$: 
\begin{equation}
    w_j(\fa) := \prod_{i=1}^k w_{ij}(\fa),\quad \text{where }w_{ij}(\fa) = (\log X)^{-1}\sum_{\substack{x,y\in \Z \\ \fa=(T_x(x) + T_y(y)\alpha)}} \begin{cases}
        \Lambda_{\Cramer}(L_i(x,y)), & i<j,\\
        \Lambda-\Lambda_{\Cramer}(L_j(x,y)), & i=j, \\
        \Lambda(L_i(x,y)), & i>j.
    \end{cases}
\end{equation}
Since $Q(x,y) = N_{K/\Q}(T_x(x)+T_y(Y)\alpha)$, we know that if $Q(x,y)$ is prime then $(T_x(x)+T_y(y)\alpha)$ is a prime ideal. Hence, 
\begin{equation}
    S_j(X) = (\log X)^k\sum_{\substack{N(\fp)\leq X}} w_{j}(\fp) + O_\epsilon(X^{1/2+\epsilon})
\end{equation}
Here, the error term comoes from a classical bound on the inert prime ideals.

Our strategy to bound $S_j(X)$ is now to show that $w_j:\calI(\msO_K)\rightarrow\C$ satisfies Type I bounds (with savings $(\log X)^{-B}$) at level $X^{1/2}(\log X)^{-C}$ and Type II bounds (with savings $(\log X)^{-B}$) between $[(\log X)^C,X^{3/8}]$ for $B$ arbitrarily large. Let us record our Type I and Type II estimates.

The following are adaptations of \cite[Proposition 4.2]{Green_Sawhney} and \cite[Proposition 4.3]{Green_Sawhney}; however, we must recreate the proofs in the setting where we have $k$ linear forms.
\begin{prop}[Type I]\label{prop: type I}
    Fix $j\in \{1,\dots,k\}$. Let $\delta\in (0,1)$ be a parameter. There is an explicit constant $c_k>0$ such that if for $D\leq \delta^{O(1)} X^{1/2},$ we have that 
    \begin{equation}
        \sum_{\fd: N(\fd)\sim D} \left|\sum_{\substack{N(\fa)\leq X \\ \fd\mid \fa}} w_j(\fa)\right|\geq \delta X,
    \end{equation}
    then for any $i\in \{1,\dots,k\}$, 
    \begin{equation}
        \|\lambda_{ij}\|_{U^{c_k}[X^{1/2}]} \gg \left(\frac{\delta}{\log X}\right)^{O(1)},
    \end{equation}
    where $\lambda_{ij}\in \{\Lambda_{\Cramer},\Lambda-\Lambda_{\Cramer}, \Lambda\}$ is the function satisfying 
    \begin{equation}
        w_{ij}(\fa) = (\log X)^{-1} \sum_{\substack{x,y\in \Z \\ \fa= (T_x(x)+T_y(y)\alpha )}} \lambda_{ij}(L_i(x,y)).
    \end{equation}
\end{prop}
\begin{prop}[Type II]\label{prop: type II}
    Fix $j\in \{1,\dots,k\}$. Let $\delta\in (0,1)$ be a parameter. There is an explicit constant $r_k>0$ such that if for $D\in [\delta^{O(1)},X^{3/8}]$ and some 1-bounded complex sequences $\alpha,\beta:\calI(\msO_K)\rightarrow\C$ we have that 
    \begin{equation}
        \left|\sum_{\substack{N(\fa\fb)\leq X \\ N(\fa)\sim D}} \alpha_\fa\beta_{\fb} w_j(\fa\fb)\right| \geq \delta X,
    \end{equation}
    then for any $i\in \{1,\dots,k\}$, 
    \begin{equation}
        \|\lambda_{ij}\|_{U^{r_k}[X^{1/2}]} \gg \left(\frac{\delta}{\log X}\right)^{O(1)},
    \end{equation}
    where $\lambda_{ij}\in \{\Lambda_{\Cramer},\Lambda-\Lambda_{\Cramer}, \Lambda\}$ is the function satisfying 
    \begin{equation}
        w_{ij}(\fa) = (\log X)^{-1} \sum_{\substack{x,y\in \Z \\ \fa= (T_x(x)+T_y(y)\alpha)}} \lambda_{ij}(L_i(x,y)).
    \end{equation}
\end{prop}

While neither of the above propositions are exactly of the shape \eqref{eq: type I def} and \eqref{eq: type II def}, they immediately give the requisite Type I and Type II bounds upon noting that $\|w_j\|_\infty\leq 1$ and applying the recent breakthroughs of Leng \cite{Leng} and Leng, Sah, and Sawhney \cite{LengSahSawhney}. 
Specifically, Leng \cite{Leng}, building on the work of Leng, Sah, and Sawhney \cite{LengSahSawhney} on quantitative improvements for the inverse theorem for $U^s$, established that for any $s\geq 1$ and any $B>0$,   
\begin{equation}
    \|\Lambda-\Lambda_{\Cramer}\|_{U^s[X^{1/2}]}\ll_C (\log X)^{-B}.
\end{equation}
This immediately allows us to take $\delta = (\log X)^{-C}$ for $C$ sufficiently large in Proposition \ref{prop: type I} and Proposition \ref{prop: type II} above. Thus, the following corollary holds. 
\begin{cor}\label{cor: type I and type II}
    Fix $j\in \{1,\dots,k\}$. Let $B\geq 1$. Then there exists $C\geq 1$ such that for $D\leq X^{1/2}(\log X)^{-C}$, 
    \begin{equation}
        \sum_{N(\fd)\sim D} \left|\sum_{\fd\mid \fa}w_j(\fa)\right|\ll_B X(\log X)^{-B}. 
    \end{equation}
    There exists $C'>0$ such that for $D\in [(\log X)^C,X^{3/8}]$, for any 1-bounded sequences $\alpha,\beta:\calI(\msO_K)\rightarrow \C$, 
    \begin{equation}
        \sum_{\substack{N(\fa\fb)\leq X \\ N(\fa)\sim D}} \alpha_\fa \beta_\fb w_j(\fa\fb) \ll_B X(\log X)^{-B}. 
    \end{equation}
\end{cor}

Corollary \ref{cor: type I and type II} now, in addition to Theorem \ref{thm: DFI sieve}, implies that for any $A>0$,
\begin{equation}
    \sum_{N(\fp)\leq X} w_j(\fp) + \sum_{\substack{N(\fp\fq)\leq X \\ N(\fp),N(\fq)\geq X^{1/2} \exp(-(\log\log X)^2)}} w_j(\fp\fq) \ll X (\log X)^{-A}.
\end{equation}
To upper bound $S_j(X)$, it remains to establish bounds on the sums over semiprimes. 

\begin{lemma}\label{lem: bound on semiprimes}
    Fix $j\in \{1,\dots,k\}$. Then as $X\rightarrow\infty$, 
    \begin{equation}
       \sum_{\substack{N(\fp\fq)\leq X \\ N(\fp),N(\fq)\geq X^{1/2}\exp(-(\log \log X)^2)}} w_j(\fp\fq) \ll \frac{X(\log\log X)^2}{(\log X)^{k+2}}.
    \end{equation}
\end{lemma}
\begin{rem}
    Recall that 
    \begin{equation}
        S_j(X) = \sum_{N(\fp)\leq X} (\log X)^k w_j(\fp).
    \end{equation}
    So, Lemma \ref{lem: bound on semiprimes} will eventually give us an error term of the form $O(X(\log\log X)^2/(\log X)^2).$
\end{rem}
\begin{proof}[Proof of Lemma \ref{lem: bound on semiprimes}]

    Let $\theta:\Z \rightarrow \R$ denote the function: 
    \begin{equation}
        \theta(n) = \begin{cases}
            \log |n|, & n\in \scrP,\\
            0, & \text{otherwise.}
        \end{cases}
    \end{equation}
    Observe that $\theta$ and $\Lambda$ only differ on positive and negative prime powers (which is a subset of $\{\pm n^d:n\in \Z, d>1\}$, which itself is a very sparse set). Hence, we can write
    \begin{equation}
        (\log X)^k\sum_{\substack{N(\fp\fq)\leq X \\ N(\fp),N(\fq)\geq X^{1/2}\exp(-(\log \log X)^2)}} w_j(\fp\fq) = \sum_{\substack{N(\fp\fq)\leq X \\ N(\fp),N(\fq)\geq X^{1/2}\exp(-(\log \log X)^2)}} \prod_{i=1}^k g_{ij}(\fp\fq) + O_\epsilon(X^{1/2+\epsilon}),
    \end{equation}
    where we define now 
    \begin{equation}
        g_{ij}(\fa) = \sum_{\substack{x,y\in \Z \\ \fa=(T_x(x)+T_y(y)\alpha)}} \begin{cases}
            \Lambda_{\Cramer}(L_i(x,y)), & i<j, \\ 
            (\theta-\Lambda_{\Cramer})(L_j(x,y)), & i=j, \\
            \theta(L_i(x,y)), & i>j.
        \end{cases}
    \end{equation}

    Next, fix a prime $\fp$ with $N(\fp)\in [X^{1/2}\exp(-(\log\log X)^2), X^{1/2}].$ We can assume that $\fp$ is not a rational prime, as the contribution from when either $\fp$ or $\fq$ is a rational prime is bounded by $O_\epsilon(X^{1/2+\epsilon})$. If we knew that for each such non-rational prime $\fp$ the following holds
    \begin{equation}\label{eq:upp bound sieve goal for each p}
        \sum_{N(\fq)\leq X/N(\fp)} \prod_{i=1}^k |g_{ij}(\fp\fq)| \ll \frac{X}{N(\fp)\log X},
    \end{equation}
    then we can immediately see that 
    \begin{align}
        \sum_{\substack{N(\fp\fq)\leq X \\ N(\fp),N(\fq)\geq X^{1/2}\exp(-(\log\log X)^2)}} \prod_{i=1}^{k} g_{ij}(\fp \fq) &\leq \sum_{N(\fp)\in [X^{1/2}\exp(-(\log\log X)^2),X^{1/2}]} \sum_{N(\fq)\leq X/N(\fp)} \prod_{i=1}^k|g_{ij}(\fp\fq)| \\
        &\ll \sum_{N(\fp)\in [X^{1/2}\exp(-(\log\log X)^2),X^{1/2}]} \frac{X}{N(\fp)\log X} \\
        &\ll \frac{X(\log\log X)^2}{(\log X)^2}. 
    \end{align}
    The bound \eqref{eq:upp bound sieve goal for each p} then suffices to establish the lemma.

    Our new goal is now to establish \eqref{eq:upp bound sieve goal for each p}. Since $g_{ij}:\calI(\msO_K)\rightarrow \C$ is supported on principal ideals (for each $i$), we only need to consider when $\fp\fq$ is a principal ideal. We can write that
    \begin{equation}
        \fp = (a+b\alpha)\mfc_\fp, \quad \fq = (c+d\alpha)\mfc_\fq,
    \end{equation}
    where $\mfc_\fp$ (resp. $\mfc_\fq$) are representatives of the ideal class of $\fp$ (resp. $\fq$); these will satisfy that $N(\mfc_\fp) = O_n(1)$ (resp. $N(\mfc_\fq) = O_n(1)$). Then we only must consider when $\mfc_\fp\mfc_\fq$ is principal; write that $\mfc_\fp\mfc_\fq = (\gamma)$ and $\mfc_\fq = \mfc_{\fp^{-1}}$. Note that we can scale $\mfc_\fp$ and $\mfc_{\fp^{-1}}$ such that $N(\mfc_p), N(\mfc_{p^{-1}}),\gamma\in \Z$. For more details on ideal classes, we refer the reader to \cite[Lemma 2.1]{Green_Sawhney}. 

    In our new notation, $g_{ij}:\calI(\msO_K)\rightarrow \C$ is supported on ideals where there is some unit $u\in \msO_K^*$ such that 
    \begin{equation}
        u\gamma (a+b\alpha) (c+d\alpha) = T_x(x)+T_y(y)\alpha = t_{1,x} x + t_{0,x} +(t_{1,y}y+t_{0,y})\alpha.    \end{equation}
    Since $t_{1,x},t_{1,y}\neq 0$, there exist functions
    \begin{equation}
        \psi_{1,u}(a,b),\psi_{2,u}(a,b), \phi_{1,u}(a,b), \phi_{2,u}(a,b) \in \Q[a,b]
    \end{equation}
    that are linear in $a$ and $b$ with coefficients depending on $T_x$, $T_y$ and $K$, and such that 
    \begin{equation}
        x+y\alpha = (\psi_{1,u}(a,b)c + \psi_{2,u}(a,b)d) + (\phi_{1,u}(a,b)c + \phi_{2,u}(a,b)d)\alpha.
    \end{equation}
    Thus, we have that 
    \begin{equation}
        \prod_{i=1}^k |g_{ij}(\fp\fq)| \ll \sum_{u\in \msO_K^*} \prod_{i=1}^k (\theta+\Lambda_{\Cramer})(L_i(\psi_{1,u}(a,b)c+\psi_{2,u}(a,b)d, \phi_{1,u}(a,b)c+\phi_{2,u}(a,b) d)), 
    \end{equation}
    where we use the convention that $\theta$ and $\Lambda_\Cramer$ are zero on non-integers. Since $N(\fq)\leq X/N(\fp)$, we can assume that $|c|,|d|\ll (X/N(\fp))^{1/2}$. To establish \eqref{eq:upp bound sieve goal for each p}, it suffices to prove the upper bound for each unit $u\in \msO_K^*$. 

    From now on, we fix a unit $u\in \msO_K^*$ and we drop the dependency on $u$ in our notation. Recall that $w = \exp((\log X)^{1/10k})$ and set $w'=X^{1/100}$. Observe that for $n\gg w'$, we have that 
    \begin{equation}
        \theta(n)\ll \Lambda_{\text{Cram\'er},w'}(n), \quad \mathbf{1}_{\scrP}(n) \ll \frac{\Lambda_{\text{Cram\'er},w'}(n)}{\log X }.
    \end{equation}
    Since the pairs $(c,d)$ such that $|N(\mfc_{\fp^{-1}}) N_{K/\Q}(c+d\alpha)|\leq w'$ will contribute negligibly to the final sum (as they contribute a factor of $O(X^{1/100}(\log X)^k)$). It suffices to consider the sum: 
    \begin{multline}
        \sum_{|c|,|d|\ll (X/N(\fp))^{1/2}} \prod_{i=1}^k (\Lambda_{\Cramer,w}+\Lambda_{\Cramer,w'})(L_i(\psi_1(a,b)c+\psi_2(a,b)d, \phi_1(a,b)c+\phi_2(a,b)d)) \\ \times \frac{\Lambda_{\Cramer,w'}(N_{K/\Q}(c+d\alpha)N(\mfc_{\fp^{-1}}))}{\log X}
    \end{multline}
    Hence, let $w_0 =w'$ and $w_i\in \{w,w'\}$ for $i=1,\dots,k$. It is therefore enough to show that 
    \begin{equation}\label{eq:upp bound sieve goal removed log}
        \sum_{|c|,|d|\ll Y} \Lambda_{\text{Cram\'er},w_0}((c^2+nd^2)N(\mfc_{\fp^{-1}}))\prod_{i=1}^k \Lambda_{\text{Cram\'er},w_i}(L_i(\psi_1(a,b)c+\psi_2(a,b)d,\phi_1(a,b)c+\phi_2(a,b)d)) \ll Y^2.
    \end{equation}
    (In the application above, take $Y = (X/N(\fp))^{1/2}$). 

    Suppose that $w_k\leq w_{k-1}\leq \dots \leq w_1\leq w_0=w'$ (all other permutations of $w_1,\dots,w_k$ can be handled analogously). We will use the large sieve to bound \eqref{eq:upp bound sieve goal removed log}. Let $p$ be a prime satisfying $w_{j+1}< p\leq w_j$. Define: 
    \begin{equation}
        \Omega_p := \{(u,v)\in \F_p^2: N(\mfc_{\fp^{-1}})N_{K/\Q}(u+v\alpha) \prod_{i=1}^j L_i(\psi_1(a,b)u+\psi_2(a,b)v,\phi_1(a,b)u+\phi_2(a,b)v) = 0\}.
    \end{equation}
    For all other primes $p$, define $\Omega_p = \emptyset$. 
    To apply \cite[Lemma C.3]{Green_Sawhney}, we observe that
    \begin{multline}
        \Lambda_{\text{Cram\'er},w_0}(N_{K/\Q}(u+v\alpha)N(\mfc_{\fp^{-1}}))\prod_{i=1}^k \Lambda_{\text{Cram\'er},w_i}(L_i(\psi_1(a,b)c+\psi_2(a,b)d,\phi_1(a,b)c+\phi_2(a,b)d)) \\ \ll \prod_{i=0}^k \log w_i \cdot \prod_{p} \mathbf{1}_{(c,d)\in \Omega_p^c}.
    \end{multline}
    Following \cite[p.16]{Green_Sawhney} line-by-line, we see that if $w_{j+1}\leq p\leq w_j$ and $p\nmid \disc(K)\gamma$ and if $L_i(\psi_1(a,b)u+\psi_2(a,b)v,\phi_1(a,b)u+\phi_2(a,b)v)\in \F_p[c,d]$ is a nonconstant linear form, then
    \begin{equation}
        \beta_p = \frac{|\omega_p|}{p^2} = \frac{1}{p}\left(j+1+\left(\frac{\disc(K)}{p}\right)\right)+ O\left(\frac{1}{p^2}\right).
    \end{equation}
    This comes from the fact that there will be $j$ nontrivial forms in the defining equation for $\Omega_p$ and a quadratic form of the shape $N_{K/\Q}(u+v\alpha)$. Consequently, following again \cite[p.16-17]{Green_Sawhney}, we have that 
    \begin{equation}
         \prod_{p\leq w_0} \left(1-\beta_p\right)\asymp \left(\prod_{i=0}^k\log(w_i)\right)^{-1}.
    \end{equation}
    Plugging this into \eqref{eq:upp bound sieve goal removed log}, we see that we have a bound of 
    \begin{equation}
        \sum_{\substack{|c|,|d|\ll Y}} \prod_{i=0}^k \log(w_i) \prod_{p} \mathbf{1}_{(c,d)\in \Omega_p}  \ll \left(\prod_{i=0}^k \log w_i\right) Y^2 \prod_{p\leq w_0} (1-\beta_p) \ll Y^2
    \end{equation}
    Thus, this establishes \cref{eq:upp bound sieve goal removed log} as desired and completes the proof of the lemma. 
\end{proof}

\section{Type I estimates}\label{sec: primes type I}
In this section, we prove Proposition \ref{prop: type I}. For this proof, we will follow the proof of \cite[Proposition 4.2]{Green_Sawhney} until the last step. However, before proceeding with the proof, let us establish some notation and lemmas that will be useful for both this section and \S\ref{sec: primes type II}.

Recall that 
\begin{equation}
    L_i(x,y) = L_i^0(x,y) + l_i^0 = \beta_{i,x}^0 x + \beta_{i,y}^0 y + l_i^0.
\end{equation}
We define the modified linear form: 
\begin{equation}
    L_{i,T}(x,y) = L_i\left(\frac{x-t_{0,x}}{t_{1,x}}, \frac{y-t_{0,y}}{t_{1,y}}\right) = \beta_{i,x} x + \beta_{i,y} y + l_i \in \Q[x,y],
\end{equation}
where $\beta_{i,x},\beta_{i,y},l_i\in \frac{1}{M}\Z$ where $M = O_{K,L_1,\dots,L_k,T_x,T_y}(1).$

For each pair $i\neq j$, define: 
\begin{equation}
    \delta_{ij} = \beta_{i,x}\beta_{j,y} - \beta_{i,y}\beta_{j,x} = \det\begin{pmatrix}
        \beta_{i,x} & \beta_{i,y} \\
        \beta_{j,x} & \beta_{j,y}
    \end{pmatrix} \in \frac{1}{M^2}\Z.
\end{equation}
Since we assumed that the homogeneous parts of $L_i$ and $L_j$ are pairwise linearly independent, we have that $\delta_{ij}$ is always nonzero. This will be important for our proof of both Proposition \ref{prop: type I} and Proposition \ref{prop: type II}.

Next, we recall that $K=\Q(\alpha)$ is an imaginary quadratic field. For the sake of sanity, let us assume that $\alpha = \sqrt{-n}$ for $n\in \mathbb{N}$ squarefree and hence $N_{K/\Q}(x+y\alpha) = x^2 + ny^2.$ The following computation can be recreated for any choice of integral $\alpha$ by changing the linear transformations $T_x$ and $T_y$ (ex: if $\alpha = 2+i$, this $N_{K/\Q}(x+y\alpha) = N_{K/\Q}((x+2)+yi)$; hence, taking $T_x = x+2$ and $T_y = y$ would have the same effect). 

Let $i,j\in \{1,\dots,k\}$ (not necessarily distinct). Define: 
\begin{equation}
    \delta_{ij}' = n\beta_{i,x}\beta_{j,x} + \beta_{i,y}\beta_{j,y}.
\end{equation}
Observe that if $i=j$, $\delta_{ii}' = \beta_{j,y}^2 + n\beta_{i,x}^2$. This constant will only be used in the proof of Proposition \ref{prop: type II}, but we record it now anyways. 

\begin{lemma}\label{lem: differencing and algebraic manipulations}
    Assume that for $a,b\in \Z$, and for $j\in \{1,\dots,k\}$, 
    $$\gcd(M (\beta_{j,x} a+\beta_{j,y} b), M(\beta_{j,y} a - n\beta_{j,x} b))= e_j.$$
Then if for $x,y,x',y'\in \Z$,  $$L_{j,T}(ax-nby,ay+bx) = L_{j,T}(ax'-nby',ay'+bx'),$$ there exists some $h_{j}\in \Z$ such that for any $i\neq j$, 
\begin{equation}
    L_{i,T}(ax-nby,ay+bx)-L_{i,T}(ax'-nby',ay'+bx') = \frac{M h_{j}\delta_{ij}}{e_j} \cdot (a^2+nb^2). 
\end{equation}
\end{lemma}
\begin{proof}
    We start with the relation that $L_{j,T}(ax-nby,ay+bx) = L_{j,T}(ax'-nby',ay'+bx').$ Expanding both sides gives us that 
    \begin{equation}
        (y-y') (\beta_{j,y} a-n\beta_{j,x} b) = (x'-x) (\beta_{j,x} a+ \beta_{j,y} b).
    \end{equation}
    Since $\gcd(M (\beta_{j,x} a+\beta_{j,y} b), M(\beta_{j,y} a - n\beta_{j,x} b))= e_j$, this gives us that for some $h_{j}\in \Z$, 
    \begin{equation}
        y-y' = \frac{M(\beta_{j,x} a+\beta_{j,y} b)}{e_j} \cdot h_{j}, \quad x'-x = \frac{M(\beta_{j,y} a -n\beta_{j,x} b)}{e_j}\cdot h_{j}.
    \end{equation}
    Now, we expand out the difference: 
    \begin{align}
        L_{i,T}(ax-nby,ay+bx)-L_{i,T}&(ax'-nby',ay'+bx') = (x'-x)(\beta_{i,x} a+\beta_{i,y} b) + (y-y')(n\beta_{i,x}b - \beta_{i,y} a)\\
        & = \frac{M h_{j}}{e_j} \left((\beta_{j,y} a - n\beta_{j,x} b)(\beta_{i,x} a + \beta_{i,y} b) + (\beta_{j,x} a+\beta_{j,y} b)(n\beta_{i,x} b-\beta_{i,y} a)\right)\\
        &= \frac{M h_{j}}{e_j} (\beta_{i,x}\beta_{j,y} -\beta_{j,x}\beta_{i,y})(a^2+nb^2).
    \end{align}
    By the definition of $\delta_{ij}$, we get the desired expression. 
\end{proof}
Finally, before embarking on the proofs of Proposition \ref{prop: type I} and Proposition \ref{prop: type II}, we need to introduce the concatenation result. The first concatenation results were obtained by Tao and Ziegler \cite{TZ:polynomial-progressions}, and were first made quantitative by Peluse and Prendiville \cite{PP:nonlinear-Roth} and Peluse \cite{Peluse:sets_wo_poly_prog}. The statement that we use is the formulation by Green and Sawhney in \cite{Green_Sawhney} of the result of Kravitz, Kuca, and Leng \cite{KravitzKucaLeng}.
\begin{lemma}[{\cite[Lemma 5.5]{Green_Sawhney}}]\label{lem: concatenation}
    Let $\delta\in (0,1/2)$ and $s\geq 1$. Let $I$ be an indexing set and $m$ a power of $2$. Suppose that for $i\in I$ and $j\in \{1,\dots,s\}$, we have symmetric probability measures $\mu_{ji}$ supported on $[\pm N]$. Let $f:\Z\rightarrow\C$ be a 1-bounded function such that 
    \begin{equation}
        \E_{i\in I} \|f\|_{U_{\textrm{GP}}[N;\mu_{1i},\dots,\mu_{si}]}^{2^s} \geq \delta.
    \end{equation}
    Then there exists $t=t_0(m,s)$ such that 
    \begin{equation}
        \E_{i_1,\dots,i_t\in I}\|f\|_{U_{\textrm{GP}}[N;\Omega]}^{2^d} \geq \delta^{O_{m,s}(1)},
    \end{equation}
    where $\Omega$ is the collection of probability measures 
    \begin{equation}
        \Omega = (\mu_{ji_{k_1}}*\dots * \mu_{ji_{k_m}})_{j\in [s], 1\leq k_1<\dots<k_m\leq t}
    \end{equation}
    and $d = s\cdot \binom{t}{m}$.
\end{lemma}

\begin{proof}[Proof of Proposition \ref{prop: type I}]
We will prove the statement of the Proposition for $i=1$; all other cases can be handled analogously. Consider the lower bound of our assumption:
\begin{equation}
    \sum_{ N(\fd)\sim D} \left|\sum_{\substack{N(\fa)\leq X \\ \fd\mid \fa}} w_j(\fa)\right|\geq \delta X;
\end{equation}
by the pigeonhole principle, we can fix an ideal class for $\fd$, i.e. for some $g\in C$, the class group of $K$, we know that 
\begin{equation}
    \sum_{\substack{N(\fd)\sim D\\ [\fd]=g}} \left|\sum_{\substack{N(\fa)\leq X \\ \fd\mid \fa}} w_j(\fa)\right|\gg \delta X.
\end{equation}
Since $w_j$ is supported on principal ideals, we can assume that $[\fa]=g^{-1}$. By \cite[Lemma 2.1]{Green_Sawhney}, we can choose representatives $\mfc_g$ and $\mfc_{g^{-1}}$ such that 
\begin{equation}
    \fd = (a+b\sqrt{-n})\mfc_g, \quad \fb = (x+y\sqrt{-n})\mfc_{g^{-1}},
\end{equation}
where $\fb= \fa/\fd.$ Consequently, we know that 
\begin{equation}
    \fa = (\|g\|((ax-nby)+(ay+bx)\sqrt{-n})),
\end{equation}
where $(\|g\|)=\mfc_g \mfc_{g^{-1}}$. (Note that we can choose these representatives $\mfc_g$ and $\mfc_{g^{-1}}$ such that $\|g\|\in \Z$ and $\|g\| = O_K(1)$. We can also choose that $N(\mfc_g),N(\mfc_{g^{-1}})\in \Z$ by scaling the representatives.)

Let $\Gamma_g = \{(u,v): u+v\sqrt{-n}\in \Pi_g\subset \Z[\sqrt{-n}]\}\subset \Z^2$ be the lattice that is the projection of the lattice $\Pi_g$ (defined in \cite[Lemma 2.1 (iii)]{Green_Sawhney}) to $\Z^2$. There exists an integer $M_g \in \Z$ such that $M_g = O_K(1)$ and $\Gamma_g \subset \frac{1}{M_g}\Z^2$. Therefore, we have that 
\begin{equation}
    \sum_{\substack{a,b\in 1/M_g \Z, |a|,|b|\ll D^{1/2}}} \left|\sum_{\substack{(x,y)\in \Gamma_{g^{-1}} \\ |x|,|y|\ll (X/D)^{1/2}}} w_j(\|g\|((ax-nby)+(ay+bx)\sqrt{-n}))\right| \gg \delta X.
\end{equation}
Next, for a tuple $(e_1,\dots,e_k)\in \N^k,$ define the set 
\begin{equation}
    \mathscr{S}_{(e_1,\dots,e_k)} := \left\{(a,b)\in \frac{1}{M_g}\Z^2: \forall i\in \{1,\dots,k\}, \gcd(M_gM(\beta_{i,x} a+\beta_{i,y} b), M_gM(\beta_{i,y} a - n\beta_{i,x} b))=e_i\right\}.
\end{equation}
(Recall that $M$ is the integer such that $\beta_{i,x},\beta_{i,y},l_i\in \frac{1}{M}\Z$ for each $i$). Since there are $\ll \sum_{e\geq E} D/e^2 \ll D/E$ pairs where for at least one $i$, we have $e_i>E$, we can restrict to the situation where $e_i\ll \delta^{-O(1)}$ for all $i$. By further pigeonholing, there is some tuple $\Vec{e}=(e_1,\dots,e_k)$ such that $e_i\ll \delta^{-O(1)}$ for each $i$ and 
\begin{equation}
    \sum_{\substack{a,b\in 1/M_g\Z \\ |a|,|b|\ll D^{1/2}\\ (a,b)\in \mathscr{S}_{\Vec{e}} }} \left|\sum_{\substack{(x,y)\in \Gamma_{g^{-1}} \\ |x|,|y|\ll (X/D)^{1/2}}} w_j(\|g\|((ax-nby)+(ay+bx)\sqrt{-n}))\right| \gg \delta^{O(1)} X.
\end{equation} 

It is now time to expand the definition of $w_j$: recall that 
\begin{equation}
    w_j(\fa) = \prod_{i=1}^k w_{ij}(\fa), \quad w_{ij}(\fa) = (\log X)^{-1} \sum_{\substack{x,y\in \Z \\ \fa = (T_x(x)+T_y(y)\sqrt{-n}))}} \lambda_{ij}(L_i(x,y)).
\end{equation}
Let $f_{ij} = (\log X)^{-1} \lambda_{ij}.$ Expanding the definition of $w_j$ gives us 
\begin{equation}
    \sum_{\substack{a,b\in 1/M_g\Z, |a|,|b|\ll D^{1/2}\\ (a,b)\in \mathscr{S}_{\Vec{e}} }} \left|\sum_{\substack{(x,y)\in \Gamma_T \\ |x|,|y|\ll (X/D)^{1/2}}} \prod_{i=1}^k f_{ij}(\|g\| \cdot L_{i,T}(ax-nby, ay+bx))\right| \gg \delta^{O(1)} X,
\end{equation}
where $\Gamma_T$ is the sublattice of $\Gamma_{g^{-1}}$ where $\|g\|((ax-nby)+(ay+bx)\sqrt{-n})$ is contained in $\{T_x(u)+T_y(v)\sqrt{-n}: u+v\sqrt{-n}\in \Z[\sqrt{-n}]\}$. Observe that $\Gamma_T$ is a lattice of determinant $O_{K,T_x,T_y}(1).$

We now apply many iterations of Cauchy-Schwarz. By the 1-boundedness of $f_{kj}$, we know that 
\begin{equation}
    \sum_{\substack{a,b\in 1/M_g\Z, |a|,|b|\ll D^{1/2}\\ (a,b)\in \mathscr{S}_{\Vec{e}} }} \sum_{\substack{t_k\in \frac{1}{M_gM }\Z\\ |t_k|\ll X^{1/2}}}\left|\sum_{\substack{(x,y)\in \Gamma_{T} \\ |x|,|y|\ll (X/D)^{1/2}\\ L_{k,T}(ax-nby,ay+bx)=t_k}} \prod_{i=1}^{k-1}f_{ij}(\|g\| \cdot L_{i,T}(ax-nby, ay+bx))\right| \gg \delta^{O(1)} X.
\end{equation}
After one application of Cauchy-Schwarz, we have the bound:
\begin{multline}
    \sum_{\substack{a,b\in 1/M_g\Z, |a|,|b|\ll D^{1/2}\\ (a,b)\in \mathscr{S}_{\Vec{e}} }} \sum_{\substack{(x,y),(x',y')\in \Gamma_{T} \\ |x|,|y|, |x'|,|y'|\ll (X/D)^{1/2}}}\mathbf{1}_{L_{k,T}(ax-nby,ay+bx)=L_{k,T}(ax'-nby',ay'+bx')} \\\times \prod_{i=1}^{k-1}f_{ij}(\|g\| \cdot L_{i,T}(ax-nby, ay+bx))f_{ij}(\|g\| \cdot L_{i,T}(ax'-nby', ay'+bx')) \gg \delta^{O(1)} \frac{X^{3/2}}{D}.
\end{multline}
By Lemma \ref{lem: differencing and algebraic manipulations}, we know that since $L_{k,T}(ax-nby,ay+bx) = L_{k,T}(ax'-nby',ay'+bx')$, we have that for some $h_k\in \Z$,
\begin{equation}
    L_{i,T}(ax'-nby',ay'+bx')-L_{i,T}(ax-nby,ay+bx) = \frac{M_g M \delta_{ik}}{e_i} (a^2+nb^2) \cdot h_k.
\end{equation}
Let $M_{ik} = M_g M \delta_{ik} \|g\| e_i^{-1}$. We can therefore collapse the terms:
\begin{equation}
     \sum_{\substack{a,b\in 1/M_g\Z, \\ |a|,|b|\ll D^{1/2}\\ (a,b)\in \mathscr{S}_{\Vec{e}} }} \sum_{|h_k| \ll X^{1/2}/D} \left|\sum_{\substack{(x,y)\in \Z^2 \\ |x|,|y|\ll (X/D)^{1/2}}} \prod_{i=1}^{k-1} (\Delta_{M_{ik}(a^2+nb^2)h_k}f_{ij})(\|g\| \cdot L_{i,T}(ax-nby,ay+bx))\right| \gg \delta^{O(1)} \frac{X^{3/2}}{D}.
\end{equation}
(We remark that we can replicate the arguments of \cite[p.29]{Green_Sawhney} to also simultaneously remove the condition that $(x,y)\in \Gamma_T$ during this stage via expanding $\mathbf{1}_{\Gamma_T}$ as in \cite[(6.16)]{Green_Sawhney}, applying Cauchy-Schwarz with the phase, and eliminating the phase with the final set of absolute values). Since we assumed that the homogeneous parts of $L_i$ and $L_k$ are linear independent, we know that $\delta_{ik}\neq 0$ and hence $M_{ik}\neq 0.$

Note that since $f_{ij}$ is 1-bounded, $\Delta_{M_{ik}(a^2+nb^2)h_k}f_{ij}$ is 2-bounded. We can therefore iterate the above argument $k-1$ times to achieve the lower bound:
\begin{multline}
    \E_{\substack{a,b\in \frac{1}{M_g}\Z \\ |a|,|b|\ll D^{1/2}\\ (a,b)\in \mathscr{S}_{\Vec{e}}}} \E_{|h_2|\ll X^{1/2}/D}\dots \E_{|h_k|\ll X^{1/2}/D} \\ \left|\E_{\substack{|x|,|y|\ll (X/D)^{1/2}}} (\Delta_{M_{12}(a^2+nb^2)h_2}\dots \Delta_{M_{k2}(a^2+nb^2)h_k} f_{1j})(\|g\|\cdot L_{1,T}(ax-nby,ay+bx)) \right|\gg \delta^{O(1)}.
\end{multline}
Observe that since the homogeneous parts of $L_1$ and $L_i$ are linearly independent for all $i\neq 1$, we know that $M_{1i}\neq 0$ for each $i=2,\dots,k$. 

This tells us that there are $\gg \delta^{O(1)} D (X^{1/2}/D)^{k-1}$ tuples $(a,b,h_2,\dots,h_k)$ such that $(a,b)\in \mathscr{S}_{\Vec{e}}$ such that 
\begin{equation}
    \E_{\substack{|x|,|y|\ll (X/D)^{1/2}}} (\Delta_{M_{12}(a^2+nb^2)h_2}\dots \Delta_{M_{k2}(a^2+nb^2)h_k} f_{1j})(\|g\|\cdot L_{1,T}(ax-nby,ay+bx)) \gg \delta^{O(1)}.
\end{equation}
We can further assume that $(a,b)$ satisfies that if 
\begin{equation}
    L_{1,T}(ax-nby,ay+bx) = \psi_1(a,b) x + \psi_2(a,b) y + l',
\end{equation}
then $|\psi_1(a,b)|,|\psi_2(a,b)|\gg \delta^{O(1)} D^{1/2}.$ We can therefore apply \cite[Lemma 5.4]{Green_Sawhney} as in \cite[p. 30]{Green_Sawhney} to see that 
\begin{equation}\label{eq: type I pre linearization}
    \E_{\substack{a,b\in \frac{1}{M}\Z \\|a|,|b|\ll D^{1/2} \\ (a,b)\in \mathscr{S}_{\Vec{e}}}} \E_{|h_2|,\dots,|h_k|\ll X^{1/2}/D} \|\Delta_{M_{12}(a^2+nb^2)h_2}\dots \Delta_{M_{1k}(a^2+nb^2)h_k} f_{1j}\|^2_{U^2[X^{1/2}]} \gg \delta^{O(1)}.
\end{equation}
(Note that although \cite[Lemma 5.4]{Green_Sawhney} is stated for linear forms rather linear polynomials such as $L_{1,T}$, this can remedied by shifting the intervals $I_1$ and $I_2$).

We now proceed with our ``linearization'' step. Let us now define for $t\in \Z$: 
\begin{equation}
    \sigma(t) := \#\{(a,b)\in \mathscr{S}_{\Vec{e}}: |a|,|b|\ll D^{1/2}, a^2+nb^2 = t\}.
\end{equation}
We can rewrite \eqref{eq: type I pre linearization} as: 
\begin{equation}
     \E_{|h_2|,\dots,|h_k|\ll X^{1/2}/D} \E_{|t|\ll D} \sigma(t) \|\Delta_{M_{12} t h_2}\dots \Delta_{M_{1k}th_k} f_{1j}\|_{U^2[X^{1/2}]}^2 \gg \delta^{O(1)}.
\end{equation}
Apply Cauchy-Schwarz and see that 
\begin{equation}
     \left(\E_{|h_2|,\dots,|h_k|\ll X^{1/2}/D} \E_{|t|\ll D} \sigma(t)^2 \right)\left(\E_{|h_2|,\dots,|h_k|\ll X^{1/2}/D} \E_{|t|\ll D}  \|\Delta_{M_{12} t h_2}\dots \Delta_{M_{1k}th_k} f_{1j}\|_{U^2[X^{1/2}]}^4\right)\gg \delta^{O(1)}.
\end{equation}
Since $\#\{(a,b)\in \frac{1}{M_g}\Z: a^2+nb^2 = t\} \ll_{K} \tau(t)$ for any $t\in \Z_{\geq 0}$, it is straightforward to see that 
\begin{equation}
    \E_{|h_2|,\dots,|h_k|\ll X^{1/2}/D} \E_{|t|\ll D} \sigma(t)^2  \ll \E_{|t|\ll D} \tau(t)^2 \ll (\log X)^3.
\end{equation}
Here we have used that $D\leq X^{1/2}.$ Thus, this tells us that 
\begin{equation}
    \E_{|h_2|,\dots,|h_k|\ll X^{1/2}/D} \E_{|t|\ll D}  \|\Delta_{M_{12} t h_2}\dots \Delta_{M_{1k}th_k} f_{1j}\|_{U^2[X^{1/2}]}^4 \gg \frac{\delta^{O(1)}}{(\log X)^3}. 
\end{equation}
At this point, we no longer have any quadratic forms so we have, in essence, linearized our problem. 

After expanding the $U^2[X^{1/2}]$ norm, we have that 
\begin{equation}
    \E_{y,y'\sim X^{1/2}}\E_{|t|\ll D} \|\Delta_y \Delta_{y'}f_{1j}\|_{U_{\textrm{GP}}[X^{1/2};\mu_{2,t},\dots,\mu_{k,t}]} \gg \frac{\delta^{O(1)}}{(\log X)^3},
\end{equation}
where our probability measures are defined as $\mu_{j,t}$ is the uniform probability measure on $M_{1j}t [\pm X^{1/2}/D].$ (Recall the definition of $U_{\textrm{GP}}$ in Definition \ref{defn: gowers peluse norm}).

\begin{lemma}\label{lem: type I convolution of measures}
    The following bound holds for any $j\in \{2,\dots,k\}$: 
    \begin{equation}
        \E_{t_1,t_2,t_3\ll D} \|\mu_{j,t_1}*\mu_{j,t_2} * \mu_{j,t_3}\|_2^2 \ll X^{-1/2}.
    \end{equation}
\end{lemma}
\begin{proof}
    First, note that 
    \begin{equation}
        \|\mu_{j,t_1}*\mu_{j,t_2} * \mu_{j,t_3}\|_2^2 \ll (\tilde{\mu}_{j,t_1}*\tilde{\mu}_{j,t_2} * \tilde{\mu}_{j,t_3})(0),
    \end{equation}
    where $\tilde{\mu}_{j,t_i}$ is the uniform probability measure on $M_{1j}t[\pm 2X^{1/2}/D]$ (see \cite[p.35]{Green_Sawhney}). Let us expand what this expression means: 
    \begin{multline}\E_{t_1,t_2,t_3\ll D} (\tilde{\mu}_{j,t_1}*\tilde{\mu}_{j,t_2} * \tilde{\mu}_{j,t_3})(0) \\ = \frac{1}{D^{3}} \cdot \frac{1}{(X^{1/2}/D)^3} \#\{(t_1,t_2,t_3) \in [\pm D]^3, (h_1,h_2,h_3)\in [\pm 2X^{1/2}/D]^3: M_{1j}(t_1h_1 +  t_2 h_2 + t_3 h_3 )= 0\}.
    \end{multline}
    The expression $t_1h_1+t_2h_2+t_3h_3$ is a bilinear form -- that is, it is both linear in $(t_1,t_2,t_3)$ and linear in $(h_1,h_2,h_3)$. As such, by geometry of numbers estimates, we have that both of the following estimates hold: 
    \begin{multline}
        \#\{(t_1,t_2,t_3) \in [\pm D]^3, (h_1,h_2,h_3)\in [\pm 2X^{1/2}/D]^3: M_{1j}(t_1h_1 +  t_2 h_2 + t_3 h_3 )= 0\} \\
        \ll  \begin{cases}\sum_{|t_1|,|t_2|,|t_3|\leq D} \frac{X}{D^2} \cdot \frac{1}{\|(t_1,t_2,t_3)\|_2} + \frac{X^{1/2}}{D \lambda_1((t_1,t_2,t_3))} + 1 \\ 
        \sum_{|h_1|,|h_2|,|h_3|\leq 2X^{1/2}/D} \frac{D^2}{\|(h_1,h_2,h_3)\|_2} + \frac{D}{\lambda_1((h_1,h_2,h_3))} +1.
        \end{cases}
    \end{multline}
    Here $\lambda_1(\Vec{v})$ denotes the first successive minima of the lattice $\{\Vec{x}\in \Z^3: \Vec{v}\cdot \Vec{x}=0\}$. For our purposes, we will only use that $\lambda_1(\Vec{v})\geq 1$. Consequently, we have that 
    \begin{equation}
        \E_{t_1,t_2,t_3\ll D} (\tilde{\mu}_{j,t_1}*\tilde{\mu}_{j,t_2} * \tilde{\mu}_{j,t_3})(0) \ll X^{-3/2} \min\left(X+X^{1/2}D^2+D^3, X + X^{3/2}D^{-2}+X^{3/2}D^{-3}\right).
    \end{equation}
    Note that the above is always bounded by $\ll X^{-1/2}$ as claimed. 
    
\end{proof}

By Lemma \ref{lem: type I convolution of measures} and Lemma \ref{lem: concatenation} with $s=3$ (the concatenation lemma), we have that 
\begin{equation}
    \E_{y,y'\sim X^{1/2}} \E_{|t_1|,|t_2|,|t_3|\ll D} \|\Delta_y \Delta_{y'} f_{1j} \|_{U_{\textrm{GP}}[X^{1/2}; \Omega]} \gg \left(\frac{\delta}{\log X}\right)^{O(1)},
\end{equation}
where $\Omega = (\mu_{i,t_1}*\mu_{i,t_2}*\mu_{i,t_3})_{i\in \{2,\dots,k\}, |t_1|,|t_2|,|t_3|\ll D}$.
Now by Corollary \ref{cor: Gowers Peluse to Gowers}, since $$\E_{|t_1|,|t_2|,|t_3|\ll D}\|\mu_{i,t_1}*\mu_{i,t_2}*\mu_{i,t_3}\|_2^2 \ll X^{-1/2}$$ for all $i\in \{2,\dots,k\}$, this gives us that for some constant $r>0$,
\begin{equation}
    \E_{y,y'\sim X^{1/2}} \|\Delta_y \Delta_{y'} f_{1j}\|_{U^r[X^{1/2}]}^{2^r} \gg \left(\frac{\delta}{\log X}\right)^{O(1)}. 
\end{equation}
We can now expand the Gowers norm and see that the above is equivalent to 
\begin{equation}
    \|f_{1j}\|_{U^{c_k}[X^{1/2}]}^{2^{c_k}} \gg \left(\frac{\delta}{\log X}\right)^{O(1)},
\end{equation}
where $c_k=r+2$. Recalling that $f_{1j} = (\log X)^{-1} \lambda_{1j}$, this gives us the desired conclusion that 
\begin{equation}
    \|\lambda_{1j}\|_{U^{c_k}[X^{1/2}]} \gg \left(\frac{\delta}{\log X}\right)^{O(1)}. 
\end{equation}
This completes the proof of Proposition \ref{prop: type I}.
\end{proof}

\section{Type II estimates}\label{sec: primes type II}
In this section we prove Proposition \ref{prop: type II}; the initial steps of this proof will be very similar to the previous section. It is mostly at Lemma \ref{lem: type I convolution of measures} that the proofs diverge; in this proof we are not able to easily ``linearize'' the problem as in the Type I computations.

\begin{proof}[Proof of Proposition \ref{prop: type II}]
We prove the statement of the Proposition for $i=1$; all other cases can be handled analogously. As in the Type I setting, we can fix the ideal class for $\fa$ and $\fb$; that is, we can take some $g\in C$ (the class group of $K$) and assume that $[\fa]=g$ and $[\fb] =g^{-1}$. 
Since 
\begin{equation}
    \sum_{\substack{N(\fa)\leq X \\ [\fa]=g,[\fb]=g^{-1}\\ N(\fa)\sim D}} \alpha_\fa\beta_\fb w_j(\fa\fb) \geq \delta X,
\end{equation}
and since $\alpha,\beta:\calI(\msO_K)\rightarrow \C$ are 1-bounded functions, after an application of Cauchy-Schwarz we have the lower bound:
\begin{equation}
    \sum_{\substack{ N(\fa),N(\fa')\sim D \\ [\fa]=[\fa']=g}} \left|\sum_{\substack{N(\fa\fb),N(\fa'\fb) \leq X \\ [\fb] = g^{-1}}} \prod_{i=1}^k w_{ij}(\fa\fb) w_{ij}(\fa'\fb)\right| \gg \delta^{O(1)}XD.
\end{equation}
From \cite[Lemma 2.1]{Green_Sawhney}, we again choose representatives $\mfc_g,\mfc_{g^{-1}}$ such that $N(\mfc_g),N(\mfc_{g^{-1}})\in \Z$ and $(\|g\|)= \mfc_g\mfc_{g^{-1}}$ for $\|g\|\in \Z$. We can therefore write:
\begin{equation}\label{eq: type II expansion of ideal in coord}
    \fa = (a+b\sqrt{-n})\mfc_g, \quad \fa' = (a'+b'\sqrt{-n})\mfc_g, \quad \fb = (x+y\sqrt{-n}) \mfc_{g^{-1}}
\end{equation}

Let $\Gamma_g = \{(u,v): u+v\sqrt{-n}\in \Pi_g \subset \Z[\sqrt{-n}]\}$ be the lattice that is the projection of the lattice $\Pi_g$ (defined in \cite[Lemma 2.1]{Green_Sawhney}) to $\Z^2$. Again, let $f_{ij} = (\log X)^{-1}\lambda_{ij}$. Expanding the definition of $w_{ij}$ and using \eqref{eq: type II expansion of ideal in coord}, we have:
\begin{multline}
    \sum_{\substack{(a,b),(a',b')\in \Gamma_g \\ |a|,|b|,|a'|,|b'|\sim D^{1/2}}} \left|\sum_{\substack{(x,y)\in \Gamma_{g^{-1}}\\ |x|,|y|\ll (X/D)^{1/2}}} \prod_{i=1}^k f_{ij}(\|g\| L_{i,T}(ax-nby, ay+bx)) f_{ij}(\|g\| L_{i,T}(a'x-nb'y, a'y+b'x))\right|\\  \gg \delta^{O(1)} XD.
\end{multline}

Let $\calQ'$ denote a set of tuples $(a,b,a',b')$ with $|a|,|b|,|a'|,|b'| \ll D^{1/2}$ and such that 
\begin{equation}
    \left|\sum_{\substack{(x,y)\in \Gamma_{g^{-1}}\\ |x|,|y|\ll (X/D)^{1/2}}}\prod_{i=1}^k f_{ij}(\|g\| L_{i,T}(ax-nby, ay+bx)) f_{ij}(\|g\| L_{i,T}(a'x-nb'y, a'y+b'x))\right| \gg  \frac{\delta^{O(1)} X}{D}. 
\end{equation}
It is straightforward to see that $|\calQ'|\gg D^{2}$. 
We can further restrict to $\calQ\subset \calQ'$ where we add the additional assumptions that for each $i$, 
\begin{equation}
    a\beta_{i,x} + b\beta_{i,y}\neq 0, \quad a'\beta_{i,x}+b'\beta_{i,y}\neq 0, \quad a\beta_{i,y} - nb\beta_{i,x} \neq 0, \quad a'\beta_{i,y} - nb'\beta_{i,x}\neq 0,
\end{equation}
where we recall that $L_{i,T}(x,y) = \beta_{i,x} x + \beta_{i,y} y + l_i\in \frac{1}{M_g}\Z[x,y]$. We also add the assumptions that $aa'+nbb'\neq 0$ and $ab'-a'b\neq 0$ and that 
\begin{equation}
    (aa'+nbb')\delta_{1i} + (a'b-ab')\delta_{1i}'\neq 0
\end{equation} for $(a,b,a',b')\in \calQ$. We make a final assumption that if 
\begin{equation}
    L_{1,T}(ax-nby,ay+bx) = \psi_1(a,b) x + \psi_2(a,b) y + l',
\end{equation}
then $|\psi_1(a,b)|,|\psi_2(a,b)|\gg \delta^{O(1)}D^{1/2}.$
Again, since all of these conditions are generic, we can see that $|\calQ|\gg D^2$.

For now, fix a tuple $(a,b,a',b')\in \calQ$. Let $H = [X^{1/2}/D]$. We know that for any $|h_k|\ll H$, 
    \begin{align}
        L_{k,T}(ax-nby,ay+bx) &= L_{k,T}(a\tilde{x}-nb\tilde{y}, a\tilde{y}+b\tilde{x})\\
        &= (a\beta_{k,x}+b\beta_{k,y}) \tilde{x}+ (a\beta_{k,y}-nb\beta_{k,x})\tilde{y} + l_k,
    \end{align}
for $\tilde{x} = x + (a\beta_{k,y}-nb\beta_{k,x})h_k$ and $\tilde{y} = y - (a\beta_{k,x}+b\beta_{k,y})h_k$. Since $(a,b,a',b')\in \calQ$, we know that $a\beta_{k,y}-nb\beta_{k,x}$ and $a\beta_{k,x}+b\beta_{k,y}$ are nonzero.  
Hence, we have that there are elements $c_{h_k} \in \R$ with $|c_{h_k}|= 1$ such that 
\begin{multline}\label{eq: type II added average over h_k}
    \sum_{\substack{(x,y)\in \Gamma_{g^{-1}}\\ |x|,|y|\ll (X/D)^{1/2}}}\left|\sum_{|h_k|\ll H} c_{h_k} \prod_{i=1}^{k-1} f_{ij}(\|g\| L_{i,T}(a\tilde{x}-nb\tilde{y}, a\tilde{y}+b\tilde{x})) \prod_{i=1}^k f_{ij}(\|g\| L_{i,T}(a'\tilde{x}-nb'\tilde{y}, a'\tilde{y}+b'\tilde{x}))\right| \\ \gg \frac{\delta^{O(1)}XH}{D}.
\end{multline}
Observe that we can drop the condition that $(x,y)\in \Gamma_{g^{-1}}$; we can do so using by detecting $(x,y)\in \Gamma_{g^{-1}}$ as in \cite[(6.16)]{Green_Sawhney}. 

In a similar vein to Lemma \ref{lem: differencing and algebraic manipulations}, we know that for $\tilde{x} = x + (a\beta_k-nb\alpha_k)h_k$ and $\tilde{y} = y - (a\alpha_k+b\beta_k)h_k$,
\begin{multline}
    L_{i,T}(a\tilde{x}-nb\tilde{y}, a\tilde{y}+b\tilde{x}) - L_{i,T}(ax-nby, ay+bx) \\
    = h_k \left(a\beta_{i,x}(a\beta_{k,y}-nb\beta_{k,x}) + nb\beta_{i,x}(a\beta_{k,x}+b\beta_{k,y}) - a\beta_{i,y}(a\beta_{k,x}+b\beta_{k,y}) + b\beta_{i,y}(a\beta_{k,y}-nb\beta_{k,x})\right)\\
    = h_k (a^2 +nb^2) \delta_{ik}
\end{multline}
On the other hand, 
\begin{multline}
    L_{i,T}(a'\tilde{x} - nb'\tilde{y}, a'\tilde{y} + b'\tilde{x}) - L_{i,T}(a'x-nb'y, a'y+b'x) \\
    = \beta_{i,x} (a'(a\beta_{k,y}-nb\beta_{k,x})h_k +nb'(a\beta_{k,x}+b\beta_{k,y})h_k) + \beta_{i,y}(-a'(a\beta_{k,x}+b\beta_{k,y})h_k + b'(a\beta_{k,y}-nb\beta_{k,x})h_k) \\
    = ((aa'+nbb')\delta_{ik} + (ab'-a'b) \delta_{ik}')h_k.
\end{multline}
(Recall that $\delta_{ij}' = n\beta_{i,x}\beta_{j,x}+\beta_{i,y}\beta_{j,y}$.)
By our choice of $(a,b,a',b')\in \calQ$, we know that $a^2+nb^2, aa'+nbb', ab'-a'b\neq 0$. 

Thus, we can rewrite \eqref{eq: type II added average over h_k} as 
\begin{multline}
     \sum_{|x|,|y|\ll (X/D)^{1/2}}\Bigg|\sum_{|h_k|\ll H} c_{h_k} \prod_{i=1}^{k-1} f_{ij}(\|g\| (L_{i,T}(ax-nby, ay+bx)+ \delta_{ik}(a^2+nb^2)h_k)) \\ \prod_{i=1}^k f_{ij}(\|g\| (L_{i,T}(a'x-nb'y, a'y+b'x) + ((aa'+nbb')\delta_{ik} + (ab'-a'b)\delta_{ik}')h_k))\Bigg|\gg \frac{\delta^{O(1)}XH}{D}.
\end{multline}
After an application of Cauchy-Schwarz, we achieve 
\begin{multline}
    \sum_{|h_k|,|h_k'|\ll H} \Bigg|\sum_{|x|,|y|\ll (X/D)^{1/2}} \prod_{i=1}^{k-1} \Delta_{\delta_{ik}(a^2+nb^2)(h_k,h_k')} f_{ij}(\|g\| L_{i,T}(ax-nby,ay+bx)) \\ \times \prod_{i=1}^k \Delta_{((aa'+nbb')\delta_{ik} + (ab'-a'b)\delta_{ik}')(h_k,h_k')}f_{ij}(\|g\| L_{i,T}(a'x-nb'y, a'y+b'x))\Bigg| \gg \frac{\delta^{O(1)}XH^2}{D}.
\end{multline}
(Recall the definition of Gowers-Peluse norms and $\Delta_{(h,h')}$ in Definition \ref{defn: gowers peluse norm}.)

We can now iterate the procedure $k-1$ times and obtain the following bound:
\begin{multline}
    \sum_{\substack{|h_2|,|h_2'|\ll H}}\dots \sum_{|h_k|,|h_k'|\ll H} \Bigg|  \sum_{|x|,|y|\ll (X/D)^{1/2}}\prod_{i=1}^k \Delta_{((aa'+nbb')\delta_{12} + (ab'-a'b)\delta_{12}')(h_2,h_2')} \\ \cdots \Delta_{((aa'+nbb')\delta_{ik} + (ab'-a'b)\delta_{ik}')(h_k,h_k')}f_{ij}(\|g\| L_{i,T}(a'x-nb'y, a'y+b'x)) \\  \Delta_{\delta_{12}(a^2+nb^2)(h_2,h_2')} \dots \Delta_{\delta_{1k}(a^2+nb^2)(h_k,h_k')} f_{1j}(\|g\| L_{1,T}(ax-nby,ay+bx)) \Bigg| \gg \frac{\delta^{O(1)}XH^{2(k-1)}}{D}.
\end{multline}

Next, let $\eta_k$ range in $|\eta_k|\ll H$. Then for $\tilde{x} = x+(a'\beta_{k,y}-nb'\beta_{k,x})\eta_k$ and $\tilde{y} = y -(a'\beta_{k,x}+b'\beta_{k,y})\eta_k$, we similarly see that 
\begin{equation}
    L_{k,T}(a'x-nb'y, a'y+b'x) = L_{k,T}(a'\tilde{x}-nb'\tilde{y}, a'\tilde{y} + b'\tilde{x}).
\end{equation}
Therefore, by the analogous computation, we have that 
\begin{multline}
    \sum_{|\eta_k|,|\eta_k'|\ll H} \sum_{\substack{|h_2|,|h_2'|\ll H}}\dots \sum_{|h_k|,|h_k'|\ll H} \Bigg|  \sum_{|x|,|y|\ll (X/D)^{1/2}}\prod_{i=1}^{k-1} \Delta_{\delta_{ik}(a'^2+nb'^2)(\eta_k,\eta_k')}\Delta_{((aa'+nbb')\delta_{12} + (ab'-a'b)\delta_{12}')(h_2,h_2')} \\ \cdots \Delta_{((aa'+nbb')\delta_{ik} + (ab'-a'b)\delta_{ik}')(h_k,h_k')}f_{ij}(\|g\| L_{i,T}(a'x-nb'y, a'y+b'x)) \\  \Delta_{((aa'+nbb')\delta_{1k} + (a'b-ab')\delta_{1k}')(\eta_k,\eta_k')} \Delta_{\delta_{12}(a^2+nb^2)(h_2,h_2')} \dots \Delta_{\delta_{1k}(a^2+nb^2)(h_k,h_k')} f_{1j}(\|g\| L_{1,T}(ax-nby,ay+bx)) \Bigg| \\ \gg \frac{\delta^{O(1)}XH^{2k}}{D}.
\end{multline}
We now iterate this procedure $k$ times to see that for this choice of $(a,b,a',b')\in \calQ$, 
\begin{multline}
    \sum_{|\eta_1|,|\eta_1'|\ll H} \dots \sum_{|\eta_k|,|\eta_k'|\ll H} \sum_{\substack{|h_2|,|h_2'|\ll H}}\dots \sum_{|h_k|,|h_k'|\ll H} \\
    \Bigg| \sum_{|x|,|y|\ll (X/D)^{1/2}}
    \Delta_{((a'b-ab')\delta_{11}')(\eta_1,\eta_1')}\Delta_{((aa'+nbb')\delta_{12} + (a'b-ab')\delta_{12}')(\eta_2,\eta_2')} \cdots \Delta_{((aa'+nbb')\delta_{1k} + (a'b-ab')\delta_{1k}')(\eta_k,\eta_k')}  \\
      \Delta_{\delta_{12}(a^2+nb^2)(h_2,h_2')} \dots \Delta_{\delta_{1k}(a^2+nb^2)(h_k,h_k')} f_{1j}(\|g\| L_{1,T}(ax-nby,ay+bx))
    \Bigg|\gg \frac{\delta^{O(1)}XH^{2(2k-1)}}{D}.
\end{multline}

For $\gg \delta^{O(1)}H^{2(2k-1)}$ tuples $(\eta_1,\eta_1',\dots,\eta_k,\eta_k',h_2,h_2',\dots, h_k,h_k')$, this internal summation is $\gg \delta^{O(1)} X/D.$ By our assumption on $(a,b,a',b')\in \calQ$, we can apply \cite[Lemma 5.4]{Green_Sawhney} and see that for such values of $(\eta_1,\eta_1',\dots,\eta_k,\eta_k',h_2,h_2',\dots, h_k,h_k')$, 
\begin{multline}
    \|\Delta_{((a'b-ab')\delta_{11}')(\eta_1,\eta_1)'}\Delta_{((aa'+nbb')\delta_{12} + (a'b-ab')\delta_{12}')(\eta_2,\eta_2')} \dots \Delta_{((aa'+nbb')\delta_{1k} + (a'b-ab')\delta_{1k}')(\eta_k,\eta_k')}  \\
      \Delta_{\delta_{12}(a^2+nb^2)(h_2,h_2')} \dots \Delta_{\delta_{1k}(a^2+nb^2)(h_k,h_k')} f_{1j}\|_{U^2[X^{1/2}]}\gg \delta^{O(1)}.
\end{multline}
Define the measures $\mu_{i,q}$ to be the uniform probability measure on $((aa'+nbb')\delta_{1i} + (a'b-ab')\delta_{1i}')[\pm H]$ and $\nu_{i,q}$ denotes the uniform probability measure on $(a^2+nb^2)\delta_{1i} [\pm H].$ 
By the definition of the Gowers-Peluse norm (Definition \ref{defn: gowers peluse norm}), and noting that this holds for $\gg \delta^{O(1)}H^{2(2k-1)}$ tuples $(\eta_1,\eta_1',\dots,\eta_k,\eta_k',h_2,h_2',\dots, h_k,h_k')$ and $\gg \delta^{O(1)} D^2$ tuples $(a,b,a',b')$ with $|a|,|b|,|a'|,|b'|\ll D^{1/2}$, this implies that 
\begin{equation}
    \E_{y,y'\sim X^{1/2}} \E_{\substack{q=(a,b,a',b')\in \calQ\\ |a|,|b|,|a'|,|b'|\sim D^{1/2}}} \|\Delta_y\Delta_{y'} f_{1j}\|_{U_{\text{GP}}[X^{1/2};\mu_{1,q},\dots,\mu_{k,q},\nu_{2,q},\dots,\nu_{k,q}]} \gg \delta^{O(1)}. 
\end{equation}
(We refer to \cite[(7.9)-(7.12)]{Green_Sawhney} for a more detailed explanation of this process). 

We claim that the measures $\mu_{i,q}$ and $\nu_{i,q}$ satisfy that a sufficient number of convolutions with themselves will be close to uniform measures. One can think of the below lemma as the analog of Lemma \ref{lem: type I convolution of measures} in this setting and note that this is a weaker statement. 

\begin{lemma}\label{lem: type II convolution of measures}
Let $m\geq 1$. Let $i\in \{1,\dots k\}$, then 
\begin{equation}
    \E_{q_1,\dots,q_m\in \calQ} \|\mu_{i,q_1} * \dots \mu_{i,q_m}\|_2^2 \ll \delta^{O(1)}\max(H^{-m}, X^{-1/2}) (\log X)^{O_m(1)}.
\end{equation}
Let $i\in \{2,\dots,k\}$, then 
\begin{equation}
    \E_{q_1,\dots,q_m\in \calQ} \|\nu_{i,q_1} * \dots \nu_{i,q_m}\|_2^2 \ll \delta^{O(1)}\max(H^{-m},X^{-1/2})(\log X)^{O_m(1)}.
\end{equation}
\end{lemma}
Before completing the proof of the lemma, let us now finish the proof of Proposition \ref{prop: type II}. 
Recall that $H=\lfloor X^{1/2}/D\rfloor$ and $D\in [\delta^{O(1)},X^{3/8}]$. Therefore, for $m>4$, we have that 
\begin{equation}
    \E_{q_1,\dots,q_m\in \calQ} \|\mu_{i,q_1} * \dots \mu_{i,q_m}\|_2^2 \ll_m \delta^{O(1)} X^{-1/2}(\log X)^{O_m(1)}
\end{equation}
\begin{equation}
    \E_{q_1,\dots,q_m\in \calQ} \|\nu_{i,q_1} * \dots \nu_{i,q_m}\|_2^2 \ll_m \delta^{O(1)} X^{-1/2}(\log X)^{O_m(1)}.
\end{equation}
We therefore apply Lemma \ref{lem: concatenation} with $m=5$ to see that 
\begin{equation}
    \E_{y,y'\sim X^{1/2}} \E_{q_1,\dots,q_5\in \calQ} \|\Delta_y\Delta_{y'} f_{1j}\|_{U_{\text{GP}}[X^{1/2};\Omega]} \gg \left(\frac{\delta}{\log X}\right)^{O(1)}.
\end{equation}
By Corollary \ref{cor: Gowers Peluse to Gowers} and Lemma \ref{lem: type II convolution of measures}, this gives us that for some $r>0$ (depending on $k$), 
\begin{equation}
    \E_{y,y'\sim X^{1/2}} \|\Delta_y \Delta_{y'} f_{1j}\|_{U^r[X^{1/2}]}^{2^r}\gg \left(\frac{\delta}{\log X}\right)^{O(1)}. 
\end{equation}
Expanding the Gowers norm and collapsing again gives for $r_k=r+2$ that 
\begin{equation}
    \|f_{1j}\|_{U^{r_k}[X^{1/2}]}^{2^{r_k}} \gg \left(\frac{\delta}{\log X}\right)^{O(1)}. 
\end{equation}
Finally, recall that $f_{1j} = (\log X)^{-1}\lambda_{1j}$, we achieve the desired lower bound on $\|\lambda_{1j}\|_{U^{r_k}[X^{1/2}]}.$
This completes the proof of Proposition \ref{prop: type II}. 
\end{proof}

\begin{proof}[Proof of Lemma \ref{lem: type II convolution of measures}]
For both statements, we follow the proof of \cite[Lemma 7.1]{Green_Sawhney}
Define the following measures: 
\begin{equation}
    \mu' = \E_{q\in \calQ} \mu_{i,q},\quad \nu' = \E_{q\in \calQ} \nu_{i,q}.
\end{equation}
Then since for any measure $\mu$, $\|\mu\|_2^2 \ll \mu * \mu (0)$, and $\mu_{i,q}$ and $\nu_{i,q}$ are uniform probability measures, as in \cite[Lemma 7.1]{Green_Sawhney}, the problem reduces to showing that 
\begin{equation}\label{eq: type I convolution bound on mu' nu'}
    \mu'^{(m)}(0), \nu'^{(m)} (0) \ll_m \delta^{O(1)} \max(H^{-m},X^{-1/2})(\log X)^{O_m(1)}.
\end{equation}
Now we define $\mu_*$ (resp. $\nu_*$) to satisfy that $\mu' \leq \delta_0 H^{-1} + \mu_*$ (resp. $\nu'\leq \delta_0 H^{-1} + \nu_*$), where $\mu_*$ is uniform on the set $\{M_g\gamma((aa'+nbb')\delta_{1i} + (a'b-ab')\delta_{1i}') h: h\in [\pm H]\setminus\{0\}\} $ (resp. $\nu_*$ is uniform on $\{M_g\gamma (a^2+nb^2)\delta_{1i}h: h\in [\pm H]\setminus\{0\}\}$). By our definition of $\calQ$, we note that this is never the empty set. 
Our new goal is to show that 
\begin{equation}\label{eq: type II convolution bound on * measures}
    \|\mu_*^{(2)}\|_\infty, \|\nu_*^{(2)}\|_\infty \ll \delta^{O(1)} X^{-1/2}(\log X)^{O(1)},
\end{equation}
as this will imply \eqref{eq: type I convolution bound on mu' nu'} since we can expand $\mu'^{(m)}$ and $\nu'^{(m)}$ in terms of $\delta_0$ and $\mu_*$ (resp. $\delta_0$ and $\nu_*$).

We treat the bound on $\mu_*$ first. Let $t\neq 0$. Then 
\begin{multline}
    \#\{(a,b,a',b')\in \calQ: (aa'+nbb')\delta_{1i} + (a'b-ab')\delta_{1i}' = t\} \\ \ll \#\{(a,b,a',b')\in \calQ: a(a'\delta_{1i} - b'\delta_{1i}') + b(nb'\delta_{1i} + a'\delta_{1i'})= t\}\\
    \ll \sum_{|u|\ll D} \#\left\{(a,b,a',b')\in \calQ: \begin{array}{c}
         a(a'\delta_{1i} - b'\delta_{1i}') = u \\ 
         b(nb'\delta_{1i} + a'\delta_{1i'}) = t-u 
    \end{array}\right\}\\
    \ll \sum_{|u|\ll D} \#\left\{|a'|,|b'|\sim D^{1/2}: \begin{array}{c}
         0\neq a'\delta_{1i} - b'\delta_{1i}'\mid  u \\ 
         0\neq nb'\delta_{1i} + a'\delta_{1i'}\mid t-u 
    \end{array}\right\}
\end{multline}
Let $d_1\mid u$ be a divisor of $u$ and $d_2\mid t-u$ be a divisor of $t-u$. Since 
\begin{equation}
    \det\begin{pmatrix}
        \delta_{1i} & -\delta_{1i}' \\ 
        \delta_{1i}' & n\delta_{1i}
    \end{pmatrix} = \delta_{1i}^2 + n \delta_{1i}'^2 >0
\end{equation}
(as $\delta_{1i}\neq 0$), we know that there is at most one solution to the system of equations
\begin{equation}
   0\neq a'\delta_{1i} - b'\delta_{1i}'= d_1,\quad  
         0\neq nb'\delta_{1i} + a'\delta_{1i'}= d_2.
\end{equation}
Hence, we have that 
\begin{equation}
    \#\{(a,b,a',b')\in \calQ: (aa'+nbb')\delta_{1i} + (a'b-ab')\delta_{1i}' = t\}  \ll \sum_{|u|\ll D} \tau(u)\tau(t-u) \ll D (\log D)^3.
\end{equation}
Since this bound holds uniformly for $|t|\ll D$, we know that 
\begin{equation}
    \|\mu_*\|_\infty \ll |\calQ|^{-1} H^{-1} D(\log D)^3 \tau. 
\end{equation}
As $H = \lfloor X^{1/2}/D\rfloor$ and $|\calQ|\gg \delta^{O(1)}D^2$, we see that $\|\mu_*\|_\infty\ll \delta^{O(1)} X^{-1/2}(\log X)^3 \tau$.
Thus, 
\begin{equation}
    \mu_*^{(2)}(t) \ll \delta^{O(1)}X^{-1} (\log X)^6 \sum_{|x|\ll X^{1/2}}\tau(x)\tau(t-x) \ll \delta^{O(1)} X^{-1/2}(\log X)^9.
\end{equation} 
This completes the proof of the first half of \eqref{eq: type II convolution bound on * measures}. \\ 

We apply a similar procedure to understand $\nu_*$. Let $t\neq 0$. Then 
\begin{equation}
    \#\{(a,b,a',b')\in \calQ: (a^2+nb^2)\delta_{1i} = t\}
    \ll D \#\{|a|,|b|\sim D^{1/2}: a^2+nb^2 = \delta_{1i}^{-1} t\} \ll D \tau(t).
\end{equation}
Here we have used that for any integer $N$, $\#\{(a,b): a^2+nb^2 = N\}  \ll \tau(N)$. Thus, we have the upper bound that 
\begin{equation}
    \|\nu_*\|_\infty \ll |\calQ|^{-1} H^{-1} D \tau^2\ll \delta^{O(1)} X^{-1/2} \tau^2.
\end{equation}
Finally, we see that 
\begin{equation}
    \nu_*^{(2)}(t) \ll \delta^{O(1)} X^{-1} \sum_{|x|\ll X^{1/2}} \tau^2(x)\tau^2(t-x) \ll \delta^{O(1)} X^{-1/2} (\log X)^{15}.
\end{equation}
Thus, this completes the proof of the second half of \eqref{eq: type II convolution bound on * measures} and hence the lemma. 
\end{proof}

\part{Weak approximation and the Brauer--Manin obstruction}\label{part:bmo}
\stepcounter{section}
\setcounter{thm}{0}
\setcounter{equation}{0}

The goal of this section is to establish Theorem \ref{thma:BMO}. As mentioned in the introduction, there is a longstanding history of applying new instances of prime number theorems to prove Hasse principles and weak approximation for algebraic varieties. Indeed, Hasse used Dirichlet's theorem for primes in arithmetic progressions to establish the Hasse principle for quadratic forms for three and four variables (see \cite[p. 16 and p. 87]{Hasse:Hasse_principle}). Colliot-Th\'el\`ene and Sansuc \cite{CT-Sansuc:schinzel-hasse} then used Schinzel's Hypothesis H to prove the Hasse principle and weak approximation for conic bundles of the shape $$y^2-az^2 = P(x)$$ for $a\in \Q^\times$ and $P(x)\in \Q[x]$ an irreducible polynomial. Using descent theory, they additionally establish that when $P(x)$ is reducible and separable, and not of even degree with an odd degree irreducible factor, then the Brauer-Manin obstruction is the only obstruction to the Hasse principle and weak approximation. Since then this result has been vastly generalized. In \cite{CT-SD:schinzel-hasse-severi-brauer}, Colliot-Th\'el\`ene and Swinnerton-Dyer showed that under Schinzel's Hypothesis H, if $X$ is a smooth and projective variety that admits a fibration $X/\mathbb{P}_{\Q}^1$ where the generic fibre is a generalized Severi-Brauer variety (see Definition \ref{def:generalized_Severi_Brauer}), then the Brauer-Manin obstruction is the only obstruction to weak approximation. The main result of \cite{CT-SD:schinzel-hasse-severi-brauer} completes earlier results of \cite{CT-Sansuc:schinzel-hasse,Serre:Galois_cohomology,SD:pencils-of-conics} which showed that Schinzel's Hypothesis H implies that the Brauer--Manin obstruction is the only obstruction to weak approximation on conic bundles, and quadric bundles of relative dimension $2$, over $\Q$. Prior to \cite{CT-SD:schinzel-hasse-severi-brauer}, Serre generalized \cite{CT-Sansuc:schinzel-hasse} to handle arbitrary Severi--Brauer bundles (including conic bundles, but not including quadric bundles of relative dimension $2$); see \cite[Appendix to Chapter II, Theorem 7.6]{Serre:Galois_cohomology}. The generalized Severi--Brauer bundles studied in \cite{CT-SD:schinzel-hasse-severi-brauer} subsume all the bundles considered in \cite{CT-Sansuc:schinzel-hasse,Serre:Galois_cohomology,SD:pencils-of-conics}. 

\begin{defn}[See {\cite[Definition 2.1.1]{CT-SD:schinzel-hasse-severi-brauer}}]\label{def:generalized_Severi_Brauer}
    A \define{Severi-Brauer variety} $X$ over $\Q$ is a $\Q$-scheme such that $X_{\overline{\Q}}$ is isomorphic to $\P^n_{\overline{\Q}}$ for some $n$. A \define{generalized Severi-Brauer variety} $X$ over $\Q$ is a proper smooth scheme such that there exists a finite set $I$, finite (separable) field extensions $K_i/\Q$ and Severi-Brauer varieties $V_i/\Q$ such that $X$ is isomorphic to the product of Weil restrictions $\prod_{i\in I} R_{K_i/\Q}(V_i\times_\Q K_i).$ 
\end{defn}
\begin{example}
    The Severi-Brauer curves are precisely conics in $\P^2$. A quadric surface in $\P^3$ is not a Severi-Brauer variety, but is a generalized Severi-Brauer variety; see \cite[Section 2]{CT-SD:schinzel-hasse-severi-brauer}.
\end{example}

Part of the innovation of Serre's work was the introduction of a more general, number field version of Schinzel's Hypothesis H, dubbed Hypothesis (H$_1$), and the observation that it is implied by the Schinzel's Hypothesis H; this passage from H to (H$_1$) first appeared in print in \cite{CT-SD:schinzel-hasse-severi-brauer}. More recently, one often replaces Hypothesis (H$_1$) with a homogeneous variant, called Hypothesis (HH$_1$) in \cite[Section 9.2.1]{HW:fibration-method}, which is still implied by Schinzel's Hypothesis H. The passage to a homogeneous version of Schinzel's Hypothesis H goes back to earlier work of Swinnerton-Dyer; see \cite[Section 2]{SD:Applications-of-Schinzel}. Its utility lies in its likelihood of being made unconditional through analytic methods, such as additive combinatorics.

The first major success followed breakthrough work of Green and Tao \cite{Green_Tao:Linear_equations}, and Green, Tao, and Ziegler \cite{GTZ:inverse_theorem} on simultaneous prime values of linear systems.
In \cite{browning-matthiesen-skorobogatov}, Browning, Matthiesen, and Skorobogatov used the descent method and Matthiesen's work on representation of linear polynomials by binary quadratic forms \cite{Matthiesen:correlation_pos_def_quad,Matthiesen:correlation_quad} to prove the Hasse principle and weak approximation are controlled by the Brauer-Manin obstruction for pencils of conics and higher-dimensional quadrics, whose degenerate fibres are defined over $\Q$. This was later generalized further by Browning and Matthiesen in \cite{BM:norm-forms} to varieties defined by $N_{K/\Q}(x_1,\dots,x_d) = L_1\dots L_k(z)$ for any number field $K/\Q$, i.e. to the question of representing products of linear polynomials as a norm from $K$. In Harpaz, Skorobogatov, and Wittenberg \cite{HSW:RationalPoints}, a similar approach 
was used to study varieties whose fibres are generalized Severi-Brauer varieties. 
These results (whose general strategy is labeled ``the fibration method") have since been further streamlined and extended in works of Harpaz and Wittenberg \cite{HW:fibration-method} and Harpaz, Wei, and Wittenberg \cite{HWW:fibration}. In particular, in \cite[Theorem 1.6]{HW:fibration-method}, Harpaz and Wittenberg apply Matthiesen's work on norm forms \cite{Matthiesen:norm} to unconditionally prove that the Brauer-Manin obstruction is the only obstruction to weak approximation on any fibration $X\to\P^1_\Q$ whose generic fiber is (geometrically) rationally connected, provided, for example, its \emph{nonsplit}\footnote{A variety $X$ is \define{split} if it contains a geometrically integral open subscheme $U\subset X$.} fibers lie over $\Q$-points of $\P^1$ and that the Brauer-Manin obstruction is the only obstruction to weak approximation on all fibers $X_c$ for all $c$ outside a finite subset of $\P^1_\Q$.

Our new unconditional prime number theorem (Theorem \ref{thm: prime count}) provides new unconditional cases of Hypothesis (HH$_1$) (see Lemma \ref{lem:BMO hypothesis H-1 sub}), and so results of \cite{HW:fibration-method,HWW:fibration} convert it into new unconditional instances of the fibration method (the rough expectation that if the Brauer-Manin obstruction controls the arithmetic of the fibers of a map $X\to\P^1$, then the same should be true on the total space $X$). Below, we record one such statement, Theorem \ref{thm:general-BMO}, of which Theorem \ref{thma:BMO} is a special case. In Appendix \ref{app:BMO}, we will give a simple, direct proof of Theorem \ref{thma:BMO} as well.

First, recall that if $\pi\colon X\to\P^1$ is a dominant morphism of smooth, integral $\Q$-varieties, then the corresponding \define{vertical Brauer group} is defined to be
\begin{equation}
    \vBr(X)\coloneqq\Br(X)\cap\pull\pi\Br(\Q(\P^1))\subset\Br(\Q(X)).  
\end{equation}
\begin{thm}\label{thm:general-BMO}
    Let $X$ be a smooth, proper, geometrically integral $\Q$-variety equipped with a dominant morphism $\pi\colon X\to\P^1$.
    Assume
    \begin{alphabetize}
        \item\label{item:BMO-bad fibers} there is an imaginary quadratic field $K$ and a finite list of fibers $X_0,X_1,\dots,X_r$ of $\pi$, with $X_0$ defined over $K$ and all other $X_i$'s defined over $\Q$, such that \important{every other} fiber of $\pi$ contains a geometrically integral irreducible component of multiplicity one;
        \item\label{item:BMO-abelian extensions} for each $i=0,1,\dots,r$, the fiber $X_i$ contains an irreducible component of multiplicity one such that the algebraic closure of $\Q$ in its field of functions is an abelian extension of $\Q$, if $r\ge1$, and is an abelian extension of $K$ if $r=0$.
        \item\label{item:BMO-fibral Hasse/WA} all but finitely many $\Q$-fibers of $\pi$ satisfy the Hasse principle (resp. weak approximation).
    \end{alphabetize}
    
    Then, the vertical Brauer--Manin obstruction is the only obstruction to the Hasse principle (resp. weak approximation) on $X$.
\end{thm}
\begin{rem}
    The vertical Brauer--Manin obstruction being the only obstruction to the Hasse principle means that the implication $X(\A_\Q)^{\vBr(X)}\neq\emptyset\implies X(\Q)\neq\emptyset$ holds for such $X$. It being the only obstruction to weak approximation is the statement that $X(\A_\Q)^{\vBr(X)}$ is the topological closure of $X(\Q)\subset X(\A_\Q)$.
\end{rem}
\begin{cor}
    Let $X/\Q$ be a smooth projective variety which admits a fibration $X\to\mathbb{P}_\Q^1$ whose generic fibre is a generalized Severi-Brauer variety. Assume that all of the degenerate fibres are defined over $\Q$, except one that may be defined over an imaginary quadratic field $K/\Q$. Then the Brauer-Manin obstruction is the only obstruction to weak approximation. 
\end{cor}
\begin{example}[Theorem \ref{thma:BMO}]
    Let $P(x)\in \Q[x]$ be a separable polynomial that factors as 
    \begin{equation}
        P(x) = L_1(x)\dots L_k(x) \cdot Q(x),
    \end{equation}
    where $L_1,\dots,L_k$ are linear and $Q(x)$ is a irreducible quadratic satisfying that $K = \Q[x]/Q(x)$ is an imaginary quadratic field. Then Theorem \ref{thm:general-BMO} implies that the Brauer-Manin obstruction controls weak approximation for the conic bundle over $\Q$ given by the affine equation 
    \begin{equation}
        y^2-az^2 = P(x)
    \end{equation}
    for $a\in \Q^\times.$
\end{example}

Let $L_1(t), \dots, L_r(t)\in \Q[t]$ be linear polynomials satisfying that $L_i(t)\neq \lambda L_j(t)$ for $i\neq j$ and let $Q(t)\in \Q[t]$ be an irreducible quadratic polynomial. Define  $S_{L_1,\dots,L_r,Q}$ to be the set of finite primes $p$ satisfying any of the following:
\begin{itemize}
    \item $p$ divides the denominator of a rational coefficient in $L_1(t),\dots,L_r(t)$ or $Q(t)$;
    \item $p\mid\Disc Q(t)$; or
    \item $\b{t\in\F_p\colon L_i(t),Q(t)\in\units\F_p\tcomma\forall i\in[1,r]}=\emptyset$.
\end{itemize}
In other words, $S_{L_1,\dots,L_r,Q}$ contains the set of primes for which there may be local obstructions to $L_i(t)$ and $Q(t)$ being simultaneously prime. 

\begin{lemma}\label{lem:BMO hypothesis H-1 sub}
Let $L_1(t),\dots L_r(t)\in \Q[t]$ be linear polynomials satisfying that $L_i(t)\neq \lambda L_j(t)$ for $i\neq j$, {let $K$ be an imaginary quadratic field, and set $Q(t) = N_{K/\Q}(t-e_0)$ for some $e_0\in K$}.
    Let $S$ be a finite set of primes containing $S_{L_1,\dots,L_r,Q}\cup\{\infty\}$. Given elements $\lambda_p\in \Q_p$ for $p\in S$, a number $\epsilon>0$, and integers $n_p\geq 0$ for $p\in S$, and abelian extensions $F_i/\Q$ for $i>0$ and $F_0/K$, there exists an element $\lambda\in \Q$ and primes $\{p_1,\dots,p_r,q\}\not\in S$ such that 
    \begin{enumerate}
        \item\label{item:H1 sub-close finite} $\val_p(\lambda-\lambda_p) \geq n_p$ for $p\in S\setminus \{\infty\}$,
        \item\label{item:H1 sub-close infinite} $|\lambda-\lambda_\infty|<\epsilon,$
        \item\label{item:H1 sub-S-unit L} for each $i=1,\dots,r$, for any prime $p\not\in S\cup \{p_i\}$, $\val_p(L_i(\lambda))\le 0$,
        \item\label{item:H1 sub-S-unit Q} for any prime $p\not\in S\cup \{q\}$, $\val_p(Q(\lambda))\le 0,$
        \item\label{item:H1 sub-uniformizers} $\val_{p_i}(L_i(\lambda))=1$ for all $i=1,\dots,r$, and $\val_q(Q(\lambda))=1$.
        \item\label{item:H1 sub-hilbert Q} Fix some $i\in\{1,\dots,r\}$, and a cyclic subextension $E/\Q$ of $F_i/\Q$, which is unramified outside $S$. Assume
        \[\sum_{p\in S}\opinv_p\p{E/\Q,L_i(\lambda_p)}=0\in\Q/\Z.\]
        Then, $\opinv_{p_i}(E/\Q,L_i(\lambda))=0$.
        \item\label{item:H1 sub-hilbert K} Fix a cyclic subextension $E/K$ of $F_0/K$, which is unramified outside $S$. Assume
        \[\sum_{p\in S}\opinv_p\opname{Cores}_{K/\Q}\p{E/K,\lambda_p-e_0}=0\in\Q/\Z.\]
        Then, $\opinv_q\opname{Cores}_{K/\Q}(E/K,\lambda-e_0)=0$.
    \end{enumerate}
\end{lemma}
\begin{proof}
    Set $D\coloneqq\max_{i\ge0}[F_i:\Q]$. For each $p\in S\setminus\{\infty\}$, fix $D_p\ge1$ large enough that $\val_p(x-y)\ge D_p$ implies that $x/y\in(\units\Q_p)^D$, and then set $N_p\coloneqq\max(D_p,n_p,n_p-\val_p(\lambda_p))$.
    For $p\in S\setminus\{\infty\}$, by the Chinese remainder theorem,  we can find $\alpha,\beta\in \Z_S$ such that 
    $\val_p(\alpha-\lambda_p) \geq N_p$ and $\val_p(\beta-1)\geq N_p$ for $p\in S\setminus\{\infty\}.$ We can also assume that $|\alpha/\beta - \lambda_\infty|\leq \epsilon.$

    Let $M= \prod_{p\in S} p^{N_p}.$ Assume that 
    \begin{equation}
        L_i(t) = a_i t+ b_i, \quad Q(t)=N_{K/\Q}(t-e_0) = a_0 t^2 + b_0 t + c_0.
    \end{equation}
    Let us define: 
    \begin{equation}
        L_i'(x,y) = d_i(a_i(Mx+\alpha) + b_i(My+\beta)), \quad 
        Q'(x,y) = d_0 (a_0(Mx+\alpha)^2 + b_0(Mx+\alpha)(My+\beta) + c_0(My+\beta)^2).
    \end{equation}
    Here we take $d_0,d_1,\dots,d_r\in \Q_{\geq 0}$ to satisfy that the coefficients of $L_i'(x,y)$ and $Q'(x,y)$ are integral and coprime; by definition of {$S_{L_1,\dots,L_r,Q}\subset S$}, such $d_i$ will be $S$-units. Recall that $Q(t) = N_{K/\Q}(t-e_0)$ and define $e_i = b_i/a_i$ for each $i>0$.
    We can further also choose $\alpha,\beta$ and $d_0,\dots,d_r$ such that for each $p\in S\setminus\{\infty\}$, 
    \begin{equation}\label{eq: hypothesis h_1 choosing moduli}
        \val_p(d_i(\alpha+e_i\beta))=0, \quad \val_p(d_0 N_{K/\Q}(\alpha-e_0\beta))=0. 
    \end{equation}
    (for example, at $p=3$ for $L_1=t+1,L_2=t+2,L_3=t+3$, take $n_3\geq 1$, $(\alpha,\beta)=(1,1)$ and $d_1=1$, $d_2 = 1/3$, $d_3=1$).
    
    To establish the lemma, it suffices to find $(x,y)\in\Z$ such that $|\frac{Mx+\alpha}{My+\beta}-\lambda_\infty|<\epsilon$, and such that $L_i'(x,y)=p_i$ and $Q'(x,y)=q$ for distinct primes $p_1,\dots,p_r,q\not\in S$. We then set $\lambda=\frac{Mx+\alpha}{My+\beta}$, so $L_i(\lambda)=\frac1{d_i(My+\beta)}L_i'(x,y)$ and $Q(\lambda)=\frac1{d_0(My+\beta)^2}Q'(x,y)$. It is clear that such a $\lambda$ satisfies \ref{item:H1 sub-close finite} -- \ref{item:H1 sub-uniformizers}. To see \ref{item:H1 sub-S-unit L} and \ref{item:H1 sub-S-unit Q}, it is enough to note that if $\val_{p_i}(My+\beta)>0$ then $\val_{p_i}(Mx+\alpha)>0$ and comparing when $i\neq j$ this would contradict $\{p_1,\dots,p_r,q\}$ being distinct primes. So, it remains to check properties \ref{item:H1 sub-hilbert Q} and \ref{item:H1 sub-hilbert K}.
    
    Suppose $i>0$ and suppose a cyclic subextension $E/\Q$ of $F_i/\Q$, which is unramified outside $S$, is given so that
    \[\sum_{p\in S}\opinv_p(E/\Q,L_i(\lambda_p))=0\in\Q/\Z.\]
    By construction, for any $p\in S$, we have chosen $N_p$ such that $\lambda$ is sufficiently close to $\lambda_p$ that $L_i(\lambda)/L_i(\lambda_p)$ is an $[E:\Q]$th power in $\Q_p$, so $\opinv_p(E/\Q,L_i(\lambda_p))=\opinv_p(E/\Q,L_i(\lambda))$. Thus, 
    \[\sum_{p\in S}\opinv_p(E/\Q,L_i(\lambda))=0.\]
    Next, observe that, for $p\in S$, $\val_p(My+\beta-1)\ge N_p$, so $(My+\beta)$ is a $D$th power in $\Q_p$. Hence, 
    $\opinv_p(E/\Q,L_i(\lambda))= \opinv_p(E/\Q, d_i^{-1}L_i'(x,y))$.
    Therefore, the following relation holds: 
    \[0 = \sum_{p\in S}\opinv_p(E/\Q,d_i^{-1}L_i'(x,y))=\sum_{p\not\in S}\opinv_p(E/\Q,d_i^{-1}L_i'(x,y)),\]
    with the latter equality holding by Hilbert reciprocity. Consider now those primes $p\not\in S$: for $p\neq p_i$ we have that $\val_{p}(d_i^{-1}L_i'(x,y))=0$ and hence $\opinv_p(E/\Q, d_i^{-1}L_i'(x,y))=0$. Thus, this implies that 
    \begin{equation}\opinv_{p_i}(E/\Q,d_i^{-1} L_i'(x,y)) = \opinv_{p_i}(E/\Q, L_i(\lambda)) =0\end{equation}

    Finally, suppose $i=0$ and $E/K$ is a cyclic subextension of $F_0/\Q$ which is unramified outside of $S$. Assume that 
    \[\sum_{p\in S}\opinv_p\opname{Cores}_{K/\Q}(E/K,\lambda-e_0)=0\in\Q/\Z.\]
    Indeed, by the choice of $N_p$, $(\lambda-e_0)/((Mx+\alpha) - e_0 (My+\beta))$ is a $D$th power in $K_p=\Q_p\otimes_\Q K$, so
    \[\opinv_p\opname{Cores}_{K/\Q}(E/K,\lambda-e_0)=\opinv_p\opname{Cores}_{K/\Q}\big(E/K,(Mx+\alpha)-e_0(My+\beta)\big).\] 
    Hence,
    \begin{multline}0 = \sum_{p\in S}\opinv_p\opname{Cores}_{K/\Q}(E/K,(Mx+\alpha) - e_0 (My+\beta))\\ =\sum_{p\not\in S}\opinv_p\opname{Cores}_{K/\Q}(E/K,(Mx+\alpha) - e_0 (My+\beta)).\end{multline}
    Finally, since $Q'(x,y)=q$ is prime, we know that $N_{K/\Q}((Mx+\alpha) - e_0 (My+\beta))=q d_0^{-1}$ where $\val_{p}(d_0)=0$ for $p\not\in S$. Hence, we have that 
    \begin{equation}
        \opinv_q \opname{Cores}_{K/\Q}(E/K, (Mx+\alpha) - e_0 (My+\beta)) = \opinv_q \opname{Cores}_{K/\Q}(E/K, \lambda - e_0) = 0.
    \end{equation}

    To establish the existence of such a $\lambda$, we want to apply Theorem \ref{thm: prime count} in the region $\{\left|\frac{Mx+\alpha}{My+\beta}-\lambda_\infty\right|<\epsilon\}$ which is a fixed convex region (and it is straightforward to see that the proof of Theorem \ref{thm: prime count} generalizes to a convex region). By definition of $S_{L_1,\dots,L_r,Q}\subset S$, we can see that  
    \begin{equation}
        \prod_{p\not\in S}\frac{p^{-2}\#\{x,y\in \F_p: L_i'(x,y),Q'(x,y)\in \F_p^\times, \forall i\in[1,r]\}}{1-p^{-1}}>0.
    \end{equation}
    At those primes $p\in S$, we have that $L_i'(x,y) \equiv d_i(a_i\alpha+b_i \beta) \bmod p$
    and $Q'(x,y) \equiv d_0(a_0\alpha^2+b_0\alpha\beta+c_0\beta^2)\bmod p$. So, by our assumption \eqref{eq: hypothesis h_1 choosing moduli}, we also have that 
    \begin{equation}
        \prod_{p\in S}\frac{p^{-2}\#\{x,y\in \F_p: L_i'(x,y),Q'(x,y)\in \F_p^\times, \forall i\in[1,r]\}}{1-p^{-1}}>0.
    \end{equation}
    Hence, we know that $\mathfrak{S}_{L_1',\dots,L_r',Q'}>0$ and by Theorem \ref{thm: prime count} we can find such a tuple $(x,y)\in \Z$. 
\end{proof}
\begin{rem}
    Let $L_1(\lambda,\mu),\dots,L_r(\lambda,\mu)$ be homogeneous linear polynomials, and let $Q(\lambda,\mu)$ be an irreducible homogeneous quadratic polynomial with negative discriminant. It is straightforward to deduce \cite[Hypothesis (HH$_1$)]{HW:fibration-method} for $L_1,\dots,L_r,Q$ from the proof of Lemma \ref{lem:BMO hypothesis H-1 sub} (indeed, the proof moves through the two-variable polynomials $L_i'$ and $Q'$); observe that we only need \ref{item:H1 sub-close finite}-\ref{item:H1 sub-uniformizers} and to achieve the final positivity condition (in \cite{HW:fibration-method} denoted as $\lambda_0\lambda_\infty + \mu_0\mu_\infty >0$ ) we note that Theorem \ref{thm: prime count} can be applied in this convex region.
    We have written Lemma \ref{lem:BMO hypothesis H-1 sub} in the above form because this is the form we found convenient for applying in Appendix \ref{app:BMO}.
\end{rem}

The following proof was communicated to us by Olivier Wittenberg.
\begin{proof}[Proof of Theorem \ref{thm:general-BMO}]
    In light of Lemma \ref{lem:BMO hypothesis H-1 sub}, Theorem \ref{thm:general-BMO} is ultimately a consequence of \cite[Theorem 4.1]{HWW:fibration}. In the notation of loc. cit., take $A=0$ to be the trivial subgroup and take $H=\P^1_\Q$. For $m\in M=\{\pi(X_0),\dots,\pi(X_r)\}\subset\P^1_\Q$, let $Y_m\subset X_m=\inv\pi(m)$ be the multiplicity one irreducible component guaranteed by \ref{item:BMO-abelian extensions}, let $E_m=\Q(Y_m)$ be its field of rational functions, and set $L_m=K_m$ to be algebraic closure of $\kappa(m)$ in $\Q(Y_m)$, so $L_m/\kappa(m)$ is abelian by \ref{item:BMO-abelian extensions}. Then, the hypothesis in \cite[Theorem 4.1]{HWW:fibration} that \cite[Conjecture F$_+$]{HWW:fibration} holds for the relevant parameters is verified by \cite[Corollary 6.7]{HWW:fibration} in light of Lemma \ref{lem:BMO hypothesis H-1 sub} proving Hypothesis (HH$_1$) for the relevant polynomials. Thus, \cite[Theorem 4.1]{HWW:fibration} yields that
    \[\bigcup_{c\in\P^1(\Q)-M(\Q)}X_c(\A_\Q)\,\t{ is dense in }\,X(\A_\Q)^{\vBr(X)},\]
    from which the claim follows.
\end{proof}

We end by remarking that Theorem \ref{thm:general-BMO} is not the only or most general possible instance of the fibration method made unconditional by Theorem \ref{thm: prime count}/Lemma \ref{lem:BMO hypothesis H-1 sub}. Other variations of Theorem \ref{thm:general-BMO} may be obtained, for example, by combining Lemma \ref{lem:BMO hypothesis H-1 sub} with \cite[Theorem 1.4]{HWW:fibration} or \cite[Corollaries 4.6 + 4.7]{HWW:fibration}. For example, a consequence of these methods is that the Brauer-Manin obstruction is the only obstruction to the Hasse principle and weak approximation for the varieties cut out by $N_{K/\Q}(x_1,\dots,x_d)=L_1\dots L_k(z) Q(z)$, where $K/\Q$ is a degree $d$ cyclic extension, and $Q(z)$ is a quadratic polynomial that defines an imaginary quadratic field. However, these results do not extend to all number fields $K/\Q$ without extra work.

\part{Elliptic curves of rank two}\label{part:rank2}

\section{Main statement \& a Rubin--Silverberg quadratic twist family}

As our second application of our main Theorem \ref{thm: prime count}, we wish to use it to construct new examples of elliptic curves, over $\Q$, of rank exactly two. Recently, Zywina \cite{Zywina:Rank2} has shown there are infinitely elliptic curves over $\Q$, with \important{differing $j$-invariants}, who each have rank exactly $2$. We wish to complement this result by showing that, for fixed $j$-invariant (i.e. within a single $\bar\Q$-isomorphism class), one may still find infinitely many non-$\Q$-isomorphic elliptic curves of rank exactly $2$.

\begin{thm}\label{thm:rank2}
    For an admissible $(\mu,\lambda)\in\Z^2$, as defined in Definition \ref{defn:admissible}, the elliptic curve
    \[E_{\mu,\lambda}\colon y^2=x(x-\mu)(x-\lambda)\]
    has infinitely many quadratic twists of rank exactly $2$.
\end{thm}

The starting point for our proof of Theorem \ref{thm:rank2} is a construction, due to Rubin and Silverberg  \cite{Rubin-Silverberg}, of a family of quadratic twists of $E_{\mu,\lambda}$ with generic rank $2$ (over $\Q(u)$):

\begin{prop}\label{prop:RS-family}
    Let $\mu,\lambda\in\Z$ be distinct nonzero integers. Define
    \[\alpha\coloneqq\frac{\mu\lambda}{\mu+\lambda}\in\Q\tand G(u)\coloneqq\frac{(u^2+\alpha^2)(u^2+\alpha^2-\alpha\mu)(u^2+\alpha^2-\alpha\lambda)}{\alpha^3}\in\Q[u].\]
    Then, the quadratic twist family
    \[E_{\mu,\lambda}^{G(u)}\colon G(u)y^2=x(x-\mu)(x-\lambda)\]
    has generic rank $\ge2$, i.e. it has rank $\ge2$ over $\Q(u)$.

    Furthermore, if $\mu\lambda$ is a square, say $\mu\lambda=\eta^2$, then $G$ further factors as
    \begin{equation}\label{eqn: G(u)=QLLLL}
        G(u) = \frac1{\alpha^3}(u^2+\alpha^2)\p{u+\frac{\eta\mu}{\mu+\lambda}}\p{u-\frac{\eta\mu}{\mu+\lambda}}\p{u+\frac{\eta\lambda}{\mu+\lambda}}\p{u-\frac{\eta\lambda}{\mu+\lambda}}
        .
    \end{equation}
\end{prop}
\begin{proof}
    One can check that $\mcE\coloneqq E_{\mu,\lambda}^{G(u)}$ admits the points
    \[P_1\coloneqq\p{\frac{u^2+\alpha^2}\alpha,1}\tand P_2\coloneqq\p{\frac{\alpha u^2+\alpha^3}{u^2},\frac{\alpha^3}{u^3}}.\]
    We claim that these points are $\Z$-linearly independent in $\mcE(\Q(u))$, which suffices to prove the proposition. For this, note that $\mcE$ is defined already over $\Q(u^2)\subset\Q(u)$, and so inherits an action of $\Gal(\Q(u)/\Q(u^2))\simeq\zmod2$. This action is trivial on $P_1$ but is negation on $P_2$, from which one sees that the two points are independent.
\end{proof}
\begin{rem}
    The quadratic twist family in Proposition \ref{prop:RS-family} comes from work of Rubin and Silverberg \cite{Rubin-Silverberg}. It was obtained by applying \cite[Propositions 2.7 + 2.9]{Rubin-Silverberg}, starting with the linear fractional transformation which fixes $0$ and swaps $\mu,\lambda$. \cite[Section 4]{Rubin-Silverberg} gives a method for finding quadratic twist families, of the form $E^{G(u)}$ for a suitable starting curve $E/\Q$ and polynomial $G(u)$, of generic rank $\ge3$. It would interesting to search for examples where known cases of Conjecture \ref{conj: multi BH} may be applied to the irreducible factors of $G(u)$ and then use this to construct rank exactly $3$ elliptic curves in quadratic twist families.
\end{rem}
We will ultimately prove Theorem \ref{thm:rank2} by showing that, under suitable hypotheses, there are infinitely many values of $u\in\Q$ such that $\rank E^{G(u)}_{\mu,\lambda}(\Q)=2$. The lower bound $\rank E^{G(u)}_{\mu,\lambda}(\Q)\ge2$ will follow from Proposition \ref{prop:RS-family}. The upper bound will rely on a $2$-Selmer computation. As in previous applications of additive combinatorics towards ranks of elliptic curves -- such as \cite{koymans2025ellipticcurvesranknumber,Zywina:Rank1_NumberFields,Zywina:SmallRanks,Savoie:rank2} -- the main point will be to find values of $u$ such that the factors of $G(u)$ (one quadratic and four linear; see (\ref{eqn: G(u)=QLLLL})) are all primes satisfying congruence conditions suitable for controlling the 2-descent. This will be taken care of by an application of Theorem \ref{thm: prime count}. Before this though, one needs to do a certain amount of `pre-twisting'; in essence, one needs to know an explicit basis of $\Sel_2(E)$ in order to adequately control how twisting by the primes produced by Theorem \ref{thm: prime count} (i.e. the factors of $G(u)$) will affect 2-Selmer. This `pre-twisting' is the subject of Section \ref{sect:pre-twist}, after which we apply Theorem \ref{thm: prime count} to prove Theorem \ref{thm:rank2} in Section \ref{sect:AC-twist}. The strategy of combining descent in quadratic twist families and prime number theorems like Theorem \ref{thm: prime count} is due to Koymans and Pagano, in their work on Hilbert's 10th problem \cite{Koymans-Pagano:H10}.

We end this section by defining the admissible condition referenced in Theorem \ref{thm:rank2}.
\begin{defn}\label{defn:admissible}
    A pair $(\mu,\lambda)\in\Z^2$ is \define{admissible} if it satisfies all of the following:
    \begin{enumerate}[label=\bfseries{(A$_{\arabic*}$)}]
        \item\label{item:admissible-square product} $\mu\lambda$ is a square, but $\mu$ is squarefree.

        Hence, $\mu$ is precisely the squarefree part of $\lambda$, i.e. $\lambda=\mu m^2$ for some $m\in\Z_{>0}$.
        \item\label{item:admissible-m-at-2} $m\equiv2\bmod8$, with $m$ as above, but $\mu$ is odd.

        In light of \ref{item:admissible-square product}, this forces $\lambda\equiv4\bmod8$.
        \item\label{item:admissible-p primes} there exists distinct odd primes $p_1,p_2>14$ each satisfying
        \begin{equation}
            p_i\nmid\mu\lambda,\tbut\ord_{p_i}(\mu-\lambda)\t{ is odd}.
        \end{equation}
        For each $p_i$, write $m^2-1=p_i^{\ord_{p_i}(1-m^2)}\cdot d_i=p_i^{\ord_{p_i}(\mu-\lambda)}\cdot d_i$; we further require that
        \[d_i\not\equiv\square\bmod p_i,\quad\t{i.e. that }(1-m^2,p_i)_{p_i}=(-p_i,p_i)_{p_i}(d_i,p_i)_{p_i}=-1,\]
        where $(-,-)_p$ denotes the multiplicative Hilbert symbol.
        Finally, we require that $p_1\equiv\pm1\bmod 8$ (i.e. $2\equiv\square\bmod p_1$) and $p_2\equiv7\bmod8$ (i.e. $2\equiv\square\bmod p_2$ and $p_2\equiv3\bmod4$).

        \item\label{item:admissible-l primes} there exists distinct odd primes $\l_0,\l_1$ each satisfying
        \[\l_i\nmid\mu\tbut\l_i\mid\lambda.\]
        We further require that $\l_1\equiv\pm1\bmod8$ (i.e. $2\equiv\square\bmod\l_1$).
        In light of \ref{item:admissible-square product}, $\ord_{\l_i}(\lambda)$ is even, for each $i$.
        \item\label{item:admissible-coprimality} $\gcd(\mu,m^5-m)=\gcd(\mu, m)\gcd(\mu,m^2+1)\gcd(\mu,m^2-1)=1$.
        \item\label{item:admissible-odd-num-factors} finally, we require that $m^2-1$ has an odd number of prime factors.
    \end{enumerate}
\end{defn}

One may construct infinitely many admissible pairs $(\mu,\lambda)$. To do so, we first need to establish the following lemma.
\begin{lemma}\label{lem:odd-prime-factors}
    Let $a,q\in \mathbb{N}$. Let $\mu\colon\N\to\Z$ denote the M\"obius function, so $\mu(n)=(-1)^k$ if $n$ is squarefree with $k$ prime factors, and $\mu(n)=0$ otherwise. Then there are infinitely many integers $n \equiv a \bmod q$ such that 
    \begin{equation}
        \mu(n)\mu(n+2)=-1. 
    \end{equation}
    In other words, there are infinitely many $n\equiv a \bmod q$ such that $n,n+2$ are both squarefree and the number of prime factors of $n$ and $n+2$ have opposite parity. 
\end{lemma}
\begin{proof}
    The above is a consequence of \cite[Corollary 5]{Tao:log-Chowla} of Tao, resolving the logarithmically averaged Elliott's conjecture. Indeed, taking $a_1=a_2 = q$ and $b_1=a$, $b_2=a+2$, we have that $a_1b_2-a_2b_1 \neq 0$ and thus: 
    \begin{equation}
        \sum_{n'\leq x} \frac{\mu(qn'+a)\mu(qn'+a+2)}{n'} = o_{q,a}(\log x).
    \end{equation}
    Hence, $\mu(qn'+a)\mu(qn'+a+2) = -1$ infinitely often and we have the desired statement. 
\end{proof}
\begin{construct}
    Choose distinct odd primes $\l_0,\l_1,p_1,p_2$ with $p_1,p_2>14$ for which $2\equiv\square\bmod \l_1,p_1,p_2$ and $p_2\equiv7\bmod8$. Fix some $\alpha_1,\alpha_2\in\Z$ such that $p_i\nmid\alpha_i$ and $(\alpha_ip_i,p_i)_{p_i}=-(2,p_i)_{p_i}$ (i.e. $-2\alpha_i\not\equiv\square\bmod p_i$) for each $i$. By Lemma \ref{lem:odd-prime-factors}, we may choose some $n\in\Z$ such that $m\coloneqq n+1$ satisfies
    \begin{equation}\label{eqn:m-CRT-admissible-construction}
        m\equiv2\bmod8\tcomma m\equiv\l_0\l_1\bmod(\l_0\l_1)^2\tcomma m\equiv \alpha_1p_1-1\bmod p_1^2,\tand m\equiv\alpha_2p_2-1\bmod p_2^2,
    \end{equation}
    and so that $n(n+2)=(m-1)(m+1)=m^2-1$ is squarefree with an odd number of prime factors.
    Choose any odd squarefree $\mu$ which is coprime to $m(m^4-1)$. Set $\lambda=\mu m^2$. Then, $\mu\lambda=(\mu m)^2$ and $\mu-\lambda=\mu(1-m)(1+m)$, so
    \[\ord_{\l_i}(\lambda)=2\ord_{\l_i}(m)=2,\tand\ord_{p_i}(\lambda)=\ord_{p_i}(1+m)=1.\]
    Additionally, for each $i$,
    \[(1-m^2,p_i)_{p_i}=(1-m,p_i)_{p_i}(1+m,p_i)_{p_i}\overset{(\ref{eqn:m-CRT-admissible-construction})}=(2,p_i)_{p_i}(\alpha_ip_i,p_i)_{p_i}=-1,\]
    with the last equality holding by definition of $\alpha_i$. Hence, $(\mu,\lambda)$ is admissible.
\end{construct}

\begin{ex}\label{ex:17}
    Let $(\mu,\lambda)=(17,184\,009\,700)$. Then $\lambda=\mu m^2$, where $m=3290=2\cdot5\cdot7\cdot47$ is divisible by $\l_0\coloneqq5$ and $\l_1\coloneqq7$. Hence, $\mu$ is odd, squarefree, and coprime to $m$; $\mu\lambda=(\mu m)^2=(55\,930)^2$ is a square; and $m\equiv2\bmod8$. In addition,
    \[\mu-\lambda=\mu(1-m^2)=-\mu\cdot 10\,824\,099=-\mu\cdot3\cdot11\cdot13\cdot23\cdot1097\]
    is divisible by $p_1\coloneqq1097$ and $p_2\coloneqq 23$. Additionally, observe/compute that
    \begin{align}
        \l_1,p_2&\equiv7\bmod8\tcomma p_1\equiv1\bmod8\tand \\
        (1-m^2,p_1)_{p_1}&=(-10\,824\,099,1097)_{1097}=-1=(-10\,824\,099,23)_{23}=(1-m^2,p_2)_{p_2}
        .
    \end{align}
    Finally, $m^4-1=117\,161\,140\,809\,999=3\cdot11\cdot13\cdot23\cdot1097\cdot10\,824\,101$ is coprime to $\mu=17$, and $m^2-1=3\cdot11\cdot13\cdot23\cdot1097$ has an odd number of prime factors, so $(\mu,\lambda)$ is admissible. Hence, Theorem \ref{thma:EC} is a special case of Theorem \ref{thm:rank2}.
\end{ex}

\begin{warn}
    Throughout, we consistently abuse notation by considering the primes $p_1,p_2,\l_0,\l_1$ of \ref{item:admissible-p primes}, as implicitly part of the data of an admissible pair $(\mu,\lambda)$; see, for example, Notation \ref{notn:admissible-places}.
\end{warn}
\begin{notn}\label{notn:admissible-places}
    Given an admissible pair $(\mu,\lambda)$, we write $p_1,p_2,\l_0,\l_1\in\Z_{>2}$ for the primes guaranteed by Definition \ref{defn:admissible}, and we define the following set of places of $\Q$:
    \[M\coloneqq\{2,\infty\}\cup\{p_1,p_2,\l_1\}.\]
\end{notn}

\subsection{Notational conventions}
\begin{itemize}
    \item Given a ring $R$, we write $x=\square\in R$ to signify that $x\in R$ is an invertible square, i.e. $x=y^2$ for some $y\in\units R$. Similarly, given an ideal $I\subset R$, we will write $x\equiv\square\bmod I$ to signify that $x=\square\in R/I$.
    \item Given an abelian group $A$ and an integer $n$, we write $A/n\coloneqq A/nA$. For example, we will often write $\units\Q/2$ to denote $\units\Q/(\units\Q)^2$, i.e. units in $\Q$ up to squares.
    \item Given a field $F$, we write $G_F$ for (a fixed choice of) its absolute Galois group. If $M$ is a discrete $G_F$-module, we write $\hom^i(F,M)\coloneqq\hom^i(G_F,M)$ for Galois cohomology of $G_F$ acting on $M$. If $E/F$ is an elliptic curve, we further write $\hom^i(F,E)\coloneqq\hom^i(F,E(F^s))$ and $\hom^i(F,E[n])\coloneqq\hom^i(F,E[n](F^s))$, for any $n\ge1$, where $F^s$ is (a fixed choice of) separable closure of $F$.
    \begin{itemize}
        \item Further, if $M$ is a discrete $G_{\Q_p}$-module, we write $\uhom^i(\Q_p,M)\subset\hom^i(\Q_p,M)$ for the subgroup of \define{unramified cohomology classes}, i.e. the kernel
        \[\uhom^i(\Q_p,M)\coloneqq\ker\p{\hom^i(\Q_p,M)\to\hom^i(\ur\Q_p,M)},\]
        where $\ur\Q_p/\Q_p$ is the maximal unramified extension of $\Q_p$.
    \end{itemize}
    \item Let $p$ be a (positive or negative) rational prime. We will generally say that some $x\in\units\Q$ \define{ramifies at $p$} if $\ord_p(x)$ is odd. When $p$ is odd, this is the condition that the cohomology class in $\hom^1(\Q_p,\zmod2)\simeq\units\Q_p/2$ corresponding to $x$ is ramified, i.e. does \important{not} lie in the subgroup $\uhom^1(\Q_p,\zmod2)\simeq\units\Z_p/2$.
    \item Throughout this part of the paper, we make frequent use of the $p$-adic Hilbert symbols $(-,-)_p$, including for $p=\infty$ (where we set $\Q_\infty=\R$). Unless otherwise stated, we always work with multiplicative Hilbert symbols, viewing them as a bilinear map
    \[(-,-)_p\colon\units\Q_p/2\by\units\Q_p/2\too\{\pm1\}.\]
    By definition, given $a,b\in\units\Q_p$, we have
    \[(a,b)_p=+1\iff a\in\Nm_{\Q_p(\sqrt b)/\Q_p}\units{\Q_p(\sqrt b)}\iff b\in\Nm_{\Q_p(\sqrt a)/\Q_p}\units{\Q_p(\sqrt a)}.\]
    We frequently make use, without reference, of the following further properties of Hilbert symbols:
    \begin{itemize}
        \item the \define{Hilbert product formula} (or \define{Hilbert reciprocity}):
        \[\prod_{p\le\infty}(a,b)_p=+1\twhen a,b\in\units\Q.\]
        \item when $p$ is an odd prime, $\pi\in\units\Q_p$ is a uniformizer, and $u\in\units\Z_p$, one has
        \[(u,\pi)_p=+1\iff u\equiv\square\bmod p.\]
        \item when $p$ is an odd prime, and $u_1,u_2\in\units\Z_p$, one has
        \[(u_1,u_2)_p=+1.\]
        \item for any $x\in\units\Q_p$, for any $p$, one has
        \[(x,-x)_p=+1.\]
    \end{itemize}
\end{itemize}

\section{Preliminary machinery}
\subsection{2-descent}
\begin{set}
    Fix distinct $\mu,\lambda\in\Z_{\neq0}$. Consider the elliptic curve $E=E_{\mu,\lambda}\colon y^2=x(x-\mu)(x-\lambda)$.
\end{set}
By ``2-descent,'' we mean the problem of computing/understanding the \define{2-Selmer group} 
\[\Sel_2(E)\coloneqq\ker\p{\hom^1(\Q,E[2])\too\prod_{p\le\infty}\hom^1(\Q_p,E)}.\]
This finite $\zmod2$-module sits in the fundamental exact sequence
\[0\too E(\Q)/2\xtoo\delta\Sel_2(E)\too\Sha(E)[2]\too0,\]
where $\delta\colon E(\Q)/2\into\Sel_2(E)\subset\hom^1(\Q,E[2])$ is the Kummer map and $\Sha(E)$ is the Tate-Shaferavich group of $E$. 
Note that $\dim_{\F_2}\Sel_2(E)\ge\rank E(\Q)+\dim_{\F_2}E[2](\Q)=\rank E(\Q)+2$, so upper bounds on $2$-Selmer give rise to bounds on the rank of $E(\Q)$.

We also introduce the notation
\[e\colon E[2]\by E[2]\too\mu_2\]
for Weil pairing on $E$.
\begin{rem}
    For any nonzero $P,Q\in E[2](\Q)$, one has
    \[e(P,Q)=\twocases{+1}{P=Q}{-1}\]
\end{rem}
\begin{conv}\label{conv:H1(E[2])}
    Fix the ordered basis $P_1=(\mu,0)$, $P_2=(\lambda,0)$ of $E[2]$. We use this to identify $E[2]$ with $\mu_2^2=(\zmod2)^2$ via the isomorphism $E[2]\iso\mu_2^2$, $P\mapsto\p{e(P,P_2),e(P,P_1)}$ identifying $P_1\bij(-1,+1)$ and $P_2\bij(+1,-1)$. This in turn induces an isomorphism
    \[\hom^1(\Q,E[2])\iso\hom^1(\Q,\mu_2^2)=(\units\Q/2)^2,\]
    under which the Kummer map $\delta\colon E(\Q)/2\to(\units\Q/2)^2$ is given by
    \begin{equation}\label{eqn:kum-image}
        \delta((\mu,0))=(\mu(\mu-\lambda),\mu-\lambda)\tcomma\delta((\lambda,0))=(\lambda-\mu,\lambda(\lambda-\mu)),\tand\delta((x,y))=(x-\mu,x-\lambda)\tif x\neq\mu,\lambda.
    \end{equation}
    For the above description of the Kummer map, see, e.g., \cite[Lemma 4.3]{MichaelStoll2001}.
\end{conv}
We will need to exploit further structure on these Selmer groups, related to the Tate local pairing. Fix a (possibly infinite) prime $p$ of $\Q$. The \define{local Tate pairing} is the bilinear map $\angles{-,-}_p\colon\hom^1(\Q_p,E[2])\by\hom^1(\Q_p,E[2])\to\zmod2$ induced by cup product and the Weil pairing:
\begin{equation}\label{eqn:local-Tate-pairing}
    \angles{-,-}_p\colon\hom^1(\Q_p,E[2])\by\hom^1(\Q_p,E[2])\xto\smile\hom^2(\Q_p,E[2]\by E[2])\xto{\push e}\hom^2(\Q_p,\mu_2)\xiso{\opinv_p}\frac12\Z/\Z\simeq\zmod2,
\end{equation}
where $\opinv_p\colon\hom^1(\Q_p,\mu_2)=\Br(\Q_p)[2]\xinto\sim\frac12\Z/\Z$ is the usual invariant map from class field theory. As shown in \cite[Section 4.1]{poonen-rains}, this pairing is induced by a quadratic form
\begin{equation}\label{eqn:local-quad-form}
    q_p\colon\hom^1(\Q_p,E[2])\to\Br(\Q_p)[2]\simeq\zmod2,
\end{equation}
i.e. $\angles{\alpha,\beta}_p=q_p(\alpha+\beta)-q_p(\alpha)-q_p(\beta)$ for all $x,y\in\hom^1(\Q_p,E[2])$.
\begin{rem}\label{rem:tate-local-description}
    It was shown by O'Neil and Skorobogatov \cite{period-index,period-index-erratum} that, in the present context, under Convention \ref{conv:H1(E[2])}, one has
    \[q_p(\alpha,\beta)=(\mu(\mu-\lambda)\alpha,\lambda(\lambda-\mu)\beta)_p\in\zmod2\tfor\alpha,\beta\in\units\Q_p/2,\]
    where $(-,-)_p\in\zmod2$ here denotes the additive $p$-adic Hilbert symbol. Consequently,
    \[\angles{(\alpha_1,\beta_1),(\alpha_2,\beta_2)}_p=(\alpha_1,\beta_2)_p-(\alpha_2,\beta_1)_p\in\zmod2.\]
\end{rem}
\begin{rem}\label{rem:kummer-maximal-isotropic}
    It follows from \cite[Proposition 4.10]{poonen-rains} that $\delta(E(\Q_p))\subset\hom^1(\Q_p,E[2])$, the image of the local Kummer map at $p$, is a \define{maximal isotropic} subgroup with respect to $q_p$, i.e. that $q_p(\delta(E(\Q_p)))=0$ and that $\delta(E(\Q_p))$ is its own orthogonal complement under $\angles{-,-}_p$. This makes $(\hom^1(\Q_p,E[2]),q_p)$ a `weakly metabolic quadratic module' in the sense of \cite[Section 2.1]{poonen-rains}. In particular, $\dim\delta(E(\Q_p))=\frac12\dim\hom^1(\Q_p,E[2])=\dim(\units\Q_p/2)$, which is $2$ if $p$ is an odd prime.
\end{rem}

\subsection{Quadratic twists}\label{sect:quad-twist}

In the present section, we describe how 2-Selmer groups change under quadratic twists.
\begin{defn}
    For $d\in\units\Q/2$, the associated \define{quadratic twist} of an elliptic curve $E=E_{\mu,\lambda}\colon y^2=x(x-\mu)(x-\lambda)$ is the elliptic curve
    \[E^d=E_{d\mu,d\lambda}\colon y^2=x(x-d\mu)(x-d\lambda).\qedhere\]
\end{defn}
In the above definition, there is a natural isomorphism $E^d[2]\simeq E[2]$ which allows one to view $\Sel_2(E^d)$ as a subgroup of $\hom^1(\Q,E[2])$, but potentially cut out by different local conditions than the ones defining $\Sel_2(E)$.
\begin{notn}
    Given a $G_\Q$-module $M$ and a (possibly infinite) prime $p$, we write $\oploc_p\colon\hom^1(\Q,M)\to\hom^1(\Q_p,M)$ for the natural restriction map.
\end{notn}
\begin{defn}\label{defn:selmer}
    Let $M$ be a finite $G_\Q$-module. A \define{Selmer structure} for $M$ is a collection $\mcL=(\mcL_p)_{p\le\infty}$ of subgroups $\mcL_p\subset\hom^1(\Q_p,M)$ such that $\mcL_p=\uhom^1(\Q_p,M)$ for all but finitely many primes $p$. Given $\mcL$, the corresponding \define{Selmer group} is
    \[\Sel(\mcL)\coloneqq\ker\p{\hom^1(\Q,M)\too\bigoplus_{p\le\infty}\frac{\hom^1(\Q_p,M)}{\mcL_p}}=\b{\alpha\in\hom^1(\Q,M)\colon\oploc_p(\alpha)\in\mcL_p\tforall p}.\]
    Furthermore, given a fixed place $\l$ of $\Q$, we define the \define{$\l$-relaxed Selmer group}
    \[\Sel(\mcL)^{(\l)}\coloneqq\ker\p{\hom^1(\Q,M)\to\bigoplus_{p\neq\l}\frac{\hom^1(\Q_p,M)}{\mcL_p}}=\b{\alpha\in\hom^1(\Q,M)\colon\oploc_p(\alpha)\in\mcL_p\tforall p\neq\l},\]
    as well as the \define{$\l$-strict Selmer group}
    \[\Sel(\mcL)_{(\l)}\coloneqq\ker\p{\Sel\mcL\to\hom^1(\Q_\l,M)}=\b{\alpha\in\Sel\mcL\colon\oploc_\l(\alpha)=0}.\qedhere\]
\end{defn}
\begin{notn}\label{notn:S(E)}
    For an elliptic curve $E/\Q$, we set $S(E)\coloneqq\{2,\infty\}\cup\{\t{primes }p\t{ of bad reduction for }E\}$.
\end{notn}
\begin{set}
    Fix distinct $\mu,\lambda\in\Z$, and consider the elliptic curve $E=E_{\mu,\lambda}\colon x(x-\mu)(x-\lambda)$. Fix also a sequence $(v_1,d_1),\dots,(v_r,d_r)$ of pairs consisting of a place\footnote{We assume $v_i\neq v_j$ if $i\neq j$.} $v_i$ with an associated element $d_i\in\Q$. 
    For any $i\in\{0,1,\dots,r\}$, set
    \[\mbfv_i\coloneqq((v_1,d_1),\dots,(v_i,d_i)).\]
    Finally, set $d\coloneqq d_r$.
\end{set}
Associated to this data, we construct a sequence of Selmer structures $\mcL^0,\mcL^1,\dots,\mcL^r$ interpolating between $\Sel_2(E)$ and $\Sel_2(E^d)$.
\begin{notn}
    Given $t\in\units\Q/2$ and a prime $p\le\infty$, let
    \[\delta^t\colon E^t(\Q)/2\into\Sel_2(E^t)\subset\hom^1(\Q,E[2])\tand\delta^t_p\colon E^t(\Q_p)/2\into\hom^1(\Q_p,E[2])\]
    denote the indicated Kummer maps. When $t=1$ we omit it from the notation, so $\delta=\delta^1$ and $\delta_p=\delta_p^1$.
\end{notn}
Our Selmer structures $\mcL^i=\mcL(E,\mbfv_i)$ are defined by
\begin{equation}\label{eqn:Li-def}
    \mcL^i_p=\mcL(E,\mbfv_i)_p\coloneqq\twocases[,]{\im\delta^{d_j}_{v_j}}{p=v_j\tforsome j\in\{1,\dots,i\}}{\im\delta_p}
\end{equation}
for $i=0,1,\dots,r$ and $p\le\infty$.
\begin{lemma}
    $\mcL^i$ is a Selmer structure for all $i$, and $\Sel\mcL^0=\Sel_2(E)$. In particular,
    \[\mcL^i_p=\uhom^1(\Q_p,E[2])\tif p\not\in\{v_1,\dots,v_r\}\cup S(E).\]
\end{lemma}
\begin{proof}
    It is well known that the image of the local Kummer map consists of the unramified classes when $p$ is an odd prime of good reduction, so each $\mcL^i$ is a Selmer structure. When $i=0$, $\mcL^0$ is visibly the Selmer structure defining $\Sel_2(E)$.
\end{proof}
\begin{lemma}\label{lem:Sel Lr}
    If $d\coloneqq d_r$ is only divisible by primes in $\{v_1,\dots,v_r\}$, $d=\square\in\Q_p$ for all $p\in S(E)-\{v_1,\dots,v_r\}$, and $d/d_i=\square\in\Q_{v_i}$ for all $i$, then $\Sel\mcL^r=\Sel_2(E^d)$.
\end{lemma}
\begin{proof}
    We need to show that $\im\delta_p=\im\delta^d_p$ for all $p\not\in\{v_1,\dots,v_r\}$ and that $\im\delta_{v_i}^{d_i}=\im\delta_{v_i}^d$ for all $i$. The latter of these holds simply because $d/d_i=\square\in\Q_{v_i}$, so $E^d\cong E^{d_i}$ over $\Q_{v_i}$ and hence they have the same local Kummer image. Fix $p\not\in\{v_1,\dots,v_r\}$. If $p\in S(E)$, then $d=\square\in\Q_p$, so $E\cong E^d$ over $\Q_p$ and hence $\im\delta_p=\im\delta^d_p$. If $p\not\in S(E)$, then $p$ is an odd prime of good reduction for both $E$ and $E^d$, so $\im\delta_p=\uhom^1(\Q_p,E[2])=\im\delta_p^d$.
\end{proof}
\begin{lemma}\label{lem:twisted-kummer-image}
    Adopt Convention \ref{conv:H1(E[2])}. For any $i\in\{1,\dots,r\}$ with $v_i\not\in S(E)$ and $d_i$ a uniformizer at $v_i$,
    \[\im\delta_{v_i}^{d_i}=\big\langle\p{-\mu d_i,-\lambda d_i},\,\p{\mu(\mu-\lambda), (\mu-\lambda)d_i}\big\rangle\subset(\units\Q_{v_i}/2)^2\simeq\hom^1(\Q_{v_i},E[2]).\]
    In particular, $\im\delta_{v_i}^{d_i}\cap\im\delta_{v_i}=0$.
\end{lemma}
\begin{proof}
    It follows from (\ref{eqn:kum-image}) that $\delta^{d_i}_{v_i}\p{E^{d_i}[2](\Q_{v_i})}$ is generated by $(-\mu d_i,-\lambda d_i)$ and $(\mu(\mu-\lambda), (\mu-\lambda)d_i)$. Because $d_i$ is a uniformizer at $d_i$ (while $\mu,\mu-\lambda$ are units at $v_i$, as $v_i\not\in S(E)$), this space is visibly $2$-dimensional. Since $\im\delta_{v_i}^{d_i}$ is also 2-dimensional (see Remark \ref{rem:kummer-maximal-isotropic}), they must generate the entire image. Every nonzero element of this space is visibly ramified, so it has trivial intersection with $\im\delta_{v_i}=\uhom^1(\Q_{v_i},E[2])$.
\end{proof}
\begin{lemma}\label{lem:naive-relaxed-dimension-bound}
    Let $p$ be any (possibly infinite) prime, fix $i\in\{0,\dots,r-1\}$, and set $\mcL\coloneqq\mcL^i$. Then, $\dim\Sel(\mcL)^{(p)}\le\dim\Sel(\mcL)+\frac12\dim\hom^1(\Q_p,E[2])$.
\end{lemma}
\begin{proof}
    Let $V\coloneqq\oploc_p\p{\Sel(\mcL)^{(p)}}\subset\hom^1(\Q_p,E[2])$. By definition, $\Sel\mcL=\ker\p{\oploc_p\colon\Sel(\mcL)^{(p)}\to V/(V\cap\mcL_p)}$, so $\dim\Sel(\mcL)^{(p)}-\dim\Sel\mcL\le\dim V$. Hence, it suffices to show that $\dim V\le\frac12\dim\hom^1(\Q_p,E[2])$, so it suffices to show that $V$ is isotropic for the local Tate pairing at $p$. Indeed, given $\alpha_p,\beta_p\in V$, lift these to $\alpha,\beta\in\Sel(\mcL)^{(p)}$. Then,
    \[0\overset{\t{Brauer reciprocity}}=\sum_{\l\le\infty}\angles{\alpha,\beta}_\l\overset{\mcL_\l\t{ is an isotropic, Remark \ref{rem:kummer-maximal-isotropic}}}=\angles{\alpha,\beta}_p=\angles{\alpha_p,\beta_p}_p.\qedhere\]
\end{proof}
\begin{prop}\label{prop:Selmer-change}
    Fix $i\in\{0,\dots,r-1\}$. Set $\mcL\coloneqq\mcL^i$, $\mcM\coloneqq\mcL^{i+1}$, and $v\coloneqq v_{i+1}$ to ease notation. Assume $v$ is odd and $\mcL_v\cap\mcM_v=0$. Also set $V\coloneqq\oploc_v\p{\Sel(\mcL)^{(v)}}$. Then, $V$ is a maximal isotropic subspace of $\hom^1(\Q_p,E[2])$; furthermore
    \begin{enumerate}
        \item If $\mcM_v\subset V$ (equivalently, $\mcM_v=V$), then
        \[\dim_{\F_2}\Sel\mcM-\dim_{\F_2}\Sel\mcL=+2\tcomma\Sel\mcM=\Sel(\mcL)^{(v)},\tand\Sel\mcL=\Sel(\mcL)_{(v)}.\]
        \item If $\dim_{\F_2}\oploc_v\p{\Sel\mcL}=2$, then
        \[\dim_{\F_2}\Sel\mcM-\dim_{\F_2}\Sel\mcL=-2\tcomma\Sel\mcM=\Sel(\mcL)_{(v)},\tand\Sel\mcL=\Sel(\mcL)^{(v)}.\]
        \item Otherwise,
        \[\dim_{\F_2}\Sel\mcL=\dim_{\F_2}\Sel\mcM\tand\dim_{\F_2}\p{V\cap\mcM_v}=\dim_{\F_2}(V\cap\mcL_v).\]
    \end{enumerate}
    Note that, by Lemma \ref{lem:twisted-kummer-image}, the assumption $\mcL_v\cap\mcM_v=0$ holds if $v=v_{i+1}\not\in S(E)$ and $d_{i+1}$ is a uniformizer at $v_{i+1}$.
\end{prop}
\begin{proof}[Proof Sketch]
    This all follows from \cite[Proposition 7.2]{Klagsbrun_Mazur_Rubin_2014}. Spelled out a bit, consider the commutative diagram
    \begin{equation}\label{cd:markov}
        \begin{tikzcd}
            \Sel\mcL\ar[d]\ar[r, hook]&\Sel(\mcL)^{(v)}\ar[d, "\oploc_v"]&\Sel\mcM\ar[l, hook']\ar[d]\\
            \mcL_v\ar[r, hook]&\hom^1(\Q_v,E[2])&\mcM_v\ar[l, hook']
            ,
        \end{tikzcd}
    \end{equation}
    both of whose squares are Cartesian. As mentioned in the proof of \cite[Proposition 7.2]{Klagsbrun_Mazur_Rubin_2014}, it is a consequence of global Poitou-Tate duality that $V=\oploc_v\p{\Sel(\mcL)^{(v)}}$ is a maximal isotropic subspace of $\hom^1(\Q_v,E[2])$ with respect to the quadratic form $q_v$ of (\ref{eqn:local-quad-form}), and so $\dim V=2$. Further, $\mcL_v\cap\mcM_v=0$ by assumption. With this in mind, \bp1 and \bp2 follow easily from (\ref{cd:markov}). For \bp3, one uses \cite[Proposition 2.4]{Klagsbrun-Mazur-Rubin:disparity} which implies that $\dim(V\cap\mcL_v)\equiv\dim(V\cap\mcM_v)\bmod2$, so these intersections either both have dimension $1$ or both have dimension $0$ if we are not in cases \bp1 or \bp2.
\end{proof}
The following variant of Proposition \ref{prop:Selmer-change} will allow us to control how Selmer groups change when twisting at certain bad primes.
\begin{prop}\label{prop:Selmer-change-bad}
    Fix $i\in\{0,\dots,r-1\}$. Set $\mcL\coloneqq\mcL^i$, $\mcM\coloneqq\mcL^{i+1}$, and $v\coloneqq v_{i+1}$ to ease notation. Assume $\dim(\mcL_v\cap\mcM_v)=\dim\mcL_v-1$.\footnote{If $v$ is an odd prime, $\dim\mcL_v-1=1$. In the remaining cases, $\dim\mcL_\infty-1=0$ and $\dim\mcL_2-1=2$.} Also set $V\coloneqq\oploc_v\p{\Sel(\mcL)^{(v)}}$. Then, $V$ is a maximal isotropic subspace of $\hom^1(\Q_v,E[2])$; furthermore, exactly one of the following holds:
    \begin{enumerate}
        \item $\Sel\mcL\subsetneq\Sel\mcM$ and $\dim\Sel\mcM-\dim\Sel\mcL=+1$; or
        \item $\Sel\mcL\supsetneq\Sel\mcM$ and $\dim\Sel\mcM-\dim\Sel\mcL=-1$.
    \end{enumerate}
\end{prop}
\begin{proof}
    As in Proposition \ref{prop:Selmer-change}, $V=\oploc_v\p{\Sel(\mcL)^{(v)}}$ being a maximal isotropic is a consequence of global Poitou--Tate duality. Thus, \cite[Proposition 2.4]{Klagsbrun-Mazur-Rubin:disparity} shows that
    \[\dim(V\cap\mcL_v)+\dim(V\cap\mcM_v)\equiv\dim(\mcL_v\cap\mcM_v)+\dim\mcL_v\equiv1\bmod2,\]
    so exactly one of $\dim(V\cap\mcL_v)$ or $\dim(V\cap\mcM_v)$ is odd. Since also
    \[\dim(V\cap\mcL_v\cap\mcM_v)\le\dim(V\cap\mcL_v),\dim(V\cap\mcM_v)\le\dim(V\cap\mcL_v\cap\mcM_v)+1,\]
    we see that either $V\cap\mcL_v=V\cap\mcL_v\cap\mcM_v\subsetneq V\cap\mcM_v$ or vice versa. The claim follows.
\end{proof}
\begin{lemma}\label{lem:last-step-+2}
    Assume that $d\coloneqq d_r$ satisfies the hypotheses of Lemma \ref{lem:Sel Lr}, that $v_r\not\in S(E)$, and the $d=d_r$ is a uniformizer at $v_r$.  Then,
    \[\dim\Sel_2(E^d)=\dim\Sel\mcL^r=\dim\Sel\mcL^{r-1}+2.\]
    In fact, $\Sel_2(E^d)=\Sel\mcL^{r-1}+\delta^d(E^d(\Q))\subset\hom^1(\Q,E[2])$.
\end{lemma}
\begin{proof}
    The first claimed equality holds because Lemma \ref{lem:Sel Lr} tells us that $\Sel_2(E^d)=\Sel\mcL^r$. Let $v=v_r$. The rest of the claim follows from Proposition \ref{prop:Selmer-change}\bp1 since 
    \[V\coloneqq\oploc_v\p{\Sel(\mcL^{r-1})^{(v)}}\supset\oploc_v\p{\Sel\mcL^r}=\oploc_v\p{\Sel_2(E^d)}\supset\oploc_v\delta^d(E^d(\Q)[2])=\delta^d(E^d(\Q_v)[2]),\]
    this last equality being a consequence of Lemma \ref{lem:twisted-kummer-image}.
\end{proof}

\section{Pre-additive combinatorics twists}\label{sect:pre-twist}
\subsection{Rank zero twist}\label{sect:rank-zero-twist}

\begin{set}\label{set:rank-zero}
    Let $\mu,\lambda\in\Z$ be distinct nonzero integers whose product is a square in $\Z$, say $\mu\lambda=\eta^2$. Further assume that $\mu(\mu-\lambda)$ and $\lambda(\lambda-\mu)$ are non-squares in $\Z$. Finally, assume also that $\mu$ is squarefree (so $\mu$ is precisely the squarefree part of $\lambda$). Consider the elliptic curve $E=E_{\mu,\lambda}\colon y^2=x(x-\mu)(x-\lambda)=x(x^2-(\mu+\lambda)x+\mu\lambda)$.
\end{set}
We will construct a twist $E^d$ of $E$ with minimal possible 2-Selmer rank, i.e. $\dim\Sel_2(E^d)=2$ (so $\rank E^d(\Q)=0$). The existence of such a twist will ultimately follow from work of Alex Smith \cite{Smith:BSD-Goldfeld}, so we begin by introducing the necessary context for applying \cite[Theorem 1.12]{Smith:BSD-Goldfeld}.

Let $\phi\colon E\to A$ denote the $2$-isogeny with kernel $E[\phi]$ generated by $(0,0)\in E(\Q)$. Explicitly (see \cite[Example III.4.5]{silverman}), recalling that $\mu\lambda=\eta^2$, we let $A$ be the elliptic curve
\[A\colon Y^2=X\p{X^2+2(\mu+\lambda)X+(\mu-\lambda)^2}=X\p{X+(\mu+\lambda)+2\eta}\p{X+(\mu+\lambda)-2\eta},\]
and let $\phi\colon E\to A$ be the morphism
\[\phi\p{x,y}=\p{\frac{y^2}{x^2},\frac{y(\mu\lambda-x^2)}{x^2}}=\p{\frac{y^2}{x^2},\frac{y(\eta-x)(\eta+x)}{x^2}}.\]
Note that $\phi$ is a balanced isogeny in the sense of \cite[Definition 1.6]{Smith:BSD-Goldfeld}. In fact, using the assumptions of Setup \ref{set:rank-zero}, one can check that $\phi$ is the unique balanced 2-isogeny on $E$, so $E$ is an elliptic curve in Case IV of \cite{Smith:BSD-Goldfeld}. Following \cite[Definition 3.1]{Smith:BSD-Goldfeld}, given any nonzero integer $d$, we define
\[u(d)=u(\phi,d)\coloneqq\sum_{p\le\infty}\p{\dim_{\F_2}W_p(d)-1},\twhere W_p(d)\coloneqq\ker\p{\hom^1(\Q_p,E^d[\phi])\to\hom^1(\Q_p,E^d)}.\]

Above, $\Q_\infty\coloneqq\R$.
As remarked in \cite[Definition 3.1]{Smith:BSD-Goldfeld}, because $\phi$ is a balanced isogeny, the summand $u_p(d)\coloneqq\dim W_p(d)-1$ above vanishes unless $p\in S(E)$ (recall Noation \ref{notn:S(E)}), so $u(d)$ above depends only on the image of $d$ in $\prod_{p\in S(E)}\units\Q_p/2$. To apply \cite[Theorem 1.12]{Smith:BSD-Goldfeld} to produce a rank $0$ twist of $E$, we need a $d_0$ with $u(d_0)\le0$. We construct one using results of \cite{DokchitserDokchitser:local-invariants-isogeneous} and Tate's algorithm.
\begin{lemma}\label{lem:inv-diff-pullback}
    Let $\omega_E\coloneqq \dx/y$ and $\omega_A\coloneqq\dX/Y$ denote the N\'eron differentials on $E,A$, respectively. Then,
    \[\pull\phi\omega_A=-\omega_E.\]
\end{lemma}
\begin{proof}
    We prove this in greater generality. For this proof, fix $a,b\in\Q$ with $b\neq0$ and $r\coloneqq a^2-4b\neq0$. As in \cite[Example III.4.5]{silverman}, consider the elliptic curves
    \[E\colon y^2=x^3+ax^2+bx\tand A\colon Y^2=X^3-2aX+rX,\]
    along with the below 2-isogeny $\phi\colon E\to A$ killing $(0,0)$, 
    \[\phi(x,y)=\p{\frac{y^2}{x^2},\frac{y(b-x^2)}{x^2}}.\]
    Set $\omega_E\coloneqq\dx/y$ and $\omega_A\coloneqq\dX/Y$. Then, $\pull\phi\omega_A=-\omega_E$, as we compute below:
    \begin{align}
        \pull\phi\omega_A
        =\frac{\d\pfrac{y^2}{x^2}}{\frac{y(b-x^2)}{x^2}}
        =\frac{x^2}{y(b-x^2)}\d\p{x+a+b\inv x}
        =\frac{x^2}{y(b-x^2)}\p{1-b\inv[2]x}\dx
        =-\frac\dx y
        =-\omega_E
        .
    \end{align}
    Taking $a=-(\mu+\lambda)$ and $b=\mu\lambda$ proves the original claim.
\end{proof}

\begin{lemma}\label{lem:tamagawa-ratio}
    Let $p$ be a (finite) prime. Let $c_p(E)$ denote the local Tamagawa factor at $p$, i.e. the index $c_p(E)=[E(\Q_p)\colon E_0(\Q_p)]$ of the subgroup $E_0(\Q_p)$ of points whose reduction mod $p$ is nonsingular. Then,
    \[\frac{\#W_p(d)}2=\frac{c_p\p{A^d}}{c_p\p{E^d}}.\]
\end{lemma}
\begin{proof}
    Let $|\cdot|_p$ denote the usual $p$-adic absolute value. Then,
    \cite[Lemma 4.2(2)]{DokchitserDokchitser:local-invariants-isogeneous} applied to $\phi\colon E\to A$ and Lemma \ref{lem:inv-diff-pullback} compute that
    \[\frac{c_p(A^d)}{c_p(E^d)}\overset{\t{Lemma \ref{lem:inv-diff-pullback}}}=\abs{\frac{\pull\phi\omega_A}{\omega_E}}_p^{-1}\frac{c_p(A^d)}{c_p(E^d)}\overset{\t{\cite[Lemma 4.2(2)]{DokchitserDokchitser:local-invariants-isogeneous}}}=\frac{\#\coker\p{\phi\colon E^d(\Q_p)\to A^d(\Q_p)}}{E^d[\phi](\Q_p)}=\frac12\cdot\#\frac{A^d(\Q_p)}{\phi(E^d(\Q_p))}.\]
    Finally, taking cohomology of the exact sequence $0\to E^d[\phi]\to E^d\xto\phi A^d\to0$ induces an isomorphism
    \[\frac{A^d(\Q_p)}{\phi\p{E^d(\Q_p)}}\iso\ker\p{\hom^1(\Q_p,E^d[\phi])\to\hom^1(\Q_p,E^d)}=W_p(d),\]
    from which the claim follows.
\end{proof}
\begin{rec}
    Recall below that, as part of Setup \ref{set:rank-zero}, we are assuming that $\mu$ is squarefree and $\mu\lambda$ is a square. In particular, $\lambda=\mu m^2$ for some $m\in\Z$.
\end{rec}
\begin{lemma}\label{lem:tamagawa-computation}
    Fix some squarefree integer $d$ and a (finite) prime $p$. 
    \begin{romanize}
        \item If $p\neq2$, $p\mid\mu$, and $p\nmid d$, then $E^d$ (so also $A^d$) has additive reduction at $p$, and $c_p(E^d)=4=c_p(A^d)$.
        \item If $p\neq2$, $p\mid\lambda$, and $p\nmid d\mu$, then $E^d$ (so also $A^d$) has multiplicative reduction at $p$, and
        \[c_p(E^d)=\twocases{2\ord_p(\lambda)}{-d\mu\equiv\square\bmod p}2\twhile c_p(A^d)=\twocases{\ord_p(\lambda)}{-d\mu\equiv\square\bmod p}{2}\]

        \item If $p\neq2$, $p\mid(\mu-\lambda)$, and $p\nmid d\mu$, then $E^d$ (so also $A^d$) has multiplicative reduction at $p$, and
        \[c_p(E^d)=\twocases{2\ord_p(\mu-\lambda)}{d\mu\equiv\square\bmod p}{2}\twhile c_p(A^d)=\twocases{4\ord_p(\mu-\lambda)}{d\mu\equiv\square\bmod p}2\]
        

        \item If $p=2$, $d\mu$ is odd, and $d\lambda\equiv4\bmod8$, then $E^d$ (so also $A^d$) has additive reduction at $2$, $c_2(E^d)\ge2$ and $c_2(A^d)\le4$.

    \end{romanize}
\end{lemma}
\begin{proof}
    These all follow from applying Tate's algorithm \cite[Section IV.9]{Silverman:Advanced} to the explicit Weierstrass models
    \[E^d\colon y^2=x\p{x^2-d(\mu+\lambda)x+d^2\mu\lambda}\tand A^d\colon Y^2=X\p{X^2+2d(\mu+\lambda)X+d^2(\mu-\lambda)^2}.\]
    Note that, when $E^d$ has additive reduction at $p\neq2$, one can conclude that $c_p(E^d)=4$ using the facts that $\#E^d[2](\Q_p)=4$ and that $E_0^d(\Q_p)[2]\into\wt E^d_{\t{sm}}(\F_p)[2]\simeq\F_p[2]=0$, where $\wt E^d_{\t{sm}}$ denotes the smooth locus of the reduction mod $p$ of $E^d$. The same argument shows that $c_p(A^d)=4$ in this case as well. 
    
    In cases \bp{ii,iii}, $E^d$ will have multiplicative reduction at $p$; depending on if it is split or not, the local Tamagawa factor $c_p(E^d)$ is either equal to $\ord_p(\Delta(E^d))$ or to the unique integer in $\{1,2\}$ with the same parity as $\ord_p(\Delta(E^d))$; see \cite[Step 2 in Section IV.9]{Silverman:Advanced}. The analogue holds also for $A^d$. One may compute that, in either case (note that $p\nmid d$),
    \[\ord_p(\Delta(E^d))=\ord_p\p{\mu^2\lambda^2(\mu-\lambda)^2}\tand\ord_p(\Delta(A^d))=\ord_p\p{(\mu-\lambda)^4\mu\lambda}.\]
    In particular, $\ord_p(\Delta(E^d))$ is always even, so $c_p(E^d)=2$ when $E^d$ has non-split multiplicative reduction. Similarly, note that $\ord_p(\lambda)$ is even if $p\nmid\mu$ in case \bp{ii}, so also $c_p(A^d)=2$ in this case. 

    With the above observations in mind, cases \bp{i,ii,iii} are straightforward to verify. Let us expand on case \bp{iv}. In this case, upon reaching Step 6 of Tate's algorithm \cite[Section IV.9]{Silverman:Advanced} applied to $E^d$, one must consider the polynomial
    \[P(T)=T\p{T^2-\frac{d\mu+d\lambda+1}2T+\frac{d^2\mu\lambda}4}.\]
    Since $d\lambda\equiv4\bmod8$, $d^2\mu\lambda/4$ is odd, so $0$ is a simple root of $P(T)$. Thus, Tate's algorithm terminates either at Step 6 or Step 7. If it terminates at Step 6, then
    \[c_2(E^d)=1+\#\b{\alpha\in\F_2\colon P(\alpha)=0}\ge2.\]
    Else, $c_2(E^d)$ is computed by a sub-procedure in Step 7. This sub-procedure terminates whenever some auxiliary quadratic polynomial has distinct roots in $\bar\F_2$.
    However, because $E^d$ has no constant term in its Weierstrass form, these auxiliary polynomials always have $0\in\F_2$ as a root, and so whenever one has distinct roots in $\bar\F_2$, both of its roots will be defined over $\F_2$. It follows from this that one must have $c_2(E^d)=4$. Since $A^d$ has additive reduction at $2$ as well (since having additive reduction is an isogeny invariant), it automatically follows that $c_2(A^d)\le4$; see \cite[Remark IV.9.2.2]{Silverman:Advanced}.
\end{proof}

\begin{lemma}\label{lem:u(d0) <= 0}
    Assume $\ord_2(\lambda)=2$ (so $\lambda\equiv4\bmod8$) and $\mu$ is odd. Fix a squarefree integer $d_0$ satisfying all of
    \begin{itemize}
        \item $\gcd(d_0,\mu\lambda(\mu-\lambda))=1$;
        \item $d_0\mu\not\equiv\square\bmod p$ for all primes $p$ which divide $(\mu-\lambda)$ but do not divide $\mu$; and
        \item there exists some prime $\l$ dividing $\lambda$, but not dividing $\mu$, for which $-d_0\mu\equiv\square\bmod\l$.
    \end{itemize}
    Then $u(d_0)\le0$.
\end{lemma}
\begin{proof}
    Recall that
    \[u(d_0)=\sum_{p\le\infty}(\dim_{\F_2}W_p(d_0)-1)=\sum_{p\in S(E)}(\dim_{\F_2}W_p(d_0)-1)=\log_2\prod_{p\in S(E)}\frac{\#W_p(d_0)}2,\]
    so it suffices to show that $\prod_{p\in S(E)}\#W_p(d_0)/2\le1$. For finite $p\in S(E)$ (i.e. $p\mid\mu\lambda(\mu-\lambda)$), the assumptions on $d_0$ are chosen to ensure that Lemmas \ref{lem:tamagawa-ratio} and \ref{lem:tamagawa-computation} together imply that
    \[\frac{\#W_p(d_0)}2\le1\tif p\neq\l,2\tcomma \frac{\#W_{\l}(d_0)}2=\frac12,\tand\frac{\#W_2(d_0)}2\le2.\]
    Hence, $\prod_{p\in S(E)}\#W_p(d_0)/2\le\#W_\infty(d_0)$, so it only remains to handle the real place. Here, we simply have $\#W_\infty(d_0)\le\#\hom^1(\R,E^{d_0}[\phi])=\#\hom^1(\R,\zmod2)=\#\Hom(\Gal(\C/\R),\zmod2)=2$.
\end{proof}
\begin{cor}\label{cor:rank-0-twist}
    Assume $\ord_2(\lambda)=2$ and $\mu$ is odd, and fix $d_0$ satisfying the hypotheses of Lemma \ref{lem:u(d0) <= 0}. Then, there exists infinitely many squarefree integers $d$ such that $d=d_0\in\prod_{p\in S(E)}\units\Q_p/2$ and $\dim\Sel_2(E^{d})=2$ (i.e. $\Sel_2(E^{d})=\kumtwist{d}$).
\end{cor}
\begin{proof}
    As remarked earlier, the elliptic curve $E$ is in Case IV of \cite{Smith:BSD-Goldfeld}, so \cite[Theorem 1.12]{Smith:BSD-Goldfeld} applies. By assumption, $d_0$ satisfies the hypotheses of Lemma \ref{lem:u(d0) <= 0}, so $u(d_0)\le0$.
    It then follows from \cite[Theorem 1.12]{Smith:BSD-Goldfeld} applied to $E$ with this choice of $d_0$ that there exists infinitely many squarefree $d$ such that $d=d_0\in\prod_{p\in S(E)}\units\Q_p/2$ and $r_2(E^d)\coloneqq\dim\Sel_2(E^d)-\dim E[2](\Q)=0$. 
\end{proof}

\subsection{Padding out $2$-Selmer}

Starting from a suitable rank zero elliptic curve $E$, the goal of the present section is to construct a twist $E^d$ along with an explicit basis for $\Sel_2(E^d)$ designed to be easily manipulated by the additive combinatorics primes we will be twisting by in the final section (see Corollary \ref{cor:4-prepared} and Definition \ref{defn:n-prepared}). The construction of such a basis is where we will make use of the primes alluded to in the definition of admissibility (Definition \ref{defn:admissible}), and the ``shape'' of said basis will be dictated by the ramified elements appearing in the local Selmer conditions at admissible primes, along with the primes $2$ and $\infty$.

\begin{rec}
    For an elliptic curve of the form $E_{\mu,\lambda}\colon y^2=x(x-\mu)(x-\lambda)$, we adopt Convention \ref{conv:H1(E[2])} in identifying $\hom^1(\Q,E[2])\simeq(\units\Q/2)^2$.
\end{rec}
\begin{lemma}\label{lem:Selmer-condition-at-special-p}
    Fix a pair of distinct $\mu,\lambda\in\Z_{\neq0}$. Let $p$ be an odd prime such that $p\nmid\mu\lambda$, but $\ord_p(\mu-\lambda)$ is odd, and let $E=E_{\mu,\lambda}\colon y^2=x(x-\mu)(x-\lambda)$. Then,
    \[\delta(E(\Q_p))\cap\uhom^1(\Q_p,E[2])\t{ is }1\t{-dimensional}.\]
    Furthermore,
    \begin{itemize}
        \item if $\mu\equiv\square\bmod p$, then $\delta(E(\Q_p))\subset\hom^1(\Q_p,E[2])\simeq(\units\Q_p/2)^2$ is the subgroup of diagonal classes.
        \item if $\mu\not\equiv\square\bmod p$, then $\delta(E(\Q_p))$ is generated by $(\mu,\mu)$ and $(\mu p,p)$ in $(\units\Q_p/2)^2$.
    \end{itemize}
\end{lemma}
\begin{proof}
    Let $U\coloneqq\delta(E(\Q_p))\cap\uhom^1(\Q_p,E[2])\subset\hom^1(\Q_p,E[2])\simeq(\units\Q_p/2)^2$. We first show that $\dim U=1$. Observe that $\delta(E(\Q_p))\ni\delta((\mu,0))\overset{(\ref{eqn:kum-image})}=(\mu(\mu-\lambda),\mu-\lambda)\in(\units\Q_p/2)^2$ which is ramified at $p$. Thus, $\dim U\le 1=\dim\delta(E(\Q_p))-1$. To show that $\dim U\ge1$, consider the N\'eron model $\mcE/\Z_p$ of $E$. Let $\mcE^0/\Z_p$ denote its identity component. It is well known that the image of the composition
    \[\mcE^0(\Z_p)/2\to\mcE(\Z_p)/2=E(\Q_p)/2\xinto\delta\hom^1(\Q_p,E[2])\]
    lands in $\uhom^1(\Q_p,E[2])$, so it suffices to show that $\mcE^0(\Z_p)\to\mcE(\Z_p)/2$ is nonzero. For this, note first that $\Phi\coloneqq\mcE(\Z_p)/\mcE^0(\Z_p)$ is cyclic because $E$ has multiplicative reduction at $p$ (see Lemma \ref{lem:tamagawa-computation}). Consider the exact sequence
    \[\mcE^0(\Z_p)/2\too\mcE(\Z_p)/2\too\Phi/2\too0.\]
    Above, $\mcE(\Z_p)=E(\Q_p)$ contains a finite-index pro-$p$ subgroup, so $\#\mcE(\Z_p)/2=\#\mcE(\Z_p)[2]=4$. At the same time, $\Phi$ is cyclic, so $\#\Phi/2\le\#\zmod2=2$, so $\im(\mcE^0(\Z_p)/2\to\mcE(\Z_p)/2)$ must be nonzero and hence $\dim U\ge1$.

    To compute $\delta(E(\Q_p))\subset(\units\Q_p/2)^2$, we combine the facts: $\dim U=1$, $(\mu(\mu-\lambda),\mu-\lambda)\in\delta(E(\Q_p))$, and $\delta(E(\Q_p))$ is a maximal isotropic for the quadratic form 
    \[\mapdesc q{(\units\Q_p/2)^2}{\zmod2}{(\alpha,\beta)}{\p{\mu(\mu-\lambda)\alpha,\lambda(\lambda-\mu)\beta}_p}\]
    (see Remark \ref{rem:tate-local-description}). Let $\notsquare\in\units\Q_p/2$ denote the integral non-square class. One may compute that the only unramified class in $(\units\Q_p/2)^2$ orthogonal to $(\mu(\mu-\lambda),\mu-\lambda)$ under the above quadratic form is $(\notsquare,\notsquare)$. Thus, we must have $\delta(E(\Q_p))=\angles{(\mu(\mu-\lambda),\mu-\lambda),(\notsquare,\notsquare)}$, from which the claim follows.
\end{proof}

\begin{lemma}\label{lem:Selmer-condition-at-special-ell}
    Fix a pair of distinct $\mu,\lambda\in\Z_{\neq0}$. Let $\l$ be an odd prime such that $\l^2\mid\lambda$, but $\l\nmid\mu$, and let $E=E_{\mu,\lambda}\colon y^2=x(x-\mu)(x-\lambda)$. Assume also that $\ord_\l(\lambda)$ is even.
    \begin{itemize}
        \item if $-\mu\equiv\square\bmod\l$, then $\delta(E(\Q_\l))\subset\hom^1(\Q_\l,E[2])\simeq(\units\Q_\l/2)^2$ is the subgroup of classes $(\alpha,\beta)$ for which $\alpha=\square\in\Q_\l$.
        \item if $-\mu\not\equiv\square\bmod\l$, then $\delta(E(\Q_\l))=\uhom^1(\Q_\l,E[2])$ is the subgroup of unramified classes.
    \end{itemize}
\end{lemma}
\begin{proof}
    First assume $-\mu\not\equiv\square\bmod\l$. Note that $\delta(E(\Q_\l))\ni\delta((0,0))\overset{(\ref{eqn:kum-image})}=(-\mu,-\lambda)\in(\units\Q_\l/2)^2$, which has nontrivial first coordinate and is unramified by assumption. To finish this case, it suffices to construct a nontrivial unramified class in $\im\delta$ with trivial first coordinate. For this, we claim there exists some $a\in\F_\l$ such that $a^2+\mu$ is a nonsquare mod $\l$. Indeed, there are $(\l+1)/2$ choices of $a^2$, but only $(\l-1)/2$ possible nonzero squares mod $\l$ (note $a^2+\mu\neq0\in\F_\l$ because $-\mu\neq\square\in\F_\l$), so there must be some $a$ such that $a^2+\mu$ is a nonsquare. Fix $\alpha\in\units\Z_\l$ so that $\alpha^2+\mu\not\equiv\square\bmod\l$. Note that $(\alpha^2+\mu)(\alpha^2+\mu-\lambda)$ is a square in $\Z_\l$, so set $y\coloneqq\alpha\sqrt{(\alpha^2+\mu)(\alpha^2+\mu-\lambda)}$. Then, $(\alpha^2+\mu,y)\in E(\Q_\l)$, and (\ref{eqn:kum-image}) computes that $\delta((\alpha^2+\mu,y))=(\alpha^2,\alpha^2+\mu-\lambda)\in(\units\Q_\l/2)^2$; this class is unramified with trivial first coordinate and nontrivial second coordinate.

    Now assume instead that $-\mu\equiv\square\bmod\l$. Let $U\coloneqq\delta(E(\Q_\l))\cap\uhom^1(\Q_\l,E[2])\subset\hom^1(\Q_\l,E[2])\simeq(\units\Q_\l/2)^2$. We claim that $\dim U=1$. Since $E$ here has multiplicative reduction at $\l$ (see Lemma \ref{lem:tamagawa-computation}), the same argument as in the proof of Lemma \ref{lem:Selmer-condition-at-special-p} shows that $\dim U\ge1$. To show that $\dim U\le1$, we first observe that $(\l-\mu)(1-\lambda/\l)\in\units\Z_\l$ and is congruent to $-\mu\bmod\l$. Hence, we may set $y\coloneqq \l\sqrt{(\l-\mu)(1-\lambda/\l)}\in\Q_\l$ and note that $(\l,y)\in E(\Q_\l)$. By (\ref{eqn:kum-image}), $\delta((\l,y_\l))=(\l-\mu,\l-\lambda)\in(\units\Q_\l/2)^2$. This class is ramified (in its second coordinate), so $\dim U\le1$. We also note that $\l-\mu,1-\lambda/\l=\square\in\Z_\l$, so $(\l-\mu,\l-\lambda)=(1,\l)\in(\units\Q_\l/2)^2$. Now, recall that $\delta(E(\Q_\l))$ is a maximal isotropic for the quadratic form
    \begin{equation}
        \mapdesc q{(\units\Q_\l/2)^2}{\zmod2}{(\alpha,\beta)}{\p{\mu(\mu-\lambda)\alpha,\lambda(\lambda-\mu)\beta}_\l}
    \end{equation}
    (see Remark \ref{rem:tate-local-description}). The only unramified class in $(\units\Q_\l/2)^2$ orthogonal to $(1,\l)$ under $q$ is $(1,\notsquare)$. Thus,
    \[\delta(E(\Q_\l))=\angles{(\l-\mu,\l-\lambda),(1,\notsquare)}=\angles{(1,\l),(1,\notsquare)},\]
    the subspace of classes with trivial first coordinate, in this case.
\end{proof}

\begin{lemma}\label{lem:Selmer-condition-at-2}
    Fix distinct $\mu,\lambda\in\Z_{\neq0}$, and set $E=E_{\mu,\lambda}\colon y^2=x(x-\mu)(x-\lambda)$. Assume $\mu\equiv1\bmod4$, $\lambda=\mu m^2$ for some $m\in\Z$ with $m\equiv2\bmod8$. Then, after identifying $(\units\Q_2/2)^2\simeq\hom^1(\Q_2,E[2])$, as in Convention \ref{conv:H1(E[2])},
    \[\delta(E(\Q_2))=\angles{(5,5),(3,3),(5,1)}\subset(\units\Q/2)^2\]
    Explicitly, this is the subgroup of $(\units\Q_2/2)^2$ with elements
    \[\begin{matrix}
        (1,1) & (5,5) & (7,7) & (3,3) \\
        (5,1) & (1,5) & (3,7) & (7,3)
    \end{matrix}\]
\end{lemma}
\begin{proof}
    From Remark \ref{rem:kummer-maximal-isotropic}, we know that $\dim\delta(E(\Q_2))=\dim(\units\Q_2/2)=3$, so it suffices to show that the three claimed generators belong to this group. Equation (\ref{eqn:kum-image}) immediately tells us that $\delta(E(\Q_2))$ contains
    \[\delta((\mu,0))=(\mu(\mu-\lambda),\mu-\lambda)=\Twocases{(5,5)}{\mu\equiv1\bmod8}{(5,1)}{\mu\equiv5\bmod8}\tand\delta((0,0))=(-\mu,-\lambda)=\Twocases{(7,7)}{\mu\equiv1\bmod8}{(3,3)}{\mu\equiv5\bmod8}.\]
    Suppose now that $\mu\equiv5\bmod8$. Set $x_5\coloneqq\mu+\lambda$ and note that
    \begin{equation}\begin{split}
        x_5(x_5-\mu)(x_5-\lambda) = (\mu+\lambda)\mu\lambda = m^2\mu^2&(\mu+\lambda) \tand \\
        &(\mu+\lambda)\equiv1\bmod8,
    \end{split}\end{equation}
    so we may set $y_5\coloneqq m\mu\sqrt{\mu+\lambda}\in\units\Q_2$ and compute $\delta((x_5,y_5))\overset{(\ref{eqn:kum-image})}=(x_5-\mu,x_5-\lambda)=(5,5)\in(\units\Q_2/2)^2$, finishing the proof in this case.

    Suppose instead that $\mu\equiv1\bmod8$. Set $x_1\coloneqq\mu+5m^2$, and note that
    \begin{equation}\begin{split}
        x_1(x_1-\mu)(x_1-\lambda) = (\mu+5m^2)5m^2(\mu+5m^2-\lambda) = m^2&5(\mu+5m^2)(\mu+5m^2-\lambda) \tand \\
        &5(\mu+5m^2)(\mu+5m^2-\lambda)\equiv1\bmod8,
    \end{split}\end{equation}
    so we may set $y_1\coloneqq m\sqrt{5(\mu+5m^2)(\mu+5m^2-\lambda)}$ and compute $\delta((x_1,y_1))\overset{(\ref{eqn:kum-image})}=(5,1)\in(\units\Q_2/2)^2$, finishing the proof.
\end{proof}

\begin{lemma}\label{lem:Selmer-condition-at-infty}
    Fix distinct $\mu,\lambda\in\Z_{\neq0}$ with $\abs\mu<\abs\lambda$. Let $E=E_{\mu,\lambda}\colon y^2=x(x-\mu)(x-\lambda)$.
    \begin{itemize}
        \item If $\mu,\lambda>0$, then $\delta(E(\R))\subset\hom^1(\R,E[2])\simeq(\units\R/2)^2$ is spanned by $(-1,-1)$.
        \item If $\mu,\lambda<0$, then $\delta(E(\R))\subset\hom^1(\R,E[2])\simeq(\units\R/2)^2$ is spanned by $(-1,+1)$.
    \end{itemize}
\end{lemma}
\begin{proof}
    It follows from Remark \ref{rem:kummer-maximal-isotropic} that $\dim\delta(E(\R))=1=\dim(\units\R/2)$, so it suffices to exhibit a single nonzero element. By (\ref{eqn:kum-image}), this subspace contains both $(-\mu,-\lambda)$ and $(\mu(\mu-\lambda),\mu-\lambda)$, from which the claim follows.
\end{proof}

\begin{defn}\label{defn:n-prepared}
    Let $(\mu,\lambda)\in\Z^2$ be an admissible pair, in the sense of Definition \ref{defn:admissible}. Adopt Notation \ref{notn:admissible-places}. Fix $n\in\{0,1,2,3,4\}$.
    We say that $(\mu,\lambda)$ -- or, equivalently, $E=E_{\mu,\lambda}\colon y^2=x(x-\mu)(x-\lambda)$ -- is \define{$n$-prepared} if all of the following hold.
    \begin{enumerate}[label=\bfseries{(P$_{\arabic*}$)}]
        \item\label{item:n-prepared-odd factor} there exists an odd prime $\varpi$ dividing $\mu$;
        \item\label{item:n-prepared-hilbert-cond} \[\prod_{p\mid(1-m^2)}(\mu,1-m^2)_p=-1\iff n\neq2;\]
        \item\label{item:n-prepared-mu negative} $\phantom-\mu<0$ (so also $\lambda<0$) $\iff n\ge1$;
        \item\label{item:n-prepared-square at p} $\phantom-\mu\equiv\square\bmod p_i\iff n\ge i+1$, for $i\in\{1,2\}$;
        \item\label{item:n-prarared-square at l} $-\mu\equiv\square\bmod\l_1\iff n=4$;
        \item\label{item:n-prepared-mu at 2} $\phantom-\mu\equiv1\bmod8$;
        \item\label{item:n-prepared-gammas} there exists $\{\Gamma_i\}_{i=1}^4\subset\units\Q$ such that $\Sel_2(E)=\kumtwist{}\oplus\angles{\gamma_1,\dots,\gamma_n}$, where
        \[\gamma_1=(\Gamma_1,1)\quad\gamma_2=(\Gamma_2,\Gamma_2)\quad\gamma_3=(\Gamma_3,\Gamma_3)\quad\gamma_4=(1,\Gamma_4).\]
        In particular, $\dim\Sel_2(E)=2+n$; and
        \item\label{item:n-prepared-gammas property} 
        \begin{equation}\begin{aligned}
            &\Gamma_1<0, &&\tbut\Gamma_1=\square\in\Q_v\tforall v\in\{\varpi\}\cup M-\{\infty\},\\
            &\Gamma_2\t{ is a uniformizer at }p_1, &&\tbut\Gamma_2=\square\in\Q_v\tforall v\in\{\varpi\}\cup M-\{p_1\},\\
            &\Gamma_3\t{ is a uniformizer at }p_2\t{ and }\Gamma_3\equiv7\bmod8, &&\tbut\Gamma_3=\square\in\Q_v\tforall v\in\{\varpi\}\cup M-\{p_2,2\},\\
            &\Gamma_4\t{ is a uniformizer at }\l_1, &&\tbut\Gamma_4=\square\in\Q_v\tforall v\in\{\varpi\}\cup M-\{\l_1\}
            .
        \end{aligned}\end{equation}
        \item\label{item:n-prepared-q primes} More specifically, we require that there are primes $q_1,q_2,q_3,q_4$ -- with $q_1,\dots,q_n$ each dividing $\mu$ -- such that
        \begin{equation}\begin{split}
            \Gamma_1&=-q_1\quad\Gamma_2=p_1q_2\quad\Gamma_3=p_2q_3\quad\Gamma_4=\l_1q_4\\
            1-m^2&\equiv\square\bmod q_2,q_3,q_4\\
            1-m^2&\not\equiv\square\bmod q_1\\
            2&\equiv\square\bmod q_2,q_3,q_4,
        \end{split}\end{equation}
        where $m=\sqrt{\lambda/\mu}>0$, as in \ref{item:admissible-square product}.
    \end{enumerate}
\end{defn}
\begin{notn}
    In everything below, given $(\mu,\lambda)\in\Z^2$, we (implicitly) set $E=E_{\mu,\lambda}\colon y^2=x(x-\mu)(x-\lambda)$.
\end{notn}
\begin{lemma}\label{lem:n=0}
    Let $(\mu,\lambda)$ be an admissible pair. Then, there exists some squarefree $t\in\Z$ such that $E^t$ is $0$-prepared; equivalently, $(t\mu, t\lambda)$ is $0$-prepared.
\end{lemma}
\begin{proof}
    Fix $d_0$ as in Lemma \ref{lem:u(d0) <= 0} with the extra conditions that $d_0\mu\equiv1\bmod8$, $-d_0\mu\equiv\square\bmod\l_0$, $-d_0\mu\not\equiv\square\bmod\l_1$, and $d_0\mu>0$; such a $d_0$ exists by the Chinese remainder theorem. Now, Corollary \ref{cor:rank-0-twist} produces some squarefree integer $d$ such that $d\mu\neq1$, $d=d_0\in\prod_{p\in S(E)}\units\Q_p/2$, and $\dim\Sel_2(E^d)=2$. If $\gcd(d,1+m^2)=1$, then we win with $t=d$. Else, we inductively apply the following claim (at each step, replacing $d$ with $d'\coloneqq dt'/p^2$).

    Fix a prime $p\mid\gcd(d,1+m^2)=\gcd(d\mu,1+m^2)$. We claim there exists some squarefree $t'$ such that $p\mid t'$, $t'=1\in\units\Q_v/2$ for all $v\in S(E^d)-\{p\}$ and $\dim\Sel_2(E^{dt'})=2$. We construct this $t'$ as the end result of a sequence $(v_1,d_1),\dots,(v_r,d_r)$ as in Section \ref{sect:quad-twist}. Below, given $i\in\{0,1,\dots,r\}$, we set
    \begin{equation}
        \mbfv_i\coloneqq\p{(v_1,d_1),\dots,(v_i,d_i)}\tcomma\mcL^i\coloneqq\mcL(E^d,\mbfv_i)\t{ as in }(\ref{eqn:Li-def}),\tand S(\mcL^i)\coloneqq\{v_1,\dots,v_i\}\cup S(E^d).
    \end{equation}
    \begin{enumerate}
        \item First use Dirichlet's theorem to produce a positive prime $v_1\not\in S(E^d)=S(\mcL^0)$ such that
        \begin{equation}
            v_1\equiv1\bmod8v\tforall v\in S(\mcL^0)-\{p,\infty\},\tbut v_1\not\equiv\square\bmod p.
        \end{equation}
        Then, $-1\equiv\square\bmod v_1$, while
        \[(d\mu,v_1)_{v_1}=\prod_{v\neq v_1}(d\mu,v_1)_v=(d\mu,v_1)_p=-1\tandso d\mu\not\equiv\square\bmod v_1.\]
        Note that (\ref{eqn:kum-image}) applied to $E^d=E_{\mu d,\lambda d}$ shows that $\kumtwist d$ is generated by $(-d\mu,-d\lambda)$ and $(\mu(\mu-\lambda),d(\mu-\lambda))$ in $(\units\Q_{v_1}/2)^2$. We thus conclude that $\oploc_{v_1}\kumtwist d$ is 2-diemsnional. Set $d_1\coloneqq v_1$. Then, Proposition \ref{prop:Selmer-change} shows that $\dim\Sel\mcL^1=\dim\Sel_2(E)-2=0$, so $\Sel\mcL^1=0$.

        \item Set $v_2\coloneqq p$ and $d_2\coloneqq p$. Since $p\mid d$ and $p\nmid\mu$, we see that the $\mcL^2$-Selmer condition at $p$ is
        \[\mcL^2_p=\kumtwist[p]{dp}=\kumtwist[p]{d/p}\overset*=\uhom^1(\Q_p,E[2])=(\units\Z_p/2)^2\subset(\units\Q_p/2)^2,\]
        where the quality labeled $\overset*{}$ holds because $E^{d/p}$ has good reduction at $p$.

        \item Suppose that $\Sel\mcL^2\neq0$ (else, skip this step). Since $d\mu$ is a uniformizer at $p$, (\ref{eqn:kum-image}) shows us that
        \[\im\delta_p^d=\angles{(-d\mu,-d\lambda),(d(\mu-\lambda),\mu(\mu-\lambda)}\subset(\units\Q_p/2)^2,\]
        and in particular that no nontrivial element of $\im\delta_p^d$ is unramified at $p$. 
        It follows from Proposition \ref{prop:Selmer-change} that we must then have $\dim\Sel\mcL^2=2+\dim\Sel\mcL^1=2$, and moreover that $\Sel\mcL^2=\angles{\theta,\theta'}$ for some $\theta,\theta'\in(\units\Q/2)^2$ such that $\oploc_{v_2}\angles{\theta,\theta'}=\im\delta^{dd_2}_{v_2}=(\units\Z_{v_2}/2)^2\subset(\units\Q_{v_2}/2)^2$. Hence, we may and do choose $\theta=(\theta_1,\theta_2)$ and $\theta'=(\theta_1',\theta_2')$ in $(\units\Q/2)^2$ so that
        \begin{equation}
            \oploc_{v_2}\theta=(\notsquare,\square)\tand\oploc_{v_2}\theta'=(\square,\notsquare)\tin(\units\Q_{v_2}/2)^2,
        \end{equation}
        where, above, $\square,\notsquare\in\units\Z_{v_2}/2$ denote the square and non-square classes. Thus, Chebotarev density applied to $\Q(\sqrt{\theta_1},\sqrt{\theta_2},\sqrt{\theta_1'},\sqrt{\theta_2'})/\Q$ produces for us a positive prime $v_3\not\in S(\mcL^2)$ such that $\theta_1,\theta_2'\not\equiv\square\bmod v_3$, but $\theta_2,\theta_1'\equiv\square\bmod v_3$. Set $d_3\coloneqq v_3$. It follows that $\oploc_{v_3}\Sel\mcL^2=\oploc_{v_3}\angles{\theta,\theta'}$ is 2-dimensional, so Proposition \ref{prop:Selmer-change} shows that $\dim\Sel\mcL^3=\dim\Sel\mcL^2-2=0$.
        \item Finally, set $r\coloneqq3$ if $\Sel\mcL^2=0$ and $r\coloneqq4$ otherwise. Thus, $\Sel\mcL^{r-1}=0$. Dirichlet supplies us with a positive prime $v_r\not\in S(\mcL^{r-1})$ such that $t'=d_r\coloneqq v_1v_2\dots v_r$ satisfies
        \[t'\equiv1\bmod8v\tforall v\in S(E^d)-\{p\}\tand\frac{t'}{v_i}=\prod_{\substack{j=1\\ j\neq i}}^rv_j\equiv1\bmod v_i\tforall i\in\{1,\dots,r-1\}.\]
        It follows from Lemma \ref{lem:last-step-+2} that $\dim\Sel_2(E^{dt'})=\dim\Sel\mcL^r=\dim\Sel\mcL^{r-1}+2=2$.
    \end{enumerate}
    This completes the proof of the subclaim. By applying it inductively, replacing $d$ with $d'=dt'/p^2$ after each step (note $p\nmid d'$), we eventually arrive at a $d$ which still satisfies $d=d_0\in\prod_{p\in S(E)}\units\Q_p/2$ and $\dim\Sel_2(E^d)=2$, but now also satisfies $\gcd(d,1+m^2)=1$. Take $t=d$. Most properties of Definition \ref{defn:n-prepared} (with $n=0$) are clearly satisfied by $E^t$ (i.e. by $(t\mu,t\lambda)$). For \ref{item:n-prepared-hilbert-cond}, note that this property follows from \ref{item:admissible-odd-num-factors} and the fact that $d\mu\not\equiv\square\bmod p$ for any $p\mid(1-m^2)$, by construction.
\end{proof}
\begin{notn}
    In everything below, given a sequence $(v_1,d_1),\dots,(v_r,d_r)$ as in Section \ref{sect:quad-twist} (which will be clear from context) and $i\in\{0,1,\dots,r\}$, we set
    \begin{equation}
        \mbfv_i\coloneqq\p{(v_1,d_1),\dots,(v_i,d_i)}\tcomma\mcL^i\coloneqq\mcL(E,\mbfv_i)\t{ as in }(\ref{eqn:Li-def}),\tand S(\mcL^i)\coloneqq\{v_1,\dots,v_i\}\cup S(E).
    \end{equation}
\end{notn}
\begin{lemma}\label{lem:n+1}
    Fix $n\in\{1,2\}$. Let $(\mu,\lambda)$ be an $n$-prepared admissible pair. Then, there exists some $t$ such that $E^t$ is $(n+1)$-prepared.
\end{lemma}
\begin{proof}
    We adopt all the notation of Definition \ref{defn:n-prepared}. Set $p\coloneqq p_n$. By \ref{item:n-prepared-gammas}, we may and do write
    \[\Sel_2(E)=\kumtwist{}\oplus\angles{\gamma_1,\dots,\gamma_n},\]
    where $\gamma_1=(\Gamma_1,1)$, $\gamma_2=(\Gamma_2,\Gamma_2)$, and these data satisfy \ref{item:n-prepared-gammas property}.
    We will construct $t=d_r$ as the end result of a sequence $(v_1,d_1),\dots,(v_r,d_r)$ as in Section \ref{sect:quad-twist}.
    \begin{enumerate}
        \item By Dirichlet's theorem, there exists a positive prime $v_1\not\in S(E)=S(\mcL^0)$ such that
        \begin{equation}\label{eqn:n+1-v1}
            v_1\equiv1\bmod8v\tforall v\in S(\mcL^0)-\{\varpi,\infty\},\tbut v_1\not\equiv\square\bmod\varpi,
        \end{equation}
        with $\varpi$ as in \ref{item:n-prepared-odd factor}. Thus, $-1\equiv\square\bmod v_1$ and 
        \[(\mu,v_1)_{v_1}=\prod_{v\neq v_1}(\mu,v_1)_v=(\mu,v_1)_\varpi=-1,\]
        so $\mu\not\equiv\square\bmod v_1$. From (\ref{eqn:kum-image}), we conclude that $\oploc_{v_1}\delta(E[2](\Q))\subset(\units\Q_{v_1}/2)^2$ is 2-dimensional. Set $d_1=v_1$. Then, Proposition \ref{prop:Selmer-change} shows that $\dim\Sel\mcL^1=\dim\Sel_2(E)-2$ and that $\Sel\mcL^1=\ker\p{\oploc_{v_1}\colon\Sel_2(E)\to\hom^1(\Q_{v_1},E[2])}$. We claim that
        \begin{equation}\label{eqn:n+1-Sel L1}
            \Sel\mcL^1=\angles{\gamma_1,\dots,\gamma_n}
        \end{equation}
        By the above work, it suffices to show that $\oploc_{v_1}(\gamma_i)$ is trivial for each $i$. This follows from (\ref{eqn:n+1-v1}) and \ref{item:n-prepared-gammas property} which allow us to compute that
        \begin{align}
            (\Gamma_i,v_1)_{v_1} &= \prod_{v\neq v_1}(\Gamma_i, v_1)_v = (\Gamma_i,v_1)_\varpi = +1 \tforall i\in\{1,\dots,n\},
        \end{align}
        so $\oploc_{v_1}(\gamma_i)$ is indeed trivial, establishing (\ref{eqn:n+1-Sel L1}).

        \item Using Dirichlet, we fix a positive prime $v_2\not\in S(\mcL^1)$ such that
        \begin{equation}\label{eqn:n+1-v2}\begin{split}
            p&v_2 \equiv 1\bmod v\tforall v\in S(\mcL^1)-\{p,2,\infty\},\tbut \\
            &v_2\not\equiv\square\bmod p\tand\Twocases{pv\equiv1\bmod8}{n=1}{pv\equiv7\bmod8}{n=2}.
        \end{split}\end{equation}
        Set $d_2\coloneqq v_2$, $\Gamma_{n+1}\coloneqq pv_2$, and $\gamma_{n+1}\coloneqq(\Gamma_{n+1},\Gamma_{n+1})\in(\units\Q/2)^2$. We claim that
        \begin{equation}\label{eqn:n+1-Sel L2}
            \Sel\mcL^1\subset\Sel\mcL^2\tand \gamma_{n+1}\in\Sel(\mcL^2)^{(p)},
        \end{equation}
        where the superscript ${}^{(p)}$ denotes the $p$-relaxed Selmer group, Definition \ref{defn:selmer}.
        The first part of (\ref{eqn:n+1-Sel L2}) will follow from showing that $\oploc_{v_2}\p{\Sel\mcL^1}=0$, and, because of (\ref{eqn:n+1-v2}) and Lemma \ref{lem:Selmer-condition-at-2}, the second part will follow from showing that $\gamma_{n+1}$ satisfies the $\mcL^2$-Selmer condition at $v_2$. For the first part, we compute
        \[(\Gamma_i,\Gamma_{n+1})_{v_2}\overset{(\ref{eqn:n+1-v2})}=(\Gamma_i,\Gamma_{n+1})_p(\Gamma_i,\Gamma_{n+1})_2\overset{\t{\ref{item:n-prepared-gammas property}}}=+1\tforall i\in\{1,\dots,n\},\]
        so $\Gamma_i\equiv\square\bmod v_2$. 
        Further, we claim that $\oploc_{v_2}(\gamma_{n+1})=(-\mu d_2,-\lambda d_2)\overset{(\ref{eqn:kum-image})}=\delta^{d_2}_{v_2}((0,0))\in(\units\Q_{v_2}/2)^2$. Indeed, by \ref{item:admissible-square product}, it suffices to show that $-\mu=p\in\units\Q_{v_2}/2$, for which we compute
        \begin{align}\label{eqn:n+1-long-computation}
            (-\mu p,pv_2)_{v_2} \overset{(\ref{eqn:n+1-v2})}= (-\mu p, pv_2)_p =(\mu,pv_2)_p(-p,p)_p(-p,v_2)_p\overset{\t{\ref{item:n-prepared-square at p} and (\ref{eqn:n+1-v2})}}=+1
            ,
        \end{align}
        where the first equality, when $n=2$, additionally uses the fact that $p_2\equiv3\bmod4$.
        Finally, we remark that $1-m^2,2\equiv\square\bmod v_2$. The first of these holds by the below computations:
        \begin{align}
            (1-m^2,\Gamma_{n+1})_{v_2} &\overset{(\ref{eqn:n+1-v2})}=(1-m^2,\Gamma_{n+1})_p=(1-m^2,v_2)_p(1-m^2,p)_p\overset{\t{(\ref{eqn:n+1-v2}) + \ref{item:admissible-p primes}}}=+1 \tand\\
            (2,\Gamma_{n+1})_{v_2} &\overset{(\ref{eqn:n+1-v2})}=(2,\Gamma_{n+1})_2(2,\Gamma_{n+1})_p = (2,\Gamma_{n+1})_p\overset{\t{\ref{item:admissible-p primes}}}=+1
            .
        \end{align}
        This observation will allow us to take $q_{n+1}\coloneqq v_2$ in \ref{item:n-prepared-q primes}.

        \item Set $v_3=p$ and choose $d_3\in\Z$ to be a non-square unit at $p$. We claim that
        \begin{equation}\label{eqn:n+1-Sel L3}
            \Sel\mcL^3=\Sel\mcL^2\oplus\angles{\gamma_{n+1}}.
        \end{equation}
        First note that, because $d_3\neq\square\in\Z_p$, Lemma \ref{lem:Selmer-condition-at-special-p} implies that $\dim(\im\delta_{v_3}\cap\im\delta_{v_3}^{d_3})=1$. Thus, by Proposition \ref{prop:Selmer-change-bad}, to establish (\ref{eqn:n+1-Sel L3}), it suffices to show that $\gamma_{n+1}\in\Sel(\mcL^3)$. Since $\gamma_{n+1}\in\Sel(\mcL^2)^{(p)}=\Sel(\mcL^3)^{(p)}$, it suffices to show that $\gamma_{n+1}=(\Gamma_{n+1},\Gamma_{n+1})$ satisfies the $\mcL^3$-Selmer condition at $p=v_3$. This follows from Lemma \ref{lem:Selmer-condition-at-special-p}, noting that $d_3\mu\equiv\square\bmod p$ (by \ref{item:n-prepared-square at p}).

        \item Suppose that $\Sel\mcL^2\supsetneq\Sel\mcL^1$ (else, skip this step). In this case, e.g. as a consequence of Proposition \ref{prop:Selmer-change}, one must have $\dim\Sel\mcL^2=2+\dim\Sel\mcL^1$ and even that we may write $\Sel\mcL^2=\Sel\mcL^1\oplus\angles{\theta,\theta'}$ for some $\theta,\theta'\in(\units\Q/2)^2$ with $\oploc_{v_2}\angles{\theta,\theta'}=\im\delta^{d_2}_{v_2}$. Using the description of $\im\delta^{d_2}_{v_2}$ given by Lemma \ref{lem:twisted-kummer-image}, we may and do choose $\theta,\theta'$ so that
        \begin{align}\label{eqn:theta-theta'-at v2}
            \oploc_{v_2}\theta &=(d_2(\lambda-\mu),\lambda(\lambda-\mu))\in(\units\Q_{v_2}/2)^2 &&\tand\\
            \oploc_{v_2}\theta' &=(\mu(\mu-\lambda),(\mu-\lambda)d_2)\in(\units\Q_{v_2}/2)^2
            ,
        \end{align}
        Furthermore, it follows from Lemmas \ref{lem:Selmer-condition-at-special-p}, \ref{lem:Selmer-condition-at-special-ell}, and \ref{lem:Selmer-condition-at-infty} combined with \ref{item:n-prepared-gammas property} that, by modifying $\theta,\theta'$ by an appropriate combination of $\gamma_i$'s, we may and do assume that $\theta,\theta'$ both have positive coordinates and are unramified at $p_1,\dots,p_{n-1}$. 
        Since $\theta,\theta'\in\Sel\mcL^3$ as well (by (\ref{eqn:n+1-Sel L3})), we know from construction that
        \[\oploc_{v_3}(\theta),\oploc_{v_3}(\theta')\in\mcL^2_{v_3}\cap\mcL^3_{v_3}=\im\delta_p\cap\im\delta_p^{d_3}\overset{\t{Lemma \ref{lem:Selmer-condition-at-special-p}}}=\uhom^1(\Q_p,E[2]),\]
        so $\theta,\theta'$ are unramified at $p=p_n=v_3$. Write $\theta=(\theta_1,\theta_2),\theta'=(\theta_1',\theta_2')\in(\units\Q/2)^2$; it follows from ramification considerations that $\Q(\sqrt{\theta_1},\sqrt{\theta_2},\sqrt{\theta_1'},\sqrt{\theta_2'})$
        is linearly disjoint from 
        $\Q(\sqrt{\Gamma_1},\dots,\sqrt{\Gamma_n})$. Additionally, 
        (\ref{eqn:theta-theta'-at v2}) shows that $\theta_2,\theta_1'$ are unramified at $v_2$, while $\theta_1,\theta_2'$ are ramified at $v_2$. It follows that we may apply Chebotarev density to produce a positive prime $v_4$ satisfying
        \begin{itemize}
            \item $\Gamma_i\equiv\square\bmod v_4$ for all $i=1,\dots,n$;
            \item $\theta_1,\theta_2'\not\equiv\square\bmod v_4$; and
            \item $\theta_2,\theta_1'\equiv\square\bmod v_4$.
        \end{itemize}
        Set $d_4=v_4$. The first bullet above implies that each $\gamma_i$ satisfies the $\mcL^4$-Selmer condition at $v_4$, so $\Sel\mcL^4\supset\Sel\mcL^1\oplus\angles{\gamma_{n+1}}$. The second and third bullets imply that $\dim\oploc_{v_4}\Sel\mcL^3=2$, so Proposition \ref{prop:Selmer-change} tells us that $\dim\Sel\mcL^4=\dim\Sel\mcL^3-2=\dim\Sel\mcL^1+1$ and so
        \begin{equation}
            \Sel\mcL^4=\Sel\mcL^1\oplus\angles{\gamma_{n+1}}=\angles{\gamma_1,\dots,\gamma_{n+1}}.
        \end{equation}

        \item Finally, set $r\coloneqq4$ if $\Sel\mcL^2=\Sel\mcL^1$ and $r\coloneqq5$ otherwise. Thus, $\Sel\mcL^{r-1}=\angles{\gamma_1,\dots,\gamma_{n+1}}$. Invoking Dirichlet, we obtain a positive prime $v_r\not\in S(\mcL^{r-1})$ such that $t=d_r\coloneqq(v_1v_2\dots v_r)/v_3$ satisfies $t\equiv1\bmod8v$ for all $v\in S(E)-\{p\}$, but (recall $p=v_3$)
        \[t\equiv d_3\bmod{v_3},\tand\frac t{v_i}=\prod_{\substack{j=1\\j\neq 3,i}}^rv_j\equiv\frac{d_i}{v_i}\bmod v_i\tforall i\in\{1,\dots,r-1\}-\{3\}.\]
        It follows from Lemma \ref{lem:last-step-+2} that $\Sel_2(E^t)=\Sel\mcL^r=\Sel\mcL^{r-1}\oplus\kumtwist t$; hence, $E^t$ -- equivalently, $(t\mu,t\lambda)$ -- is $(n+1)$-prepared in the sense of Definition \ref{defn:n-prepared}.
    \qedhere        
    \end{enumerate}
\end{proof}

\begin{lemma}\label{lem:n=4}
    Set $n\coloneqq3$. Let $(\mu,\lambda)$ be an $n$-prepared admissible pair. Then, there exists some $t$ such that $E^t$ is $(n+1)$-prepared.
\end{lemma}
\begin{proof}
    Adopt all of the notation of Definition \ref{defn:n-prepared}. Set $\l\coloneqq\l_1$. By \ref{item:n-prepared-gammas}, $\Sel_2(E)=\kumtwist{}\oplus\angles{\gamma_1,\dots,\gamma_n}$, where $\gamma_1=(\Gamma_1,1)$, $\gamma_2=(\Gamma_2,\Gamma_2)$, $\gamma_3=(\Gamma_3,\Gamma_3)$, and these data satisfy \ref{item:n-prepared-gammas property}.
    We construct $t=d_r$ as the end result of a sequence $(v_1,d_1),\dots,(v_r,d_r)$ as in Section \ref{sect:quad-twist}. The argument is largely the same as the proof of Lemma \ref{lem:n+1}, so we only indicate the necessary adjustments.
    \begin{enumerate}
        \item Apply Dirichlet to produce a positive prime $v_1\not\in S(E)=S(\mcL^0)$ such that
        \begin{equation}
            v_1\equiv1\bmod8v\tforall v\in S(E)-\{\varpi,\infty\},\tbut v_1\not\equiv\square\bmod\varpi.
        \end{equation}
        Set $d_1\coloneqq v_1$. As in step \bp1 of Lemma \ref{lem:n+1}, one deduces that $\Sel\mcL^1=\angles{\gamma_1,\dots,\gamma_n}$.
        \item From Dirichlet, obtain a positive prime $v_2\not\in S(\mcL^1)$ such that
        \[\l v_2\equiv1\bmod8v\tforall v\in S(\mcL^1)-\{\l,\infty\},\tbut v_2\not\equiv\square\bmod\l.\]
        Set $d_2\coloneqq v_2$, $\Gamma_{n+1}\coloneqq\l v_2$, and $\gamma_{n+1}\coloneqq(1,\Gamma_{n+1})\in(\units\Q/2)^2$. As in step \bp2 of Lemma \ref{lem:n+1}, one checks that $\Sel\mcL^1\subset\Sel\mcL^2$ and that $\gamma_{n+1}\in\Sel(\mcL^2)^{(\l)}$. For this, one checks that $\oploc_{v_2}(\gamma_{n+1})=(\mu(\mu-\lambda),(\mu-\lambda)d_2)\overset{(\ref{eqn:kum-image})}=\delta_{v_2}^{d_2}((\mu,0))\in(\units\Q_{v_2}/2)^2$. Furthermore, here one computes that
        \[(1-m^2,\Gamma_{n+1})_{v_2}=(1-m^2,\Gamma_{n+1})_\l=(1-m^2,\l)_\l(1-m^2,v_2)_\l=+1,\]
        and similarly that $(2,\Gamma_{n+1})_{v_2}=(2,\Gamma_{n+1})_\l=+1$ (with last equality using \ref{item:admissible-l primes}),
        so $1-m^2,2\equiv\square\bmod v_2$, in line with \ref{item:n-prepared-q primes} for $q_4\coloneqq v_2$.
        \item Set $v_3=\l$ and choose $d_3$ to be an integral non-square at $\l$. As in step \bp3 of Lemma \ref{lem:n+1}, one deduces that $\Sel\mcL^3=\Sel\mcL^2\oplus\angles{\gamma_{n+1}}$, making use of Lemma \ref{lem:Selmer-condition-at-special-ell} in place of Lemma \ref{lem:Selmer-condition-at-special-p}.
        \item Suppose that $\Sel\mcL^2\supsetneq\Sel\mcL^1$ (else, skip this step). Arguing as in step \bp4 of Lemma \ref{lem:n+1}, one produces a prime $v_4=d_4$ such that $\Sel\mcL^4=\angles{\gamma_1,\dots,\gamma_n,\gamma_{n+1}}$.
        \item Finally, set $r\coloneqq4$ if $\Sel\mcL^2=\Sel\mcL^1$ and $r\coloneqq5$ otherwise. Dirichlet gives us a positive prime $v_r\not\in S(\mcL^{r-1})$ such that $t=d_r\coloneqq(v_1v_2\dots v_r)/v_3$ satisfies $t\equiv1\bmod8v$ for all $v\in S(E)-\{\l\}$, $t\equiv d_3\bmod v_3$, and $t/v_i\equiv d_i/v_i\bmod v_i$ for all $i\in\{1,\dots,r-1\}-\{3\}$. Lemma \ref{lem:last-step-+2} gives that $\Sel_2(E^t)=\Sel\mcL^r=\Sel\mcL^{r-1}\oplus\kumtwist{t}$, so $E^t$ is $(n+1)$-prepared, as desired.
        \qedhere
    \end{enumerate}
\end{proof}
\begin{lemma}\label{lem:n=1}
    Let $(\mu,\lambda)$ be a $0$-prepared admissible pair. Then, there exists some negative $t\equiv1\bmod8$ such that $E^t$ is $1$-prepared.
\end{lemma}
\begin{proof}
    Adopt all of the notation of Definition \ref{defn:n-prepared}. By \ref{item:n-prepared-gammas}, $\Sel_2(E)=\kumtwist{}$.
    We construct $t=d_r$ as the end result of a sequence $(v_1,d_1),\dots,(v_r,d_r)$ as in Section \ref{sect:quad-twist}. We proceed via a modification of the argument used to prove Lemma \ref{lem:n+1}.
    \begin{enumerate}
        \item Apply Dirichlet to produce a positive prime $v_1\not\in S(E)$ such that
        \begin{equation}\label{eqn:n=1-v1}
            v_1\equiv1\bmod 8v\tforall v\in S(E)-\{\varpi,\infty\}\tbut v_1\not\equiv\square\bmod\varpi.
        \end{equation}
        Set $d_1\coloneqq v_1$. As in step \bp1 of Lemma \ref{lem:n+1}, one deduces $\Sel\mcL^1=0$.

        \item Using Dirichlet, we fix a positive prime $v_2\not\in S(\mcL^1)$ such that
        \begin{equation}\label{eqn:n=1-v2}
            -v_2\equiv1\bmod 8v\tforall v\in S(\mcL^1)-\{\infty\}.
        \end{equation}
        Set $d_2\coloneqq{-v_2}$, $\Gamma_1\coloneqq-v_2$, and $\gamma_1\coloneqq(\Gamma_1,1)\in(\units\Q/2)^2$. We claim that $\gamma_1\in\Sel(\mcL^2)^{(\infty)}$, i.e. that $\gamma_1$ satisfies the $\mcL^2$--Selmer condition at $v_2$. We first compute that $(1-m^2,\Gamma_1)_{v_2}=(1-m^2,\Gamma_1)_\infty=-1$, so
        \begin{equation}\label{eqn:n=1-v2-squares}
            1-m^2\not\equiv\square\bmod v_2.\quad\t{Similarly, }\mu,\lambda\equiv\square\bmod v_2,\twhile{-1}\not\equiv\square\bmod v_2.
        \end{equation}
        From this, we see $\oploc_{v_2}(\gamma_1)=(v_2,1)\overset{(\ref{eqn:n=1-v2-squares})}=(d_2(\lambda-\mu),\lambda(\lambda-\mu))\overset{(\ref{eqn:kum-image})}=\delta^{d_2}_{v_2}((\lambda,0))\in(\units\Q_{v_2}/2)^2$.
        
        Finally, we claim that $\Sel\mcL^2=0$. Suppose otherwise. Then, by Proposition \ref{prop:Selmer-change} and Lemma \ref{lem:twisted-kummer-image}, we can write $\Sel\mcL^2=\angles{\theta,\theta'}$ where (using (\ref{eqn:n=1-v2-squares}) below)
        \begin{align}\label{eqn:m+1-thetas}
            \oploc_{v_2}(\theta)&=((\lambda-\mu)d_2,\lambda(\lambda-\mu))=(-v_2,+1)\in(\units\Q_{v_2}/2)^2 &&\tand\\
            \oploc_{v_2}(\theta')&=(\mu(\mu-\lambda),(\mu-\lambda)d_2)=(-1,v_2)\in(\units\Q_{v_2}/2)^2
            ,
        \end{align}
        Write $\theta=(\theta_1,\theta_2)$ and $\theta'=(\theta_1',\theta_2')$ in $(\units\Q/2)^2$. Then, (\ref{eqn:n=1-v2}) and (\ref{eqn:m+1-thetas}) together imply that $(-v_2\theta_1,\theta_2),(\theta_1',-\theta_2'v_2)\in\Sel(\mcL^1)^{(\infty)}$. Lemma \ref{lem:naive-relaxed-dimension-bound} yields $\dim\Sel(\mcL^1)^{(\infty)}\le\dim\Sel(\mcL^1)+\frac12\dim\hom^1(\R,E[2])=1$, so one of the following must hold:
        \[(-v_2\theta_1,\theta_2)=(\theta_1',-\theta_2'v_2)\tcomma (-v_2\theta_1,\theta_2)=(1,1),\tor(\theta_1',-\theta_2'v_2)=(1,1).\]
        The first of these is impossible by comparison with (\ref{eqn:m+1-thetas}); the latter two are impossible because the $\mcL^2$-Selmer condition at $\infty$ (which is spanned by $\delta_\infty((0,0))=(-\mu,-\lambda)=(-1,-1)\in(\units\R/2)^2$) forces $\sign\theta_1=\sign\theta_2$ and $\sign\theta_1'=\sign\theta_2'$. This gives a contradiction, so $\Sel\mcL^2=0$ after all.

        \item Set $v_3=\infty$ and $d_3=-1$. As in step \bp3 of Lemma \ref{lem:n+1}, one deduces that $\Sel\mcL^3=\Sel\mcL^2\oplus\angles{\gamma_1}=\angles{\gamma_1}$ making use of Lemma \ref{lem:Selmer-condition-at-infty} in place of Lemma \ref{lem:Selmer-condition-at-special-p}.

        \item Finally, invoke Dirichlet to produce a positive prime $v_4\not\in S(\mcL^3)$ such that $t=d_4\coloneqq{-v_1v_2v_4}$ satisfies $t\equiv1\bmod8v$ for all $v\in S(E)-\{\infty\}$ and $t/v_i\equiv d_i/v_i\bmod v_i$ for $i\in\{1,2\}$. Lemma \ref{lem:last-step-+2} gives that $\Sel_2(E^t)=\Sel\mcL^r=\Sel\mcL^{r-1}\oplus\kumtwist t$, so $E^t$ is $1$-prepared as desired.
        \qedhere
    \end{enumerate}
\end{proof}

\begin{cor}\label{cor:4-prepared}
    Let $(\mu,\lambda)$ be an admissible pair. Then, there exists some $d$ such that $E^d$ is $4$-prepared.
\end{cor}
\begin{proof}
    Apply Lemmas \ref{lem:n=0}, \ref{lem:n=1}, \ref{lem:n+1}, \ref{lem:n+1}, and \ref{lem:n=4}, in that order.
\end{proof}

\section{Additive combinatorics twisting}\label{sect:AC-twist}
In this final section, we construct our elliptic curves of rank $2$. Let $E=E_{\mu,\lambda}\colon y^2=x(x-\mu)(x-\lambda)$ for some pair $(\mu,\lambda)\in\Z^2$ with $\mu\lambda=\eta^2$. Let $G(u)\in\Q[u]$ be the polynomial, depending on $\mu,\lambda$, from (\ref{eqn: G(u)=QLLLL}). As in Proposition \ref{prop:RS-family}, set $\alpha=\mu\lambda/(\mu+\lambda)$. By Proposition \ref{prop:RS-family}, the quadratic twist family $E_{\mu,\lambda}^{G(u)}$ has generic rank $\ge2$. We will show that, for suitable choice of $u$, it has $2$-Selmer rank equal to $4$ (and so Mordell--Weil rank $\le2$). Observe that, for $u=X/Y$, twisting by $G(u)=G(X/Y)$ is equivalent to twisting by
\begin{equation}\label{eqn:G-denominators cleared}
    \frac{Y^6\alpha^4(\mu+\lambda)^8}{\mu^6\eta^2}G\pfrac XY=(\mu+\lambda)\p{\pfrac{\mu+\lambda}\mu^2X^2+\lambda^2Y^2}\p{\pfrac{\mu+\lambda}\mu X\pm\eta Y}\p{\pfrac{\mu+\lambda}\mu X\pm\pfrac{\eta\lambda}\mu Y},
\end{equation}
and, ignoring the $\mu+\lambda$ factor for now, the RHS above is, up to squares, a product of one quadratic form and 4 linear forms in $X,Y$. Theorem \ref{thm: prime count} will allow us to force the quadratic form and linear forms above to take on prime values -- say, $\pi_1,\pi_2,\pi_3,\pi_4,\pi_5$ -- with pre-specified congruence conditions. If $E$ is $4$-prepared in the sense of Definition \ref{defn:n-prepared}, then these congruence conditions may be chosen to ensure that $E^{G(u)}\cong E^{(\mu+\lambda)\pi_1\pi_2\pi_3\pi_4\pi_5}$ has $2$-Selmer rank $\le4$ (see Proposition \ref{prop:rank-2-crit}). Combined with the two points produced by Proposition \ref{prop:RS-family}, this will show that $\rank E^{G(u)}(\Q)=2$ for such appropriately chosen $u\in\Q$.

\subsection{Criterion for rank 2}
\begin{set}
    Fix a $4$-prepared admissible pair $(\mu,\lambda)$ -- in the sense of Definition \ref{defn:n-prepared} -- and set $E=E_{\mu,\lambda}\colon y^2=x(x-\mu)(x-\lambda)$. Adopt Notation \ref{notn:admissible-places}, and fix an odd prime $\varpi$ dividng $\mu$, as in \ref{item:n-prepared-odd factor}. Finally, by \ref{item:admissible-square product}, we write $\lambda/\mu=m^2$ for some $m\in\Z_{>0}$, so $\mu\lambda=\eta^2$ for $\eta\coloneqq m\mu$. Also recall that, by \ref{item:admissible-m-at-2}, $m\equiv2\bmod8$ and $\mu$ is odd.
\end{set}
\begin{rem}\label{rem:G-factors-simplification}
    The various factors appearing in (\ref{eqn:G-denominators cleared}) may be rewritten as
    \[\begin{matrix}
        \pfrac{\mu+\lambda}\mu^2X^2+\lambda^2Y^2 = (1+m^2)^2X^2+m^4\mu^2Y^2\\ 
        \pfrac{\mu+\lambda}\mu X\pm\eta Y=(1+m^2)X\pm m\mu Y & \pfrac{\mu+\lambda}\mu X\pm\pfrac{\eta\lambda}\mu Y=(1+m^2)X\pm m^3\mu Y
    \end{matrix}\]
    In particular, they are quadratic or linear forms (in $X,Y$) with \important{integral} coefficients.
\end{rem}
\begin{lemma}\label{lem:modified linear+quadratic forms}
    Fix $a,b\in\Z$ such that $\mu\mid(a+b)$ and $(1+m^2)\mid(a-b)$.
    Consider variables $X_0,Y_0$ and set
    \[X\coloneqq \mu^2(1+m^2)mX_0+m\mu a\tand Y\coloneqq\mu(1+m^2)^2Y_0+(1+m^2)b.\]
    Set
    \begin{align}
        \pi_1 = \pi_1(X_0,Y_0) &\coloneqq (1+m^2)(X_0+Y_0) + \frac{a+b}\mu &&\in\Z[X_0,Y_0] \\
        \pi_2 = \pi_2(X_0,Y_0) &\coloneqq \mu(X_0-Y_0) + \frac{a-b}{1+m^2} &&\in\Z[X_0,Y_0] \\
        \pi_3 = \pi_3(X_0,Y_0) &\coloneqq \mu(X_0+m^2Y_0) + \frac{a+m^2b}{1+m^2} &&\in\Z[X_0,Y_0] \\
        \pi_4 = \pi_4(X_0,Y_0) &\coloneqq \mu(1+m^2)(X_0-m^2Y_0) + \p{a-m^2b} &&\in\Z[X_0,Y_0] \\
        \pi_5 = \pi_5(X_0,Y_0) &\coloneqq \mu^2(1+m^2)(X_0^2+m^2Y_0^2) + 2\mu(X_0a+m^2Y_0b)+\frac{a^2+m^2b^2}{1+m^2} &&\in\Z[X_0,Y_0]
        .
    \end{align}
    Finally, let $G(u)$ be the polynomial from (\ref{eqn: G(u)=QLLLL}). Then,
    \[G\pfrac XY=N^2\pi_1\pi_2\pi_3\pi_4\pi_5\tforsome N\in\Q[Y].\]
\end{lemma}
\begin{proof}
    This is ultimately a matter of algebra combined with Remark \ref{rem:G-factors-simplification} and (\ref{eqn:G-denominators cleared}). First we remark that, by assumption, $1+m^2$ divides both
    \[(a-b)+b(1+m^2)=a+m^2b\tand(a-b)(a+b)+b^2(1+m^2)=a^2+m^2b^2,\]
    so $\pi_1,\dots,\pi_5$ indeed all have integral coefficients in $X_0,Y_0$.
    By (\ref{eqn:G-denominators cleared}) and Remark \ref{rem:G-factors-simplification}, $G(X/Y)$ is, up to squares, equal to $\mu+\lambda$ (i.e. $\mu(1+m^2)$) times the product of the below 5 quantities:
    \begin{alignat}{5}
        (1+m^2)X &+ m\mu Y &&= \mu^2(1+m^2)^2 m(X_0+Y_0) &&+ \mu(1+m^2)m(a+b) &&= \mu^2(1+m^2)m\cdot\pi_1, \\
        (1+m^2) X &- m\mu Y &&= \mu^2(1+m^2)^2 m (X_0-Y_0) &&+ \mu (1+m^2)m(a-b) &&= \mu(1+m^2)^2m\cdot\pi_2, \\
        (1+m^2)X &+ m^3 \mu Y &&= \mu^2(1+m^2)^2 m(X_0+m^2 Y_0) &&+ \mu(1+m^2)m(a+m^2 b) &&= \mu(1+m^2)^2m\cdot\pi_3, \\
        (1+m^2)X &- m^3 \mu Y &&= \mu^2(1+m^2)^2 m(X_0-m^2 Y_0) &&+ \mu(1+m^2)m(a -m^2 b) &&= \mu(1+m^2)m\cdot\pi_4,
    \end{alignat}
    \begin{align}
        (1+m^2)^2X^2 + m^4\mu^2 Y^2 
        &= (1+m^2)^2 (\mu^4(1+m^2)^2 m^2 X_0^2 + 2\mu^3(1+m^2)m^2 X_0 a  + \mu^2 m^2 a^2) \\
        &\phantom{=(1+} + m^4 \mu^2 (\mu^2(1+m^2)^4 Y_0^2 + 2\mu(1+m^2)^3 Y_0 b + (1+m^2)^2 b^2) \\
        &= \mu^4 (1+m^2)^4 m^2 ( X_0^2 + m^2 Y_0^2) + 
            2\mu^3 (1+m^2)^3 m^2 (X_0a+m^2 Y_0 b) \\
        &\phantom{=(1+} + \mu^2 (1+m^2)^2 m^2 (a^2+m^2b^2) \\
        &= \mu^2(1+m^2)^3m^2\cdot\pi_5
        .
    \end{align}
    Comparing with (\ref{eqn:G-denominators cleared}), we have
    \begin{align}
        G(X/Y) &= n^2\mu(1+m^2)\cdot\sq{(1+m^2)^2X^2 + m^4\mu^2 Y^2}\prod_{\substack{\eps=\pm1\\e\in\{1,3\}}}\sq{(1+m^2)X+\eps m^e\mu Y} \\
        &= {n^2\mu(1+m^2)}\cdot{\mu^7(1+m^2)^9m^6}\cdot\pi_1\pi_2\pi_3\pi_4\pi_5 \\
        &=\p{n\mu^4(1+m^2)^5m^3}^2\cdot\pi_1\pi_2\pi_3\pi_4\pi_5
        ,\quad\tforsome n\in\Q[Y].
        \qedhere
    \end{align}
\end{proof}
\begin{prop}\label{prop:rank-2-crit}
    Adopt all the notation of Definitions \ref{defn:n-prepared} and \ref{defn:admissible}. Set
    \[N\coloneqq\prod_{\substack{p\in S(E)\\p\neq2,\infty}}p\tcomma N_\mu\coloneqq\prod_{\substack{q\mid\mu\\q\neq\varpi}}q\tcomma N_1\coloneqq\l_1p_1p_2N_\mu\tcomma N_2\coloneqq\l_1p_2q_2q_3q_4,\tand N_3\coloneqq p_2q_3q_4.\]
    For $i=1,\dots,5$, let $\pi_i(X_0,Y_0)\in\Z[X_0,Y_0]$ be as in Lemma \ref{lem:modified linear+quadratic forms}. Suppose that $\xi,\upsilon\in\Z$ are nonzero integers such that the $\pi_i=\pi_i(\xi,\upsilon)\in\Z$ satisfy the following.
    \begin{enumerate}
        \item $\pi_1<0$ is a negative prime,
        \[\phantom{(1-m^2)}\pi_1\equiv\square\bmod N_1\tcomma \pi_1\equiv1\bmod4,\tand \pi_1\not\equiv\square\bmod\varpi.\phantom{\quad\t{Also},}\]
        
        \item $\pi_2>0$ is a positive prime,
        \begin{equation}\begin{split}
            \pi_2&\equiv\square\bmod N_2\tcomma \pi_2\equiv1\bmod4,\tand \pi_2\not\equiv\square\bmod p_1. \quad\t{Also},\\
            \pi_2(1-m^2)&\not\equiv\square\bmod\varpi,\tbut \pi_2(1-m^2)\equiv\square\bmod v\tforall v\mid N_\mu.
        \end{split}\end{equation}
        Note that the two lines are consistent (for $v=q_2,q_3,q_4$) because of \ref{item:n-prepared-q primes}.
        
        \item $\pi_3$ is a prime,
        \[\phantom{(1-m^2)}\pi_3\equiv\square\bmod N_3\tcomma \pi_3\equiv3\bmod4,\tand \pi_3\not\equiv\square\bmod\l_1.\phantom{\quad\t{Also},}\]
        
        \item $\pi_4$ is a prime.

        \item $\pi_5$ is a prime such that
        \[t\coloneqq\pi_1\pi_2\pi_3\pi_4\pi_5\equiv\square\bmod 8N\tand t>0.\]
    \end{enumerate}
    Then, $t=G(u)\cdot n^2$ for some $u,n\in\Q$ (so $E^t\simeq E^{G(u)}$), and $\dim\Sel_2(E^t)\le4$.
\end{prop}
\begin{proof}
    First, the equality $\pi_1\pi_2\pi_3\pi_4\pi_5=G(u)\cdot n^2$ holds, for suitable choice of $u,n\in\Q$, by Lemma \ref{lem:modified linear+quadratic forms}. For $k=1,\dots,5$, consider the Selmer structures
    \[\mcL^k\coloneqq\mcL\p{E,((\pi_1,t),\dots,(\pi_k,t))}\quad\t{ as in }(\ref{eqn:Li-def}).\]
    By assumption, $(\mu,\lambda)$ is $4$-prepared (Definition \ref{defn:n-prepared}), so we may and do choose $\Gamma_1,\dots,\Gamma_4\in\units\Q$ such that
    \begin{equation}\begin{split}
        \Sel_2&(E)= \kumtwist{}\oplus\angles{\gamma_1,\gamma_2,\gamma_3,\gamma_4},\twhere \\
        &\gamma_1\coloneqq(\Gamma_1,1)\tcomma\gamma_2\coloneqq(\Gamma_2,\Gamma_2)\tcomma\gamma_3\coloneqq(\Gamma_3,\Gamma_3),\tand\gamma_4\coloneqq(1,\Gamma_4)\tin(\units\Q/2)^2.
    \end{split}
    \end{equation}
    Additionally, the above data satisfies \ref{item:n-prepared-gammas property} and \ref{item:n-prepared-q primes}.
    We compute the Selmer groups $\Sel\mcL^k$ in a series of claims:
    \begin{romanize}
        \item Claim: $\Sel\mcL^1=\angles{(\lambda-\mu,\lambda(\lambda-\mu)),\gamma_2,\gamma_3,\gamma_4}\subset(\units\Q/2)^2$.
        First recall, from (\ref{eqn:kum-image}), that
        \begin{equation}
            \kumtwist{} = \big\langle{(-\mu,-\lambda),(\lambda-\mu,\lambda(\lambda-\mu))}\big\rangle\subset(\units\Q/2)^2\overset{\t{Convention \ref{conv:H1(E[2])}}}\simeq\hom^1(\Q,E[2]).
        \end{equation}
        Using assumption \bp1 (and $\mu$ being odd and negative), we see that $-1\not\equiv\square\bmod\pi_1$, but
        \begin{alignat}{5}
            (\mu,\pi_1)_{\pi_1} &=(\mu,\pi_1)_\infty\prod_{v\mid\mu}(\mu,\pi_1)_v &&=(\mu,\pi_1)_\infty(\mu,\pi_1)_\varpi &&=+1 &&\tand\\
            (\lambda-\mu,\pi_1)_{\pi_1} &=(\lambda-\mu,\pi_1)_\infty\prod_{v\mid\lambda-\mu}(\lambda-\mu,\pi_1)_v &&=(\lambda-\mu,\pi_1)_\infty(\lambda-\mu,\pi_1)_\varpi &&= +1
            ,
        \end{alignat}
        so $\mu,\lambda-\mu\equiv\square\bmod\pi_1$, as is $\lambda=\mu m^2$. Hence,
        \[\oploc_{\pi_1}(-\mu,-\lambda)=(-1,-1)\tand\oploc_{\pi_1}(\lambda-\mu,\lambda(\lambda-\mu))=(+1,+1)\tin(\units\Q_{\pi_1}/2)^2.\]
        Similarly computations show that $\Gamma_1\not\equiv\square\bmod\pi_1$, but $\Gamma_i\equiv\square\bmod\pi_1$ for $i>1$. It follows that $\oploc_{\pi_1}\Sel(\mcL^0)$ is 2-dimensional, generated by the images of $(-\mu,-\lambda)$ and $\gamma_1=(\Gamma_1,1)$, and that $\angles{(\lambda-\mu,\lambda(\lambda-\mu)),\gamma_2,\gamma_3,\gamma_4}\subset\Sel\mcL^1$. Proposition \ref{prop:Selmer-change} implies that $\dim\Sel\mcL^1=\dim\Sel\mcL^0-2=4$, so our first claim holds.

        \item Claim: $\Sel\mcL^2=\angles{\gamma_3,\gamma_4}$.

        We compute
        \begin{alignat}{5}
            (\Gamma_2,\pi_2)_{\pi_2} &= \prod_{v\neq\pi_2}(\Gamma_2,\pi_2)_v &&\overset{\t{\ref{item:n-prepared-q primes}+\ref{item:n-prepared-gammas property}}}= (\Gamma_2,\pi_2)_{p_1}(\Gamma_2,\pi_2)_{q_2} \overset{\bp2}={-1}
            ,
        \end{alignat}
        so $\Gamma_2\not\equiv\square\bmod\pi_2$. We also compute, recalling that $\gcd(\mu,1-m^2)=1$ by \ref{item:admissible-coprimality},
        \begin{align}
            (\mu,\pi_2)_{\pi_2} &= \prod_{v\mid\mu}(\mu,\pi_2)_v = -\prod_{v\mid\mu}(\mu,1-m^2)_v &&\t{by \bp2}\\
            &= -(\mu,1-m^2)_\infty(\mu,1-m^2)_2\prod_{v\mid(1-m^2)}(\mu,1-m^2)_v &&\t{by Hilbert reciprocity}\\
            &= -(-1)(+1)(-1) = -1 &&\t{by \ref{item:n-prepared-mu negative},\ref{item:n-prepared-mu at 2},\ref{item:n-prepared-hilbert-cond}}
            ,
        \end{align}
        so $\mu\not\equiv\square\bmod\pi_2$. It follows that $\lambda\not\equiv\square\bmod\pi_2$ and so that $(\lambda-\mu,\lambda(\lambda-\mu))$ and $\gamma_2=(\Gamma_2,\Gamma_2)$ generate a 2-dimensional subspace of $(\units\Q_{\pi_2}/2)^2\simeq\hom^1(\Q_{\pi_2},E[2])$. Proposition \ref{prop:Selmer-change} gives $\dim\Sel\mcL^2=\dim\Sel\mcL^1-2=2$. At the same time, straightforward Hilbert symbol computations, using \bp2, \ref{item:n-prepared-gammas property}, and \ref{item:n-prepared-q primes}, show that $\Gamma_3,\Gamma_4\equiv\square\bmod\pi_2$, and so $\angles{\gamma_3,\gamma_4}\subset\Sel\mcL^2$; comparing dimensions, this must be an equality.

        \item Claim: $\Sel\mcL^3=0$.

        Using \bp3, \ref{item:n-prepared-gammas property}, and \ref{item:n-prepared-q primes}, one quickly verifies that $\Gamma_3,\Gamma_4\not\equiv\square\bmod\pi_3$, so $\gamma_3=(\Gamma_3,\Gamma_3)$ and $\gamma_4=(1,\Gamma_4)$ generate a 2-dimensional subspace of $(\units\Q_{\pi_2}/2)^2\simeq\hom^1(\Q_{\pi_3},E[2])$. Hence, Proposition \ref{prop:Selmer-change} entails that $\dim\Sel\mcL^3=\dim\Sel\mcL^2-2=0$.

        \item Claim: $\dim\Sel\mcL^4\le2$ and $\dim\Sel_2(E^t)=\dim\Sel\mcL^5\le4$.
    
        By \bp 5, Lemma \ref{lem:last-step-+2} applies to $\mcL^5$, so
        \[\dim\Sel_2(E^t)=\dim\Sel\mcL^5=\dim\Sel\mcL^4+2\overset{\t{Proposition \ref{prop:Selmer-change}}}\le\dim\Sel\mcL^3+2+2=4.\qedhere\]
    \end{romanize}
\end{proof}

\subsection{Proof of Theorem \ref{thm:rank2}}\label{sect:rank2-proof}

Fix an admissible pair $(\mu,\lambda)$, and set $E=E_{\mu,\lambda}\colon y^2=x(x-\mu)(x-\lambda)$. We will prove Theorem \ref{thm:rank2} for $E$, i.e. construct infinitely many $d$ so that $\rank E^d(\Q)=2$. First, by Corollary \ref{cor:4-prepared}, by suitably twisting $E$, we may assume that $(\mu,\lambda)$ is $4$-prepared in addition to being admissible. Ultimately, we wish to apply Proposition \ref{prop:rank-2-crit}, so we show below (see Lemma \ref{lem:local conditions for linear transformation}) that its conditions are achievable.

\begin{lemma}\label{lem:Lang-Weil-curve-trick}
    Fix an odd prime power $q$. Fix $k\ge2$ along with $s_0,\dots,s_k,t_0,\dots,t_k,\eps_0,\dots,\eps_k\in\units\F_q$. Assume that the (zero-dimensional) varieties $\{s_0 X^2+t_0Z^2\},\{s_1X+t_1Z\},\dots\{s_kX+t_kZ\}\subset\P^1_{\F_q}$ (with coordinates $X,Z$) are pairwise disjoint. Consider the affine variety $V\subset\A^{k+3}_{\F_q}$ given by
    \[V\colon\left\{\begin{aligned}
        s_0x^2+t_0z^2&=\eps_0y_0^2\\
        s_1x\phantom{^2}+t_1z\phantom{^2}&=\eps_1y_1^2\\
        &\;\;\vdots\\
        s_kx\phantom{^2}+t_kz\phantom{^2}&=\eps_ky_k^2,
    \end{aligned}\right.\]
    in the variables $x,z,y_0,y_1,\dots,y_k$. Let $U\coloneqq V\cap\{y_0y_1\dots y_k\neq0\}\openset V$. Let $g\coloneqq(k-2)2^{k-2}+1$. Then,
    \[(q+1)-2g\sqrt q>(k+2)2^{k-1}\implies U(\F_q)\neq\emptyset.\]
    In particular, if $k=2$, then $g=1$ and $q>14\implies U(\F_q)\neq\emptyset$.
\end{lemma}
\begin{proof}
    Consider the projective curve $C\subset\P^{k+1}_{\F_q}$ (with coordinates $X,Z,Y_0,Y_1,Y_2,\dots,Y_{k-1}$) cut out by the $k$ equations
    \begin{equation}\label{eqn:C-lang-weil}
        C\colon\left\{\begin{aligned}
        s_0X^2 + t_0 Z^2 &= \eps_0 Y_0^2 \\
        (s_1X+t_1Z)(s_kX+t_kZ) &= \eps_1\eps_k Y_1^2 \\
        &\;\;\vdots \\
        (s_{k-1}X+t_{k-1}Z)(s_kX+t_kZ) &= \eps_{k-1}\eps_k Y_{k-1}^2
    \end{aligned}\right..
    \end{equation}
    Note that there is a map $\phi\colon V\to C$ given by
    \[\phi\p{x,z,y_0,\dots,y_k}=\sq{x,z,y_0,y_1y_k,\dots,y_{k-1}y_k}.\]
    One readily checks that any $\F_q$-point of $C\cap\{Y_0Y_1\dots Y_{k-1}\neq0\}$ lifts to an $\F_q$-point of $U\subset V$. Thus, it suffices to find a criterion guaranteeing that $C$ has an $\F_q$-point with $Y_0Y_1\dots Y_{k-1}\neq0$. Let $\wt C$ be the normalization of $C$, so $\wt C$ is the unique smooth, proper $\F_q$-curve with function field
    \[\F_q(\wt C)=\F_q(u)\p{\sqrt{\frac{s_0u^2+t_0}{\eps_0}},\sqrt{\frac{(s_1u+t_1)(s_ku+t_k)}{\eps_1\eps_k}},\dots,\sqrt{\frac{(s_{k-1}u+t_{k-1})(s_ku+t_k)}{\eps_{k-1}\eps_k}}},\]
    where $u\coloneqq X/Z$. In particular, the natural inclusion $\F_q(u)\into\F_q(\wt C)$ corresponds to a finite cover $f\colon\wt C\to\P^1$, of degree $2^k$, branched over the locus
    \[\{s_0X^2+t_0Z^2\}\sqcup\{s_1X+t_1Z\}\sqcup\dots\sqcup\{s_kX+t_kZ\}\subset\P^1_{\F_q}.\]
    Since $f$ corresponds to a multi-quadratic field extension, every ramification point in $\wt C$ has ramification degree $2$, and the total number of ramification points above each branch point is $\deg(f)/2=2^{k-1}$. Thus, Riemann--Hurwitz yields
    \[2g(\wt C)-2=2^k(-2)+2^{k-1}(2+k)(2-1),\tandso g(\wt C)=g\coloneqq(k-2)2^{k-2}+1.\]
    Finally, the Hasse--Weil bound ensures that $\#\wt C(\F_q)\ge(q+1)-2g\sqrt q$;
    since a point in $\wt C(\F_q)$ maps into $\phi(U(\F_q))\subset C(\F_q)$ if and only if it is a ramification point of $f$ and we computed above that the ramification locus of $f$ has degree $(k+2)2^{k-1}$, the statement of the lemma follows.
\end{proof}

\begin{rem}
    The bound $q>14$ when $k=2$ in Lemma \ref{lem:Lang-Weil-curve-trick} is optimal. For each odd prime power $q\le14$, there exists such a variety $V$ all of whose $\F_q$-points satisfy $y_0y_1y_2=0$. Examples of such $V$ are recorded in Table \ref{tab:lang-weil-counterexamples} below.
    \begin{table}[h]
        \centering
        \begin{tabular}{c|ccc}
            $q$ & & $V/\F_q$ \\ \hline
            $3$ & $x^2+\phantom9z^2=2y_0^2$ & $x+\phantom2z=y_1^2$ & $x-\phantom4z=\phantom2y_2^2$ \\ 
            $5$ & $x^2+\phantom9z^2=\phantom3y_0^2$ & $x+\phantom2z=y_1^2$ & $x-\phantom4z=2y_2^2$ \\
            $7$ & $x^2+\phantom9z^2=3y_0^2$ & $x+2z=y_1^2$ & $x+4z=\phantom2y_2^2$ \\
            $9$ & $x^2+\phantom9z^2=\phantom3y_0^2$ & $x+\phantom2z=y_1^2$ & $x-\phantom4z=\eps y_2^2$ \\
            $11$ & $x^2-9z^2=\phantom3y_0^2$ & $x+\phantom2z=y_1^2$ & $x+5z=2y_2^2$ \\
            $13$ & $x^2+\phantom9z^2=2y_0^2$ & $x+\phantom2z=y_1^2$ & $x-\phantom4z=2y_2^2$
        \end{tabular}
        \caption{Counterexamples to Lemma \ref{lem:Lang-Weil-curve-trick} when $k=2$ and $q\le14$. Above, $\eps\in\units\F_9$ is any fixed choice of nonsquare.}
        \label{tab:lang-weil-counterexamples}
    \end{table}

\end{rem}

\begin{rec}\label{rec:all-the-notation}
    At this point, $(\mu,\lambda)$ is a 4-prepared admissible pair. Among other things, this means that, in our current notation, all of the following hold:
    \begin{itemize}
        \item $\lambda=\mu m^2$, with $\mu<0$ (so also $\lambda<0$), $\mu\equiv1\bmod8$, and $m\equiv2\bmod8$.
        \item $p_1,p_2$ each divide $(m^2-1)$, but do \important{not} divide $\mu$ or $\lambda$. Furthermore, $2,\mu\equiv\square\bmod p_i$ for each $i$.
        \item $\l_1$ divides $m$ (and so divides $\lambda=\mu m^2$), but does \important{not} divide $\mu$. Furthermore, $-\mu,2\equiv\square\bmod\l_1$.
        \item $\varpi$ divides $\mu$, and so \important{also} divides $\lambda$.
        \item $q_1,q_2,q_3,q_4$ each divide $\mu$, so \important{also} divide $\lambda$. Furthermore, 
        \begin{align}
            1 - m^2 &\equiv\square\bmod q_2,q_3,q_4 \\
            1 - m^2 &\not\equiv\square\bmod q_1 \\
            2 &\equiv\square\bmod q_2,q_3,q_4
            .
        \end{align}
        \item A finite prime lies in the set $S(E)$ if and only if it divides $\mu m(m^2-1)$.
        \item The quantities $\mu$, $m$, $m^2-1$, and $m^2+1$ are pairwise coprime.
    \end{itemize}
    See Definitions \ref{defn:admissible} and \ref{defn:n-prepared}.
\end{rec}
\begin{lemma}\label{lem:local conditions for linear transformation}
    Set $M\coloneqq8[\mu m(m^2+1)(m^2-1)]^2$. 
    There exists $\xi,\upsilon,a,b\in\zmod M$ such that
    \begin{itemize}
        \item $\mu\mid(a+b)$ but $\mu$ is coprime to the four quantities $a-b$, $a\pm m^2b$, and $a^2+m^2b^2$; 
        \item $(1+m^2)\mid(a-b)$, but $1+m^2$ is coprime to $a+b$;
        \item letting $\pi_i(X_0,Y_0)\in\Z[X_0,Y_0]$ as in Lemma \ref{lem:modified linear+quadratic forms}, and setting $\pi_i\coloneqq\pi_i(\xi,\upsilon)\in\zmod M$, we additionally have
        \begin{alphabetize}
            \item $\pi_i$ is invertible mod $M$, for each $i$, and
            \[\pi_1\pi_2\pi_3\pi_4\pi_5\equiv\square\bmod8p\quad\t{for all odd }p\mid\mu m(m^2-1).\]
            \item For primes $p\mid\mu$, if $p\not\in\{\varpi,q_2,q_3,q_4\}$, then
            \[\pi_1\equiv\square\bmod p\tand\pi_2(1-m^2)\equiv\square\bmod p.\]
            For the prime $p=\varpi$,
            \[\pi_1\not\equiv\square\bmod\varpi\tand\pi_2(1-m^2)\not\equiv\square\bmod\varpi.\]
            For primes $p\in\{q_2,q_3,q_4\}$,
            \[\pi_1,\pi_2,\pi_3\equiv\square\bmod p.\]
            \item For the prime $p=p_1$,
            \[\pi_1\equiv\square\bmod p_1\tbut\pi_2\not\equiv\square\bmod p_1.\]
            For the prime $p=p_2$,
            \[\pi_1,\pi_2,\pi_3\equiv\square\bmod p_2.\]
            \item For the prime $p=\l_1$,
            \[\pi_1,\pi_2\equiv\square\bmod\l_1,\tbut\pi_3\not\equiv\square\bmod\l_1.\]
            \item Finally,
            \[\pi_1,\pi_2\equiv1\bmod4,\tbut\pi_3\equiv3\bmod4.\]
        \end{alphabetize}
    \end{itemize}
\end{lemma}
\begin{proof}
\noindent\textbf{Case 1 ($p$ divides $1+m^2$):}\footnote{We recommend the reader take intermittent chocolate breaks while reading the below casework.} In this case, we want to check that we can choose $\xi,\upsilon,a,b \in \Z/p^2\Z$ such that $\pi_i(\xi,\upsilon)\in \F_p^\times$ for each $i$. Since $m$ is even, $p$ must be odd and therefore $p>2$. Assume that $a-b \equiv cp \bmod p^2$, where $c\neq 0\bmod p$. Then, (mod $p$)
\begin{align}
    &\pi_1 \equiv \mu^{-1}(a+b), \quad \pi_2 \equiv \mu(\xi-\upsilon) + c, \quad  \pi_3 \equiv \mu(\xi - \upsilon) + c
    \\ & \pi_4 \equiv a+b, \quad \pi_5 \equiv 2\mu(\xi a - \upsilon b) + c(a+b). 
\end{align}
Let $\gcd(a,p)=1$ and $b=a-cp$. Then $a+b \equiv 2a-cp$ will satisfy that $\gcd(a+b,p)=1$. This is sufficient for $\pi_1$ and $\pi_4$. Next, take $\xi\equiv\upsilon\equiv 0 \bmod p$; then since $\gcd(c,p)=1$, we have that $\pi_2,\pi_3,\pi_5\in \F_p^\times$ as well. 

\noindent\textbf{Case 2 ($p$ divides $\mu$):} Assume that $a+b=cp$ (as $\mu\mid a+b$). Let us compute the reductions of $\pi_1,\dots,\pi_5$ mod $p$:
\begin{align}
    &\pi_1 \equiv  (1+m^2)(\xi+\upsilon) + c, \quad \pi_2 \equiv (1+m^2)^{-1}(a-b), \quad \pi_3 \equiv (1+m^2)(a+m^2b), \\
    &\pi_4 \equiv a-m^2 b, \quad \pi_5\equiv (1+m^2) (a^2+m^2 b^2).
\end{align}
Since $a+b\equiv 0 \bmod p$, we have that this further reduces to
\begin{align}
    &\pi_1 \equiv  (1+m^2)(\xi+\upsilon) + c, \quad \pi_2 \equiv 2(1+m^2)^{-1}a, \quad \pi_3 \equiv (1+m^2)(1-m^2)a, \\
    &\pi_4 \equiv (1+m^2)a, \quad \pi_5\equiv (1+m^2)^2 a^2.
\end{align}
Note that the variables $\xi$ and $\upsilon$, as well as the choice of $c$, only appear in the equation for $\pi_1\bmod p$. Additionally, we see that 
\begin{equation}
    \pi_1\pi_2\pi_3\pi_4\pi_5 \equiv \pi_1 \cdot 2(1+m^2)^3(1-m^2) a^5.
\end{equation}
\begin{enumerate}[label=(\Alph*)]
    \item $p\not\in \{\varpi,q_2,q_3,q_4\}$: Set $\xi,\upsilon,c$ such that $\pi_1 \equiv \square \bmod p$ (this can be done since $\pi_1$ is linear in $\xi,\upsilon,c$). Next, choose $a$ such that $2(1+m^2)(1-m^2)a\in (\F_p^\times)^2$. This forces $\pi_1\pi_2\pi_3\pi_4\pi_5\equiv\square\bmod p$ and $\pi_2(1-m^2)\equiv\square\bmod p$, so we have a solution to the desired local condition at $p$. 

    \item $p = \varpi:$ Again, since $\pi_1$ is linear in $\xi,\upsilon,c$ and the only form that depend on these variables, we can choose $\xi,\upsilon,c$ such that $\pi_1 \not\equiv \square \bmod p$. Next, choose $a$ such that $2(1+m^2)(1-m^2) a\not\equiv \square \bmod p$. Then we have a solution to the desired local condition at $p$. 


    \item $p\in\{q_2,q_3,q_4\}$: Choose $\xi,\upsilon,c$ such that $\pi_1\equiv\square\bmod p$, and choose $a$ so that $2(1+m^2)a\equiv\square\bmod p$ (so $\pi_2\equiv\square\bmod p$). Since $1-m^2,2\equiv\square\bmod p$, it follows that also $\pi_3$ and $\pi_1\pi_2\pi_3\pi_4\pi_5$ are squares mod $p$.

\end{enumerate}

\noindent\textbf{Case 3 ($p$ divides $m^2-1$):} Let us first consider the reductions of the forms $\pi_1,\dots,\pi_5\bmod p$. 
\begin{align}
    &\pi_1 \equiv 2(\xi+\upsilon) + \mu^{-1}(a+b), \quad \pi_2 \equiv \mu(\xi-\upsilon) + \frac{a-b}{2}, \\
    &\pi_3\equiv \mu(\xi+\upsilon) + \frac{a+b}{2}, \quad \pi_4 \equiv 2\mu(\xi - \upsilon) + (a-b) \\
    &\pi_5 \equiv 2\mu^2 ( \xi^2+\upsilon^2) + 2\mu(\xi a+\upsilon b) + \frac{a^2+b^2}{2} = \frac{((2\mu \xi + a)^2 + (2\mu \upsilon+b)^2}{2}.
\end{align}
Observe that $\mu \pi_1 \equiv 2 \pi_3$ and that $2\pi_2 \equiv \pi_4.$ For now, set $r = 2\mu \xi + a$ and $s = 2\mu \upsilon + b$. Then we have that 
\begin{align}
    &\pi_1 \equiv \mu^{-1}(r+s), \quad \pi_2 \equiv  2^{-1}(r-s), \quad \pi_3 \equiv 2^{-1}(r+s), \quad \pi_4 \equiv r-s, \\
    &\pi_5 \equiv 2^{-1}(r^2+s^2).
\end{align}
\begin{enumerate}[label=(\Alph*)]
    \item $p=p_1$: We claim that we can choose $r$ and $s$ such that $\pi_1 \equiv \square$, $\pi_2 \not\equiv \square$. First, observe that $\pi_1\pi_2\pi_3\pi_4\pi_5 \equiv 2^{-3} \mu^{-1} (r^2+s^2)(r^2-s^2)^2$, so we require that $2(r^2+s^2)\equiv \square$ (recall that $\mu \equiv \square$). So, let $\eps_1,\eps_2,\eps_5 \in \F_p^\times $ be chosen such that $\eps_1$ is a square, $\eps_2$ is not a square, and $2\eps_5$ is a square. It then suffices to show that there is a nontrivial solution to
    \begin{equation}
        \begin{cases}
            r+s = \eps_1 y_1^2 \\
            r-s = \eps_2 y_2^2 \\
            r^2+s^2 = \eps_5 y_5^2.
        \end{cases}
    \end{equation}
    By Lemma \ref{lem:Lang-Weil-curve-trick}, so long as $p_1>14$, we can find such $r$ and $s$. 

    \item $p=p_2$: We now want to find $r$ and $s$ such that $\pi_1,\pi_2,\pi_3\equiv\square$ (and $\pi_1\pi_2\pi_3\pi_4\pi_5\equiv\square$ as well). Since $2,\mu\equiv\square\bmod p$, by assumption, this amounts to finding a nontrivial solution to the equations
    \begin{equation}
        \begin{cases}
            r+s = y_1^2 \\
            r-s =  y_2^2 \\
            r^2+s^2 = y_5^2.
        \end{cases}
    \end{equation}
    By Lemma \ref{lem:Lang-Weil-curve-trick}, this is possible as long as $p_2>14$.
\end{enumerate}

\noindent\textbf{Case 4 ($p$ divides $m$):} We start off by reducing the equations mod $m$. 
\begin{align}
    &\pi_1 \equiv (\xi+\upsilon) + \mu^{-1}(a+b), \quad  \pi_2 \equiv \mu(\xi-\upsilon) + (a-b), \\
    &\pi_3 \equiv \mu \xi + a, \quad \pi_4 \equiv \mu \xi + a, \\
    & \pi_5 \equiv \mu^2 \xi^2 + 2\mu \xi a + a^2 = (\mu \xi+a)^2.
\end{align}
Observe that the condition $\pi_1\pi_2\pi_3\pi_4\pi_5\equiv \square$ can be reduced to the condition that $\pi_1\pi_2\equiv \square$. 

\begin{enumerate}[label=(\Alph*)]
    \item $p=\ell_1$: We see that for this condition, it is sufficient to find $\xi,\upsilon,a,b$ such that 
    \begin{equation}
        \begin{cases}
            (\xi+\upsilon) + \mu^{-1}(a+b) \equiv c_1 \\
            \mu(\xi-\upsilon) + (a-b) \equiv c_2 \\
            \mu\xi + a \equiv c_3
            ,
        \end{cases}
    \end{equation}
    where $c_1,c_2$ are fixed squares mod $p$, while $c_3$ is a fixed nonsquare mod $p$. Note that we can rewrite the above equations as: 
    \begin{equation}
        \begin{cases}
            c_3 + (\mu \upsilon + b) \equiv \mu c_1, \\
            c_3 - (\mu \upsilon +b ) \equiv c_2, \\
            \mu \xi + a \equiv c_3.
        \end{cases}
    \end{equation} 
    Writing $c_1=y_1^2$, $c_2 = y_2^2$, and $c_3= \eps y_3^2$ for $\eps\in \F_p^\times$ not a square, we see that adding the first two equations, we get $2\eps y_3^2 = \mu y_1^2 + y_2^2.$ Further, since $\mu \xi + a$ is linear, we can find a solution to $\mu \xi + a = \eps y_3^2$ for any choice of $y_3$. Similarly, given values of $y_1,y_2,y_3$, and since $2\in \F_p^\times$, we can find $\upsilon$ and $b$ such that $\mu \upsilon + b = 2^{-1}(\mu y_1^2 - y_2^2).$ 
    Thus, it suffices to find a solution to $2\eps y_3^2 = \mu y_1^2 + y_2^2$ with $y_1,y_2,y_3\in \F_p^\times$. Since every quadratic form in $\geq 3$ variables over $\F_p$ has a nontrivial solution, we can find such a solution to $y_1,y_2,y_3$, and hence to $\xi,\upsilon,a,b$. 
\end{enumerate}

\noindent\textbf{Case 5 ($p=2$):}
In this case, we want to find $\xi,\nu,a,b \in \Z/8\Z$ such that $\pi_1\equiv 1 \bmod 4$, $\pi_2\equiv1\bmod4$, $\pi_3\equiv3\bmod4$, and $\pi_1\pi_2\pi_3\pi_4\pi_5\equiv 1\bmod 8$. Recall that $m\equiv 2 \bmod 4$. Set $x'=\mu\xi+5a$ and $y'=\mu\upsilon+5b$. Reducing mod 8 gives 
\begin{align}
    &\pi_1\equiv5(\xi+\upsilon)+\inv\mu(a+b) \equiv5\mu(x'+y')\bmod8,\quad \pi_2\equiv\mu(\xi-\upsilon)+5(a-b)\equiv x'-y'\bmod8,\\
    &\pi_3\equiv\mu(\xi+4\upsilon)+5(a+4b)\equiv x'+4y'\bmod8,\quad \pi_4\equiv5\mu(\xi-4\upsilon)+(a-4b) \equiv5(x'-4y')\bmod8,\\
    &\pi_5\equiv5\sq{(\mu\xi+5a)^2+4(\mu\upsilon+5b)^2} \equiv5(x'^2+4y'^2)\bmod8
    .
\end{align}
In particular, when $x'$ is odd and $y'$ is even, we see that $\pi_1\pi_2\pi_3\pi_4\pi_5$ is congruent to $5\mu(1+y'^2)\bmod8$. Choose $\xi,\upsilon,a,b$ so that $x'\equiv3\bmod8$ and $y'\equiv2\bmod8$. Then, recalling that $\mu\equiv1\bmod8$, $\pi_1\equiv\pi_2\equiv\pi_1\pi_2\pi_3\pi_4\pi_5\equiv1\bmod8$ and $\pi_3\equiv3\bmod8$, so all the conditions are satisfied.

\end{proof}

Fix $M,a,b,\xi,\upsilon$ as in Lemma \ref{lem:local conditions for linear transformation}, and let $\pi_i(X_0,Y_0)\in\Z[X_0,Y_0]$, for $i=1,\dots,5$, be as in Lemma \ref{lem:modified linear+quadratic forms}.
Define a new set of polynomials:
\begin{equation}
    \tilde{\pi}_i(X_0,Y_0) := \pi_i(MX_0+\xi,MY_0+\upsilon).
\end{equation}
Lemma \ref{lem:local conditions for linear transformation} shows that the singular series appearing in Theorem \ref{thm: prime count} applied to $\tilde{\pi}_1(X_0,Y_0),\dots,\tilde{\pi}_5(X_0,Y_0)$ is nonzero. Thus, Theorem \ref{thm: prime count} supplies us with infinitely many $x_0,y_0\in\Z$ such that $\{\pi_i=\pi_i(Mx_0+\xi,My_0+\upsilon)\}_{i=1}^5\subset\Z$ is a set of primes satisfying the conditions of Proposition \ref{prop:rank-2-crit}. Observe that we need to restrict to when $\pi_1<0$ and $\pi_2>0$, rather than counting over the region given in the exact statement of Theorem \ref{thm: prime count}. However, it would suffice to do the prime count in the region: $\{X_0<0, m^2|Y_0|>|X_0|>|Y_0|\}$, which is convex. It is straightforward to adapt the proof of Theorem \ref{thm: prime count} to this regime (using Fourier analysis to approximate the region, as was done in \cite{Green_Sawhney} for the box). Consequently, we produce infinitely many $u\in\Q$ such that $\dim\Sel_2(E^{G(u)})\le4$. By Silverman's specialization theorem \cite[Theorem C]{Silverman:heights-specialization} and Proposition \ref{prop:RS-family}, $\rank E^{G(u)}(\Q)\ge2$ for all but finitely many $u$. Thus, there are infinitely many $u\in\Q$ satisfying both
\[\dim\Sel_2(E^{G(u)})\le4\tand\rank E^{G(u)}(\Q)\ge2.\]
Since $\dim\Sel_2E\ge\rank E(\Q)+\dim E[2](\Q)$ for any elliptic curve $E/\Q$, we conclude from this that $\rank E^{G(u)}(\Q)=2$ for any such $u$, proving Theorem \ref{thm:rank2}.
\qed

\appendix
\appendixpage
\addappheadtotoc

\section{A simple proof of Theorem \ref{thma:BMO}}\label{app:BMO}

In this appendix, for readers newer to the fibration method, we include the details of a proof of Theorem \ref{thma:BMO}, a special case of the more general \ref{thm:general-BMO}, which we restate below.
\begin{thm}[= Theorem \ref{thma:BMO}]
    Let $P(x)\in\Q[x]$ be a polynomial factoring as
    \[P(x)=L_1(x)\cdots L_k(x)\cdot Q(x)\twhere\deg L_i=1\tand\deg Q=2.\]
    Assume $K=\Q[x]/Q(x)$ is an imaginary quadratic field, and let $X$ be a smooth compactification of the affine variety
    \begin{equation}\label{eqn:affine-conic-bundle}
        y^2-az^2=P(x)\twith a\in\units\Q.
    \end{equation}
    Then, $X(\Q)$ is dense in $X(\A_Q)^{\Br}$. 
\end{thm}

\begin{proof}

Let $\pi\colon X\rightarrow \P_\Q^1$ denote projection onto the $x$ coordinate. Let $e_1,\dots,e_r\in\P^1_\Q$ denote the closed points (with residue field $\Q$) cut out by $L_1(x),\dots,L_k(x)$, and let $e_0\in\P^1_\Q$ be the closed point (with reside field $K$) cut out by $Q(x)$. For $i=0,\dots,r$, let $X_r\coloneqq X_{e_r}=\inv\pi(e_r)$. Let $U\coloneqq\P^1_\Q\sm\{e_0,\dots,e_r,\infty\}$, and set
$V\coloneqq\inv\pi (U)\openset X$.
Define the following polynomials:
\[\pull Q(t)=N_{K/\Q}(t-e_0)\tand\pull L_i(t)=t-e_i\tfor i=1,\dots,r.\]
Write $S_X = S_{\pull L_1,\dots,\pull L_r,\pull Q}$, defined as above Lemma \ref{lem:BMO hypothesis H-1 sub}.

Let $\Omega$ denote the set of all places of $\Q$ and take $(P_v)_{v\in \Omega}\in X(\A_\Q)^{\Br}$. Let $S\supset S_X\cup\{\infty\}$ be a given finite set of primes, including all primes which ramify in $\Q(\sqrt a)/\Q$.
To establish Theorem \ref{thma:BMO}, it suffices to find some $\lambda\in U(\Q)\subset\A^1(\Q)=\Q$ which is sufficiently close to each $\lambda_v\coloneqq\pi(P_v)$, $v\in S$, and such that $X_\lambda$ is locally solvable. Since $X_\lambda$ is a smooth conic, it satisfies weak approximation, so some rational point on it will approximate $(P_v)_{v\in\Omega}$.
For this, we may assume that, for each $v\in S$, $P_v\in V(\Q_v)$; if not, we can replace $P_v$ with a sufficiently close $P_v'$ which does lie in this open set.

Consider the classes
\begin{equation}
    A_0\coloneqq\opname{Cores}_{K/\Q}\p{K(\sqrt a)/K, t - e_0}\tand A_i\coloneqq\p{\Q(\sqrt a)/\Q, t-e_i},\,\t{for }i=1,\dots,r
\end{equation}
in $\Br\Q(t)$. 

Let $B \subset \Br(V)$ denote the finite group generated by $\pull\pi A_i$, for $i=0,\dots,r$.
Then, $(P_v)_{v\in\Omega}\in V(\A_\Q)^{B\cap\Br(X)}$, so, by Harari's formal lemma \cite[Theorem 13.4.3]{CT-Skorobogatov:BrauerGrothendieck}, we can find a set of primes $S_1\supset S$ and choose $P_v\in V(\Q_v)$ for $v\in S_1\setminus S$ such that 
\begin{equation}\label{eqn: A_ij sum of invariants}
    \sum_{v\in S_1} \Inv_v(A_i(\pi(P_v))) = 0,\quad \forall i=0,\dots,r.
\end{equation}
Set $\lambda_v=\pi(P_v)\in \Q_v$ for each $v\in S_1$. Applying Lemma \ref{lem:BMO hypothesis H-1 sub} with the set of primes $S_1$, we can find $\lambda \in \Q$, sufficiently close to $\lambda_v\in\Q_v$ for all $v\in S_1$, and primes $\{p_1,\dots,p_r,q\}\not\in S_1$ satisfying properties \ref{item:H1 sub-close finite} -- \ref{item:H1 sub-uniformizers} given in that lemma. 
We now claim that $X_\lambda(\A_\Q)\neq \emptyset$.

For $v\in S_1$,  the inverse function theorem gives us that $X_{\lambda}(\Q_v)\neq \emptyset$ (as we know that $\pi^{-1}(\lambda_v)\neq \emptyset$ for each $\lambda_v\in U(\Q_v)$ for $v\in S_1$). For $v\not\in S_1\cup \{q,p_1,\dots,p_r\}$,
properties \ref{item:H1 sub-S-unit L}, \ref{item:H1 sub-S-unit Q} guarantee
that $\val_v(\lambda-e_i)\le0$ and $\val_v(N_{K/\Q}(\lambda-e_0))\le0$. Hence, $\lambda\bmod v \in \P^1_{\F_v}$ lies outside the reduction of $e_0,e_1,\dots,e_r$ 
 so $X_\lambda\bmod v$ is a smooth conic over $\F_v$; it is thus isomorphic to $\P^1_{\F_v}$ and Hensel's lemma ensures that all its points lift to $\Q_v$-points on $X_\lambda$.

Next, let $v=p_i$, for $i\in\{1,\dots,r\}$. Property \ref{item:H1 sub-hilbert Q} combined with (\ref{eqn: A_ij sum of invariants}) together ensure that $0=\Inv_{p_i}(\Q(\sqrt a)/\Q,\lambda-e_i)$. Since $\val_{p_i}(\lambda-e_i)=1$ (i.e. $\lambda-e_i$ is a uniformizer at $p_i$), this means that $a$ is a square mod $p_i$. It is then clear from (\ref{eqn:affine-conic-bundle}) that $X_\lambda$ has the smooth $\F_{p_i}$-point $(x,y,z)=(\lambda,\sqrt a,1)\in\F_{p_i}^3$ and so $X_\lambda(\Q_{p_i})\neq\emptyset$.

Finally, consider $v=q$. One argues as in the previous case, with property \ref{item:H1 sub-hilbert K} in place of property \ref{item:H1 sub-hilbert Q}, to see that $a$ is again a square mod $q$, from which it easily follows that $X_\lambda$ must have a $\Q_q$-point.

Thus, $X_{\lambda}(\A_\Q)\neq \emptyset$, from which the claim follows.
\end{proof}

\printbibliography

@article{Koymans-Pagano:H10,
      AUTHOR = {Koymans, Peter and Pagano, Carlo},
     TITLE = {Hilbert's tenth problem via additive combinatorics},
   JOURNAL = {J. Amer. Math. Soc.},
  FJOURNAL = {Journal of the American Mathematical Society},
    VOLUME = {40},
      YEAR = {2027},
    NUMBER = {1},
     PAGES = {195--234},
      ISSN = {0894-0347,1088-6834},
   MRCLASS = {11G05 (03D35 11D99 11U05)},
  MRNUMBER = {5131234},
       DOI = {10.1090/jams/1081},
       URL = {https://doi.org/10.1090/jams/1081},
}

@misc{kai2025linearpatternsprimeelements,
      title={Linear patterns of prime elements in number fields}, 
      author={Wataru Kai},
      year={2025},
      eprint={2306.16983},
      archivePrefix={arXiv},
      primaryClass={math.NT},
}

@article{Klagsbrun_Mazur_Rubin_2014, title={A Markov model for Selmer ranks in families of twists}, volume={150}, DOI={10.1112/S0010437X13007896}, number={7}, journal={Compositio Mathematica}, author={Klagsbrun, Zev and Mazur, Barry and Rubin, Karl}, year={2014}, pages={1077–1106}}

@article {mitsui,
    AUTHOR = {Mitsui, Takayoshi},
     TITLE = {Generalized prime number theorem},
   JOURNAL = {Jpn. J. Math.},
  FJOURNAL = {Japanese Journal of Mathematics},
    VOLUME = {26},
      YEAR = {1956},
     PAGES = {1--42},
      ISSN = {0075-3432},
   MRCLASS = {10.1X},
  MRNUMBER = {92814},
MRREVIEWER = {S.\ Chowla},
       DOI = {10.4099/jjm1924.26.0\_1},
       URL = {https://doi.org/10.4099/jjm1924.26.0_1},
}

@book {silverman,
    AUTHOR = {Silverman, Joseph H.},
     TITLE = {The arithmetic of elliptic curves},
    SERIES = {Graduate Texts in Mathematics},
    VOLUME = {106},
   EDITION = {Second},
 PUBLISHER = {Springer, Dordrecht},
      YEAR = {2009},
     PAGES = {xx+513},
      ISBN = {978-0-387-09493-9},
   MRCLASS = {11-02 (11G05 11G20 14H52 14K15)},
  MRNUMBER = {2514094},
MRREVIEWER = {Vasil\cprime  \={I}. Andr\={\i}\u{\i}chuk},
       DOI = {10.1007/978-0-387-09494-6},
       URL = {https://doi.org/10.1007/978-0-387-09494-6},
}

@article {poonen-rains,
    AUTHOR = {Poonen, Bjorn and Rains, Eric},
     TITLE = {Random maximal isotropic subspaces and {S}elmer groups},
   JOURNAL = {J. Amer. Math. Soc.},
  FJOURNAL = {Journal of the American Mathematical Society},
    VOLUME = {25},
      YEAR = {2012},
    NUMBER = {1},
     PAGES = {245--269},
      ISSN = {0894-0347,1088-6834},
   MRCLASS = {11G10 (11G05 11G30 14G25 14K15)},
  MRNUMBER = {2833483},
MRREVIEWER = {Philipp\ Habegger},
       DOI = {10.1090/S0894-0347-2011-00710-8},
       URL = {https://doi.org/10.1090/S0894-0347-2011-00710-8},
}

@article {period-index,
    AUTHOR = {O'Neil, Catherine},
     TITLE = {The period-index obstruction for elliptic curves},
   JOURNAL = {J. Number Theory},
  FJOURNAL = {Journal of Number Theory},
    VOLUME = {95},
      YEAR = {2002},
    NUMBER = {2},
     PAGES = {329--339},
      ISSN = {0022-314X,1096-1658},
   MRCLASS = {11G05},
  MRNUMBER = {1924106},
MRREVIEWER = {Imin\ Chen},
       DOI = {10.1016/S0022-314X(01)92770-2},
       URL = {https://doi.org/10.1016/S0022-314X(01)92770-2},
}

@article {period-index-erratum,
    AUTHOR = {O'Neil, Catherine},
     TITLE = {Erratum to: ``{T}he period-index obstruction for elliptic
              curves'' [{J}. {N}umber {T}heory {\bf 95} (2002), no. 2,
              329--339; MR1924106]},
   JOURNAL = {J. Number Theory},
  FJOURNAL = {Journal of Number Theory},
    VOLUME = {109},
      YEAR = {2004},
    NUMBER = {2},
     PAGES = {390},
      ISSN = {0022-314X,1096-1658},
   MRCLASS = {11G05},
  MRNUMBER = {2106488},
       DOI = {10.1016/j.jnt.2004.07.006},
       URL = {https://doi.org/10.1016/j.jnt.2004.07.006},
}

@misc{koymans2025ellipticcurvesranknumber,
      title={Elliptic curves of rank one over number fields}, 
      author={Peter Koymans and Carlo Pagano},
      year={2025},
      eprint={2505.16910},
      archivePrefix={arXiv},
      primaryClass={math.NT},
}

@book {poonen-rat-pts,
    AUTHOR = {Poonen, Bjorn},
     TITLE = {Rational points on varieties},
    SERIES = {Graduate Studies in Mathematics},
    VOLUME = {186},
 PUBLISHER = {American Mathematical Society, Providence, RI},
      YEAR = {2017},
     PAGES = {xv+337},
      ISBN = {978-1-4704-3773-2},
   MRCLASS = {14G05 (11G35)},
  MRNUMBER = {3729254},
MRREVIEWER = {Daniel Loughran},
       DOI = {10.1090/gsm/186},
       URL = {https://doi.org/10.1090/gsm/186},
}

@article {Klagsbrun-Mazur-Rubin:disparity,
    AUTHOR = {Klagsbrun, Zev and Mazur, Barry and Rubin, Karl},
     TITLE = {Disparity in {S}elmer ranks of quadratic twists of elliptic
              curves},
   JOURNAL = {Ann. of Math. (2)},
  FJOURNAL = {Annals of Mathematics. Second Series},
    VOLUME = {178},
      YEAR = {2013},
    NUMBER = {1},
     PAGES = {287--320},
      ISSN = {0003-486X,1939-8980},
   MRCLASS = {11G05 (14G25)},
  MRNUMBER = {3043582},
MRREVIEWER = {Jan\ Nekov\'a\v r},
       DOI = {10.4007/annals.2013.178.1.5},
       URL = {https://doi.org/10.4007/annals.2013.178.1.5},
}

@article {CT-Sansuc:schinzel-hasse,
    AUTHOR = {Colliot-Th\'{e}l\`ene, Jean-Louis and Sansuc, Jean-Jacques},
     TITLE = {Sur le principe de {H}asse et l'approximation faible, et sur
              une hypoth\`ese de {S}chinzel},
   JOURNAL = {Acta Arith.},
  FJOURNAL = {Polska Akademia Nauk. Instytut Matematyczny. Acta Arithmetica},
    VOLUME = {41},
      YEAR = {1982},
    NUMBER = {1},
     PAGES = {33--53},
      ISSN = {0065-1036},
   MRCLASS = {10C02 (10B40 12B20 14G25)},
  MRNUMBER = {667708},
MRREVIEWER = {D. J. Lewis},
       DOI = {10.4064/aa-41-1-33-53},
       URL = {https://doi.org/10.4064/aa-41-1-33-53},
}

@article {CT-SD:schinzel-hasse-severi-brauer,
    AUTHOR = {Colliot-Th\'{e}l\`ene, Jean-Louis and Swinnerton-Dyer, Peter},
     TITLE = {Hasse principle and weak approximation for pencils of
              {S}everi-{B}rauer and similar varieties},
   JOURNAL = {J. Reine Angew. Math.},
  FJOURNAL = {Journal f\"{u}r die Reine und Angewandte Mathematik. [Crelle's
              Journal]},
    VOLUME = {453},
      YEAR = {1994},
     PAGES = {49--112},
      ISSN = {0075-4102},
   MRCLASS = {11G35 (11G45 14C25 14G05 14G25)},
  MRNUMBER = {1285781},
MRREVIEWER = {Wayne Raskind},
       DOI = {10.1515/crll.1994.453.49},
       URL = {https://doi.org/10.1515/crll.1994.453.49},
}

@article {browning-matthiesen-skorobogatov,
    AUTHOR = {Browning, T. D. and Matthiesen, L. and Skorobogatov, A. N.},
     TITLE = {Rational points on pencils of conics and quadrics with many
              degenerate fibers},
   JOURNAL = {Ann. of Math. (2)},
  FJOURNAL = {Annals of Mathematics. Second Series},
    VOLUME = {180},
      YEAR = {2014},
    NUMBER = {1},
     PAGES = {381--402},
      ISSN = {0003-486X},
   MRCLASS = {14G05 (14C21 14G25)},
  MRNUMBER = {3194818},
MRREVIEWER = {Cecilia Salgado},
       DOI = {10.4007/annals.2014.180.1.8},
       URL = {https://doi.org/10.4007/annals.2014.180.1.8},
}

@book {CT-Skorobogatov:BrauerGrothendieck,
    AUTHOR = {Colliot-Th\'{e}l\`ene, Jean-Louis and Skorobogatov, Alexei N.},
     TITLE = {The {B}rauer-{G}rothendieck group},
    SERIES = {Ergebnisse der Mathematik und ihrer Grenzgebiete. 3. Folge. A
              Series of Modern Surveys in Mathematics [Results in
              Mathematics and Related Areas. 3rd Series. A Series of Modern
              Surveys in Mathematics]},
    VOLUME = {71},
 PUBLISHER = {Springer, Cham},
      YEAR = {[2021] \copyright 2021},
     PAGES = {xv+453},
      ISBN = {978-3-030-74247-8; 978-3-030-74248-5},
   MRCLASS = {14F22 (14E08 14G05 14G12 14K05)},
  MRNUMBER = {4304038},
MRREVIEWER = {Thomas\ Benedict\ Williams},
       DOI = {10.1007/978-3-030-74248-5},
       URL = {https://doi.org/10.1007/978-3-030-74248-5},
}

@article{MichaelStoll2001,
author = {Michael Stoll},
journal = {Acta Arithmetica},
language = {eng},
number = {3},
pages = {245-277},
title = {Implementing 2-descent for Jacobians of hyperelliptic curves},
url = {http://eudml.org/doc/279762},
volume = {98},
year = {2001},
}

@article {Mazur-Rubin:H10,
    AUTHOR = {Mazur, B. and Rubin, K.},
     TITLE = {Ranks of twists of elliptic curves and {H}ilbert's tenth
              problem},
   JOURNAL = {Invent. Math.},
  FJOURNAL = {Inventiones Mathematicae},
    VOLUME = {181},
      YEAR = {2010},
    NUMBER = {3},
     PAGES = {541--575},
      ISSN = {0020-9910,1432-1297},
   MRCLASS = {11G05 (11G40)},
  MRNUMBER = {2660452},
MRREVIEWER = {Pete\ L.\ Clark},
       DOI = {10.1007/s00222-010-0252-0},
       URL = {https://doi.org/10.1007/s00222-010-0252-0},
}

@article {Green_Tao:Linear_equations,
    AUTHOR = {Green, Benjamin and Tao, Terence},
     TITLE = {Linear equations in primes},
   JOURNAL = {Ann. of Math. (2)},
  FJOURNAL = {Annals of Mathematics. Second Series},
    VOLUME = {171},
      YEAR = {2010},
    NUMBER = {3},
     PAGES = {1753--1850},
      ISSN = {0003-486X,1939-8980},
   MRCLASS = {11N13 (11B30 11P32)},
  MRNUMBER = {2680398},
MRREVIEWER = {Tamar\ Ziegler},
       DOI = {10.4007/annals.2010.171.1753},
       URL = {https://doi.org/10.4007/annals.2010.171.1753},
}

@article {Gowers:Gowers_norms,
    AUTHOR = {Gowers, W. T.},
     TITLE = {A new proof of {S}zemer\'{e}di's theorem},
   JOURNAL = {Geom. Funct. Anal.},
  FJOURNAL = {Geometric and Functional Analysis},
    VOLUME = {11},
      YEAR = {2001},
    NUMBER = {3},
     PAGES = {465--588},
      ISSN = {1016-443X,1420-8970},
   MRCLASS = {11B25 (11K38 11K45)},
  MRNUMBER = {1844079},
MRREVIEWER = {Hillel\ Furstenberg},
       DOI = {10.1007/s00039-001-0332-9},
       URL = {https://doi.org/10.1007/s00039-001-0332-9},
}

@book {Tao:HOFA_book,
    AUTHOR = {Tao, Terence},
     TITLE = {Higher order {F}ourier analysis},
    SERIES = {Graduate Studies in Mathematics},
    VOLUME = {142},
 PUBLISHER = {American Mathematical Society, Providence, RI},
      YEAR = {2012},
     PAGES = {x+187},
      ISBN = {978-0-8218-8986-2},
   MRCLASS = {11B30 (11L07 11U07 37A45)},
  MRNUMBER = {2931680},
MRREVIEWER = {David\ Conlon},
       DOI = {10.1090/gsm/142},
       URL = {https://doi.org/10.1090/gsm/142},
}

@ARTICLE{Green_Sawhney,
       author = {{Green}, Ben and {Sawhney}, Mehtaab},
        title = "{Primes of the form $p^2 + nq^2$}",
      journal = {arXiv e-prints},
         year = 2024,
        month = oct,
          eid = {arXiv:2410.04189},
        pages = {arXiv:2410.04189},
          doi = {10.48550/arXiv.2410.04189},
archivePrefix = {arXiv},
       eprint = {2410.04189},
 primaryClass = {math.NT},
       adsurl = {https://ui.adsabs.harvard.edu/abs/2024arXiv241004189G}
}

@article{Leng,
      title={Efficient Equidistribution of Nilsequences}, 
      author={James Leng},
      year={2024},
        journal = {Preprint arXiv:2312.10772},
      eprint={2312.10772},
      archivePrefix={arXiv},
      primaryClass={math.NT},
      %url={https://arxiv.org/abs/2312.10772}, 
}

@article {TaoTeravainen,
    AUTHOR = {Tao, Terence and Ter\"{a}v\"{a}inen, Joni},
     TITLE = {Quantitative bounds for {G}owers uniformity of the {M}\"{o}bius
              and von {M}angoldt functions},
   JOURNAL = {J. Eur. Math. Soc. (JEMS)},
  FJOURNAL = {Journal of the European Mathematical Society (JEMS)},
    VOLUME = {27},
      YEAR = {2025},
    NUMBER = {4},
     PAGES = {1321--1384},
      ISSN = {1435-9855},
   MRCLASS = {11B30 (11N37 37A44)},
  MRNUMBER = {4875606},
       DOI = {10.4171/jems/1404},
       URL = {https://doi.org/10.4171/jems/1404},
}

@article{LengSahSawhney,
  title={{Quasipolynomial bounds on the inverse theorem for the Gowers $ U^{s+ 1}[N]$-norm}},
  author={Leng, James and Sah, Ashwin and Sawhney, Mehtaab},
  journal={Preprint arXiv:2402.17994},
  year={2024}
}

@article {KravitzKucaLeng,
    AUTHOR = {Kravitz, Noah and Kuca, Borys and Leng, James},
     TITLE = {Quantitative concatenation for polynomial box norms},
   JOURNAL = {Adv. Math.},
  FJOURNAL = {Advances in Mathematics},
    VOLUME = {489},
      YEAR = {2026},
     PAGES = {Paper No. 110820},
      ISSN = {0001-8708,1090-2082},
   MRCLASS = {11B30 (05)},
  MRNUMBER = {5022880},
       DOI = {10.1016/j.aim.2026.110820},
       URL = {https://doi-org.stanford.idm.oclc.org/10.1016/j.aim.2026.110820},
}

@article {FI:coordinates_of_Gaussian_primes,
    AUTHOR = {Friedlander, John and Iwaniec, Henryk},
     TITLE = {Coordinate distribution of {G}aussian primes},
   JOURNAL = {J. Eur. Math. Soc. (JEMS)},
  FJOURNAL = {Journal of the European Mathematical Society (JEMS)},
    VOLUME = {24},
      YEAR = {2022},
    NUMBER = {3},
     PAGES = {737--772},
      ISSN = {1435-9855,1435-9863},
   MRCLASS = {11N32 (11N36)},
  MRNUMBER = {4397031},
MRREVIEWER = {Karin\ Halupczok},
       DOI = {10.4171/JEMS/1083},
       URL = {https://doi.org/10.4171/JEMS/1083},
}

@article {Peluse:sets_wo_poly_prog,
    AUTHOR = {Peluse, Sarah},
     TITLE = {Bounds for sets with no polynomial progressions},
   JOURNAL = {Forum Math. Pi},
  FJOURNAL = {Forum of Mathematics. Pi},
    VOLUME = {8},
      YEAR = {2020},
     PAGES = {e16, 55},
      ISSN = {2050-5086},
   MRCLASS = {11B30 (11B25)},
  MRNUMBER = {4199235},
MRREVIEWER = {Ben\ Joseph\ Green},
       DOI = {10.1017/fmp.2020.11},
       URL = {https://doi.org/10.1017/fmp.2020.11},
}

@article {DFI,
    AUTHOR = {Duke, W. and Friedlander, J. B. and Iwaniec, H.},
     TITLE = {Equidistribution of roots of a quadratic congruence to prime
              moduli},
   JOURNAL = {Ann. of Math. (2)},
  FJOURNAL = {Annals of Mathematics. Second Series},
    VOLUME = {141},
      YEAR = {1995},
    NUMBER = {2},
     PAGES = {423--441},
      ISSN = {0003-486X,1939-8980},
   MRCLASS = {11N64 (11L20)},
  MRNUMBER = {1324141},
MRREVIEWER = {D.\ R.\ Heath-Brown},
       DOI = {10.2307/2118527},
       URL = {https://doi-org.stanford.idm.oclc.org/10.2307/2118527},
}

@article {HSW:RationalPoints,
    AUTHOR = {Harpaz, Yonatan and Skorobogatov, Alexei N. and Wittenberg,
              Olivier},
     TITLE = {The {H}ardy-{L}ittlewood conjecture and rational points},
   JOURNAL = {Compos. Math.},
  FJOURNAL = {Compositio Mathematica},
    VOLUME = {150},
      YEAR = {2014},
    NUMBER = {12},
     PAGES = {2095--2111},
      ISSN = {0010-437X,1570-5846},
   MRCLASS = {14G05 (11D57 14F22 14G25)},
  MRNUMBER = {3292295},
MRREVIEWER = {Yasuhiro\ Goto},
       DOI = {10.1112/S0010437X14007568},
       URL = {https://doi.org/10.1112/S0010437X14007568},
}

@misc{Zywina:Rank2,
      title={There are infinitely many elliptic curves over the rationals of rank 2}, 
      author={David Zywina},
      year={2025},
      eprint={2502.01957},
      archivePrefix={arXiv},
      primaryClass={math.NT},
}

@misc{Zywina:SmallRanks,
      title={On the infinitude of elliptic curves over a number field with prescribed small rank}, 
      author={David Zywina},
      year={2026},
      eprint={2602.10865},
      archivePrefix={arXiv},
      primaryClass={math.NT},
}

@misc{Zywina:Rank1_NumberFields,
      title={Rank one elliptic curves and rank stability}, 
      author={David Zywina},
      year={2025},
      eprint={2505.16960},
      archivePrefix={arXiv},
      primaryClass={math.NT},
}

@article {Rubin-Silverberg,
    AUTHOR = {Rubin, Karl and Silverberg, Alice},
     TITLE = {Rank frequencies for quadratic twists of elliptic curves},
   JOURNAL = {Experiment. Math.},
  FJOURNAL = {Experimental Mathematics},
    VOLUME = {10},
      YEAR = {2001},
    NUMBER = {4},
     PAGES = {559--569},
      ISSN = {1058-6458,1944-950X},
   MRCLASS = {11G05 (14G25)},
  MRNUMBER = {1881757},
MRREVIEWER = {Ivica\ Gusi\'c},
       URL = {http://projecteuclid.org/euclid.em/1069855256},
}

@misc{Smith:BSD-Goldfeld,
      title={The Birch and Swinnerton-Dyer conjecture implies Goldfeld's conjecture}, 
      author={Alexander Smith},
      year={2025},
      eprint={2503.17619v1},
      archivePrefix={arXiv},
      primaryClass={math.NT},
}

@article {DokchitserDokchitser:local-invariants-isogeneous,
    AUTHOR = {Dokchitser, Tim and Dokchitser, Vladimir},
     TITLE = {Local invariants of isogenous elliptic curves},
   JOURNAL = {Trans. Amer. Math. Soc.},
  FJOURNAL = {Transactions of the American Mathematical Society},
    VOLUME = {367},
      YEAR = {2015},
    NUMBER = {6},
     PAGES = {4339--4358},
      ISSN = {0002-9947,1088-6850},
   MRCLASS = {11G07 (11G05 11G40)},
  MRNUMBER = {3324930},
MRREVIEWER = {Mahesh\ Kakde},
       DOI = {10.1090/S0002-9947-2014-06271-5},
       URL = {https://doi.org/10.1090/S0002-9947-2014-06271-5},
}

@book {Silverman:Advanced,
    AUTHOR = {Silverman, Joseph H.},
     TITLE = {Advanced topics in the arithmetic of elliptic curves},
    SERIES = {Graduate Texts in Mathematics},
    VOLUME = {151},
 PUBLISHER = {Springer-Verlag, New York},
      YEAR = {1994},
     PAGES = {xiv+525},
      ISBN = {0-387-94328-5},
   MRCLASS = {11G05 (11G07 11G15 11G40 14H52)},
  MRNUMBER = {1312368},
MRREVIEWER = {Henri\ Darmon},
       DOI = {10.1007/978-1-4612-0851-8},
       URL = {https://doi.org/10.1007/978-1-4612-0851-8},
}

@book {Hasse:Hasse_principle,
    AUTHOR = {Hasse, Helmut},
     TITLE = {Mathematische {A}bhandlungen. {B}and 1},
      NOTE = {Herausgegeben von Heinrich Wolfgang Leopoldt und Peter
              Roquette},
 PUBLISHER = {Walter de Gruyter, Berlin-New York},
      YEAR = {1975},
     PAGES = {xv+535 pp. (1 plate)},
   MRCLASS = {01A75 (12-03)},
  MRNUMBER = {465756},
MRREVIEWER = {J.\ S.\ Joel},
}

@book {Serre:Galois_cohomology,
    AUTHOR = {Serre, Jean-Pierre},
     TITLE = {Galois cohomology},
    SERIES = {Springer Monographs in Mathematics},
   EDITION = {English},
      NOTE = {Translated from the French by Patrick Ion and revised by the
              author},
 PUBLISHER = {Springer-Verlag, Berlin},
      YEAR = {2002},
     PAGES = {x+210},
      ISBN = {3-540-42192-0},
   MRCLASS = {12G05 (11R34)},
  MRNUMBER = {1867431},
}

@article {GTZ:inverse_theorem,
    AUTHOR = {Green, Ben and Tao, Terence and Ziegler, Tamar},
     TITLE = {An inverse theorem for the {G}owers {$U^{s+1}[N]$}-norm},
   JOURNAL = {Ann. of Math. (2)},
  FJOURNAL = {Annals of Mathematics. Second Series},
    VOLUME = {176},
      YEAR = {2012},
    NUMBER = {2},
     PAGES = {1231--1372},
      ISSN = {0003-486X,1939-8980},
   MRCLASS = {11B30},
  MRNUMBER = {2950773},
MRREVIEWER = {Julia\ Wolf},
       DOI = {10.4007/annals.2012.176.2.11},
       URL = {https://doi.org/10.4007/annals.2012.176.2.11},
}

@article {Matthiesen:correlation_quad,
    AUTHOR = {Matthiesen, Lilian},
     TITLE = {Correlations of representation functions of binary quadratic
              forms},
   JOURNAL = {Acta Arith.},
  FJOURNAL = {Acta Arithmetica},
    VOLUME = {158},
      YEAR = {2013},
    NUMBER = {3},
     PAGES = {245--252},
      ISSN = {0065-1036,1730-6264},
   MRCLASS = {11N37 (11E25)},
  MRNUMBER = {3040665},
MRREVIEWER = {Gihan\ Marasingha},
       DOI = {10.4064/aa158-3-4},
       URL = {https://doi.org/10.4064/aa158-3-4},
}

@article {Matthiesen:correlation_pos_def_quad,
    AUTHOR = {Matthiesen, Lilian},
     TITLE = {Linear correlations amongst numbers represented by positive
              definite binary quadratic forms},
   JOURNAL = {Acta Arith.},
  FJOURNAL = {Acta Arithmetica},
    VOLUME = {154},
      YEAR = {2012},
    NUMBER = {3},
     PAGES = {235--306},
      ISSN = {0065-1036,1730-6264},
   MRCLASS = {11N37 (11B30 11E25)},
  MRNUMBER = {2949879},
MRREVIEWER = {Tom\ Sanders},
       DOI = {10.4064/aa154-3-2},
       URL = {https://doi.org/10.4064/aa154-3-2},
}

@article {DS:multi_BH,
    AUTHOR = {Destagnol, Kevin and Sofos, Efthymios},
     TITLE = {Rational points and prime values of polynomials in moderately
              many variables},
   JOURNAL = {Bull. Sci. Math.},
  FJOURNAL = {Bulletin des Sciences Math\'{e}matiques},
    VOLUME = {156},
      YEAR = {2019},
     PAGES = {102794, 33},
      ISSN = {0007-4497},
   MRCLASS = {11N32 (11P55 14G05)},
  MRNUMBER = {3996328},
MRREVIEWER = {Bingrong Huang},
       DOI = {10.1016/j.bulsci.2019.102794},
       URL = {https://doi.org/10.1016/j.bulsci.2019.102794},
}

@article {Matthiesen:mult-functions,
    AUTHOR = {Matthiesen, Lilian},
     TITLE = {Linear correlations of multiplicative functions},
   JOURNAL = {Proc. Lond. Math. Soc. (3)},
  FJOURNAL = {Proceedings of the London Mathematical Society. Third Series},
    VOLUME = {121},
      YEAR = {2020},
    NUMBER = {2},
     PAGES = {372--425},
      ISSN = {0024-6115,1460-244X},
   MRCLASS = {11N37 (11B30 11D04 11F30)},
  MRNUMBER = {4093960},
MRREVIEWER = {Olivier\ Bordell\`es},
       DOI = {10.1112/plms.12309},
       URL = {https://doi.org/10.1112/plms.12309},
}

@article {BatemanHorn,
    AUTHOR = {Bateman, Paul T. and Horn, Roger A.},
     TITLE = {A heuristic asymptotic formula concerning the distribution of
              prime numbers},
   JOURNAL = {Math. Comp.},
  FJOURNAL = {Mathematics of Computation},
    VOLUME = {16},
      YEAR = {1962},
     PAGES = {363--367},
      ISSN = {0025-5718},
   MRCLASS = {10.42},
  MRNUMBER = {148632},
MRREVIEWER = {A. Schinzel},
       DOI = {10.2307/2004056},
       URL = {https://doi.org/10.2307/2004056},
}

@article {GT3,
    AUTHOR = {Green, Ben and Tao, Terence},
     TITLE = {The {M}\"{o}bius function is strongly orthogonal to nilsequences},
   JOURNAL = {Ann. of Math. (2)},
  FJOURNAL = {Annals of Mathematics. Second Series},
    VOLUME = {175},
      YEAR = {2012},
    NUMBER = {2},
     PAGES = {541--566},
      ISSN = {0003-486X},
   MRCLASS = {37A45 (11A25)},
  MRNUMBER = {2877066},
MRREVIEWER = {Tamar Ziegler},
       DOI = {10.4007/annals.2012.175.2.3},
       URL = {https://doi.org/10.4007/annals.2012.175.2.3},
}

@article {GT4,
    AUTHOR = {Green, Ben and Tao, Terence},
     TITLE = {The quantitative behaviour of polynomial orbits on
              nilmanifolds},
   JOURNAL = {Ann. of Math. (2)},
  FJOURNAL = {Annals of Mathematics. Second Series},
    VOLUME = {175},
      YEAR = {2012},
    NUMBER = {2},
     PAGES = {465--540},
      ISSN = {0003-486X},
   MRCLASS = {37A15},
  MRNUMBER = {2877065},
MRREVIEWER = {Tamar Ziegler},
       DOI = {10.4007/annals.2012.175.2.2},
       URL = {https://doi.org/10.4007/annals.2012.175.2.2},
}

@article {GT2,
    AUTHOR = {Green, Benjamin and Tao, Terence},
     TITLE = {Linear equations in primes},
   JOURNAL = {Ann. of Math. (2)},
  FJOURNAL = {Annals of Mathematics. Second Series},
    VOLUME = {171},
      YEAR = {2010},
    NUMBER = {3},
     PAGES = {1753--1850},
      ISSN = {0003-486X},
   MRCLASS = {11N13 (11B30 11P32)},
  MRNUMBER = {2680398},
MRREVIEWER = {Tamar Ziegler},
       DOI = {10.4007/annals.2010.171.1753},
       URL = {https://doi.org/10.4007/annals.2010.171.1753},
}

@article {Birch,
    AUTHOR = {Birch, B. J.},
     TITLE = {Forms in many variables},
   JOURNAL = {Proc. Roy. Soc. London Ser. A},
  FJOURNAL = {Proceedings of the Royal Society. London. Series A.
              Mathematical, Physical and Engineering Sciences},
    VOLUME = {265},
      YEAR = {1961/62},
     PAGES = {245--263},
      ISSN = {0962-8444},
   MRCLASS = {10.15},
  MRNUMBER = {150129},
MRREVIEWER = {John V. Armitage},
       DOI = {10.1098/rspa.1962.0007},
       URL = {https://doi.org/10.1098/rspa.1962.0007},
}

@article {FriedlanderIwaniec24,
    AUTHOR = {Friedlander, John and Iwaniec, Henryk},
     TITLE = {The polynomial {$X^2+Y^4$} captures its primes},
   JOURNAL = {Ann. of Math. (2)},
  FJOURNAL = {Annals of Mathematics. Second Series},
    VOLUME = {148},
      YEAR = {1998},
    NUMBER = {3},
     PAGES = {945--1040},
      ISSN = {0003-486X},
   MRCLASS = {11N05 (11N32 11N35 11N36)},
  MRNUMBER = {1670065},
MRREVIEWER = {Andrew Granville},
       DOI = {10.2307/121034},
       URL = {https://doi.org/10.2307/121034},
}

@article {HeathBrownx32y3,
    AUTHOR = {Heath-Brown, D. R.},
     TITLE = {Primes represented by {$x^3+2y^3$}},
   JOURNAL = {Acta Math.},
  FJOURNAL = {Acta Mathematica},
    VOLUME = {186},
      YEAR = {2001},
    NUMBER = {1},
     PAGES = {1--84},
      ISSN = {0001-5962},
   MRCLASS = {11N32 (11N36)},
  MRNUMBER = {1828372},
MRREVIEWER = {G. Greaves},
       DOI = {10.1007/BF02392715},
       URL = {https://doi.org/10.1007/BF02392715},
}

@article {HBMoroz,
    AUTHOR = {Heath-Brown, D. R. and Moroz, B. Z.},
     TITLE = {On the representation of primes by cubic polynomials in two
              variables},
   JOURNAL = {Proc. London Math. Soc. (3)},
  FJOURNAL = {Proceedings of the London Mathematical Society. Third Series},
    VOLUME = {88},
      YEAR = {2004},
    NUMBER = {2},
     PAGES = {289--312},
      ISSN = {0024-6115},
   MRCLASS = {11N32 (11N36 11R44)},
  MRNUMBER = {2032509},
MRREVIEWER = {Koichi Kawada},
       DOI = {10.1112/S0024611503014497},
       URL = {https://doi.org/10.1112/S0024611503014497},
}

@article {MaynardNormForm,
    AUTHOR = {Maynard, James},
     TITLE = {Primes represented by incomplete norm forms},
   JOURNAL = {Forum Math. Pi},
  FJOURNAL = {Forum of Mathematics. Pi},
    VOLUME = {8},
      YEAR = {2020},
     PAGES = {e3, 128},
   MRCLASS = {11N05 (11D57 11H55 11N32 11N35)},
  MRNUMBER = {4061964},
MRREVIEWER = {Tsz Ho Chan},
       DOI = {10.1017/fmp.2019.8},
       URL = {https://doi.org/10.1017/fmp.2019.8},
}

@ARTICLE{KotsovolisWoo,
       AUTHOR = {Kotsovolis, Giorgos and Woo, Katharine},
     TITLE = {Prime number theorems for polynomials from homogeneous
              dynamics},
   JOURNAL = {Geom. Funct. Anal.},
  FJOURNAL = {Geometric and Functional Analysis},
    VOLUME = {35},
      YEAR = {2025},
    NUMBER = {4},
     PAGES = {1169--1221},
      ISSN = {1016-443X,1420-8970},
   MRCLASS = {11N25 (37A17)},
  MRNUMBER = {4956619},
       DOI = {10.1007/s00039-025-00716-y},
       URL = {https://doi.org/10.1007/s00039-025-00716-y},
}

@article {GT1,
    AUTHOR = {Green, Ben and Tao, Terence},
     TITLE = {The primes contain arbitrarily long arithmetic progressions},
   JOURNAL = {Ann. of Math. (2)},
  FJOURNAL = {Annals of Mathematics. Second Series},
    VOLUME = {167},
      YEAR = {2008},
    NUMBER = {2},
     PAGES = {481--547},
      ISSN = {0003-486X,1939-8980},
   MRCLASS = {11N13 (11A41 11B25 37A45)},
  MRNUMBER = {2415379},
MRREVIEWER = {Tamar\ Ziegler},
       DOI = {10.4007/annals.2008.167.481},
       URL = {https://doi.org/10.4007/annals.2008.167.481},
}

@article {Vinogradov,
    AUTHOR = {Vinogradov, I.M.},
     TITLE = {A new method in analytic number theory (Russian)},
   JOURNAL = {Tr. Mat. Inst. Steklova},
      YEAR = {1937},
    NUMBER = {10},
     PAGES = {5–122},
}

@ARTICLE{ABHS,
       AUTHOR = {Alp\"oge, Levent and Bhargava, Manjul and Ho, Wei and
              Shnidman, Ari},
     TITLE = {Rank stability in quadratic extensions and {H}ilbert's tenth
              problem for the ring of integers of a number field},
   JOURNAL = {Invent. Math.},
  FJOURNAL = {Inventiones Mathematicae},
    VOLUME = {243},
      YEAR = {2026},
    NUMBER = {3},
     PAGES = {1129--1139},
      ISSN = {0020-9910,1432-1297},
   MRCLASS = {11G10 (11G30 11U05 14K15)},
  MRNUMBER = {5008163},
       DOI = {10.1007/s00222-025-01392-3},
       URL = {https://doi.org/10.1007/s00222-025-01392-3},
}

@article {TZ:polynomial-progressions,
    AUTHOR = {Tao, Terence and Ziegler, Tamar},
     TITLE = {The primes contain arbitrarily long polynomial progressions},
   JOURNAL = {Acta Math.},
  FJOURNAL = {Acta Mathematica},
    VOLUME = {201},
      YEAR = {2008},
    NUMBER = {2},
     PAGES = {213--305},
      ISSN = {0001-5962,1871-2509},
   MRCLASS = {11N13 (11B25 37A45)},
  MRNUMBER = {2461509},
MRREVIEWER = {Ben\ Joseph\ Green},
       DOI = {10.1007/s11511-008-0032-5},
       URL = {https://doi.org/10.1007/s11511-008-0032-5},
}

@article {FouvryIwaniec,
    AUTHOR = {Fouvry, Etienne and Iwaniec, Henryk},
     TITLE = {Gaussian primes},
   JOURNAL = {Acta Arith.},
  FJOURNAL = {Acta Arithmetica},
    VOLUME = {79},
      YEAR = {1997},
    NUMBER = {3},
     PAGES = {249--287},
      ISSN = {0065-1036,1730-6264},
   MRCLASS = {11N05 (11N32 11N35 11N36)},
  MRNUMBER = {1438827},
MRREVIEWER = {John\ B.\ Friedlander},
       DOI = {10.4064/aa-79-3-249-287},
       URL = {https://doi.org/10.4064/aa-79-3-249-287},
}

@article {HBLi,
    AUTHOR = {Heath-Brown, D. R. and Li, Xiannan},
     TITLE = {Prime values of {$a^2 + p^4$}},
   JOURNAL = {Invent. Math.},
  FJOURNAL = {Inventiones Mathematicae},
    VOLUME = {208},
      YEAR = {2017},
    NUMBER = {2},
     PAGES = {441--499},
      ISSN = {0020-9910,1432-1297},
   MRCLASS = {11M06 (11M26)},
  MRNUMBER = {3639597},
MRREVIEWER = {Xiumin\ Ren},
       DOI = {10.1007/s00222-016-0694-0},
       URL = {https://doi.org/10.1007/s00222-016-0694-0},
}

@article {LSX:quadratic-prime-var,
    AUTHOR = {Lam, Peter Cho-Ho and Schindler, Damaris and Xiao, Stanley
              Yao},
     TITLE = {On prime values of binary quadratic forms with a thin
              variable},
   JOURNAL = {J. Lond. Math. Soc. (2)},
  FJOURNAL = {Journal of the London Mathematical Society. Second Series},
    VOLUME = {102},
      YEAR = {2020},
    NUMBER = {2},
     PAGES = {749--772},
      ISSN = {0024-6107,1469-7750},
   MRCLASS = {11N32 (11E16 11N05 11N35 11N36)},
  MRNUMBER = {4171433},
MRREVIEWER = {Tsz\ Ho\ Chan},
       DOI = {10.1112/jlms.12336},
       URL = {https://doi.org/10.1112/jlms.12336},
}

@article {Merikoski,
    AUTHOR = {Merikoski, Jori},
     TITLE = {The polynomials {$X^2+(Y^2+1)^2$} and {$X^2+(Y^3+Z^3)^2$} also
              capture their primes},
   JOURNAL = {Proc. Lond. Math. Soc. (3)},
  FJOURNAL = {Proceedings of the London Mathematical Society. Third Series},
    VOLUME = {127},
      YEAR = {2023},
    NUMBER = {4},
     PAGES = {1057--1133},
      ISSN = {0024-6115,1460-244X},
   MRCLASS = {11N32 (11N36)},
  MRNUMBER = {4655347},
MRREVIEWER = {A.\ Perelli},
       DOI = {10.1112/plms.12557},
       URL = {https://doi.org/10.1112/plms.12557},
}

@article {Xiao:prime-quadratic,
    AUTHOR = {Xiao, Stanley Yao},
     TITLE = {Prime values of {$f(a,b^2)$} and {$f(a,p^2)$}, {$f$}
              quadratic},
   JOURNAL = {Algebra Number Theory},
  FJOURNAL = {Algebra \& Number Theory},
    VOLUME = {18},
      YEAR = {2024},
    NUMBER = {9},
     PAGES = {1619--1679},
      ISSN = {1937-0652,1944-7833},
   MRCLASS = {11N32 (11N35 11N36 11R45)},
  MRNUMBER = {4856605},
MRREVIEWER = {Katharine\ Woo},
       DOI = {10.2140/ant.2024.18.1619},
       URL = {https://doi.org/10.2140/ant.2024.18.1619},
}

@article {CTSSD:generalize-CTSD,
    AUTHOR = {Colliot-Th\'{e}l\`ene, J.-L. and Skorobogatov, A. N. and
              Swinnerton-Dyer, Peter},
     TITLE = {Rational points and zero-cycles on fibred varieties:
              {S}chinzel's hypothesis and {S}alberger's device},
   JOURNAL = {J. Reine Angew. Math.},
  FJOURNAL = {Journal f\"{u}r die Reine und Angewandte Mathematik. [Crelle's
              Journal]},
    VOLUME = {495},
      YEAR = {1998},
     PAGES = {1--28},
      ISSN = {0075-4102,1435-5345},
   MRCLASS = {14G25 (11G35 14C25 14D10)},
  MRNUMBER = {1603908},
MRREVIEWER = {R.\ T.\ Hoobler},
       DOI = {10.1515/crll.1998.019},
       URL = {https://doi.org/10.1515/crll.1998.019},
}

@article {SD:pencils-of-conics,
    AUTHOR = {Swinnerton-Dyer, Peter},
     TITLE = {Rational points on pencils of conics and on pencils of
              quadrics},
   JOURNAL = {J. London Math. Soc. (2)},
  FJOURNAL = {Journal of the London Mathematical Society. Second Series},
    VOLUME = {50},
      YEAR = {1994},
    NUMBER = {2},
     PAGES = {231--242},
      ISSN = {0024-6107,1469-7750},
   MRCLASS = {14G05 (11G35 14C21 14G25)},
  MRNUMBER = {1291734},
MRREVIEWER = {Philippe\ Satg\'{e}},
       DOI = {10.1112/jlms/50.2.231},
       URL = {https://doi.org/10.1112/jlms/50.2.231},
}

@misc{Savoie:rank2,
      title={Infinitely many elliptic curves over $\mathbb{Q}(i)$ with rank 2 and $j$-invariant 1728}, 
      author={Ben Savoie},
      year={2025},
      eprint={2506.17605v2},
      archivePrefix={arXiv},
      primaryClass={math.NT},
}

@article {Tao:GaussianPrimesConstellations,
    AUTHOR = {Tao, Terence},
     TITLE = {The {G}aussian primes contain arbitrarily shaped
              constellations},
   JOURNAL = {J. Anal. Math.},
  FJOURNAL = {Journal d'Analyse Math\'ematique},
    VOLUME = {99},
      YEAR = {2006},
     PAGES = {109--176},
      ISSN = {0021-7670,1565-8538},
   MRCLASS = {11B25 (11B75 37A45)},
  MRNUMBER = {2279549},
MRREVIEWER = {Ben\ Joseph\ Green},
       DOI = {10.1007/BF02789444},
       URL = {https://doi.org/10.1007/BF02789444},
}

@misc{Elkies:rank1819,
      title={An elliptic K3 surface X/Q(t) with Mordell-Weil rank 17, I: Formulas for X and base changes of rank 18 and 19}, 
      author={Noam D. Elkies},
      year={2026},
      eprint={2608.25406},
      archivePrefix={arXiv},
      primaryClass={math.NT},
      %url={https://arxiv.org/abs/2608.25406}, 
}

@misc{Elkies:lectures,
      title={Three lectures on elliptic surfaces and curves of high rank}, 
      author={Noam D. Elkies},
      year={2007},
      eprint={0709.2908},
      archivePrefix={arXiv},
      primaryClass={math.NT},
      %url={https://arxiv.org/abs/0709.2908}, 
}

@article {RS:ranks-of-ec,
    AUTHOR = {Rubin, Karl and Silverberg, Alice},
     TITLE = {Ranks of elliptic curves},
   JOURNAL = {Bull. Amer. Math. Soc. (N.S.)},
  FJOURNAL = {American Mathematical Society. Bulletin. New Series},
    VOLUME = {39},
      YEAR = {2002},
    NUMBER = {4},
     PAGES = {455--474},
      ISSN = {0273-0979,1088-9485},
   MRCLASS = {11G05 (11G40)},
  MRNUMBER = {1920278},
MRREVIEWER = {Ioannis\ Papadopoulos},
       DOI = {10.1090/S0273-0979-02-00952-7},
       URL = {https://doi.org/10.1090/S0273-0979-02-00952-7},
}

@article {Tao:log-Chowla,
    AUTHOR = {Tao, Terence},
     TITLE = {The logarithmically averaged {C}howla and {E}lliott
              conjectures for two-point correlations},
   JOURNAL = {Forum Math. Pi},
  FJOURNAL = {Forum of Mathematics. Pi},
    VOLUME = {4},
      YEAR = {2016},
     PAGES = {e8, 36},
      ISSN = {2050-5086},
   MRCLASS = {11N37},
  MRNUMBER = {3569059},
MRREVIEWER = {Y.-F.\ S.\ P\'{e}termann},
       DOI = {10.1017/fmp.2016.6},
       URL = {https://doi.org/10.1017/fmp.2016.6},
}

@article {Satge:calcul-de-Heegner,
    AUTHOR = {Satg\'e, Philippe},
     TITLE = {Un analogue du calcul de {H}eegner},
   JOURNAL = {Invent. Math.},
  FJOURNAL = {Inventiones Mathematicae},
    VOLUME = {87},
      YEAR = {1987},
    NUMBER = {2},
     PAGES = {425--439},
      ISSN = {0020-9910,1432-1297},
   MRCLASS = {11G05 (11D25 11G18)},
  MRNUMBER = {870738},
MRREVIEWER = {Andrew\ Bremner},
       DOI = {10.1007/BF01389425},
       URL = {https://doi.org/10.1007/BF01389425},
}

@article {Heegner,
    AUTHOR = {Heegner, Kurt},
     TITLE = {Diophantische {A}nalysis und {M}odulfunktionen},
   JOURNAL = {Math. Z.},
  FJOURNAL = {Mathematische Zeitschrift},
    VOLUME = {56},
      YEAR = {1952},
     PAGES = {227--253},
      ISSN = {0025-5874,1432-1823},
   MRCLASS = {10.0X},
  MRNUMBER = {53135},
MRREVIEWER = {M.\ Ward},
       DOI = {10.1007/BF01174749},
       URL = {https://doi.org/10.1007/BF01174749},
}

@article {Silverman:heights-specialization,
    AUTHOR = {Silverman, Joseph H.},
     TITLE = {Heights and the specialization map for families of abelian
              varieties},
   JOURNAL = {J. Reine Angew. Math.},
  FJOURNAL = {Journal f\"ur die Reine und Angewandte Mathematik. [Crelle's
              Journal]},
    VOLUME = {342},
      YEAR = {1983},
     PAGES = {197--211},
      ISSN = {0075-4102,1435-5345},
   MRCLASS = {14K15 (14D10 14G25)},
  MRNUMBER = {703488},
MRREVIEWER = {Gerd\ Faltings},
       DOI = {10.1515/crll.1983.342.197},
       URL = {https://doi.org/10.1515/crll.1983.342.197},
}

@article {Silverman:divisilibity-specialization,
    AUTHOR = {Silverman, Joseph H.},
     TITLE = {Divisibility of the specialization map for families of
              elliptic curves},
   JOURNAL = {Amer. J. Math.},
  FJOURNAL = {American Journal of Mathematics},
    VOLUME = {107},
      YEAR = {1985},
    NUMBER = {3},
     PAGES = {555--565},
      ISSN = {0002-9327,1080-6377},
   MRCLASS = {14K15 (11G05 14K07)},
  MRNUMBER = {789655},
MRREVIEWER = {Ernst\ Kani},
       DOI = {10.2307/2374369},
       URL = {https://doi.org/10.2307/2374369},
}

@misc{hexagon:2610.00200v2,
  title={A quasi-Riemann hypothesis for Dedekind zeta functions},
  key={hexagon:2610.00200v2},
  howpublished={Hexagon preprint hexagon:2610.00200v2},
  year={2026},
  note={Contributors: Dean Cureton; Jesse Thorner},
  url={https://hexagonmath.org/2610.00200v2},
  author={Cureton, Dean and Thorner, Jesse},
}

@misc{OAI:Uniform-exclusion-of-Landau-Siegel-zeros-October-1-2026,
  author = {{OpenAI}},
  title = {{Uniform exclusion of Landau--Siegel zeros}},
  howpublished = {OpenAI Math Release preprint
                  \href{https://github.com/openai/math/blob/main/preprints/Uniform-exclusion-of-Landau-Siegel-zeros-October-1-2026/paper.pdf}{OAI:Uniform-exclusion-of-Landau-Siegel-zeros-October-1-2026}},
  year = {2026}
}

@misc{OAI:The-Quasi-Riemann-Hypothesis-October-5-2026,
  author = {{OpenAI}},
  title = {{The Quasi-Riemann Hypothesis}},
  howpublished = {OpenAI Math Release preprint
                  \href{https://github.com/openai/math/blob/main/preprints/The-Quasi-Riemann-Hypothesis-October-5-2026/paper2.pdf}{OAI:The-Quasi-Riemann-Hypothesis-October-5-2026}},
  year = {2026}
}

@article {HW:fibration-method,
    AUTHOR = {Harpaz, Yonatan and Wittenberg, Olivier},
     TITLE = {On the fibration method for zero-cycles and rational points},
   JOURNAL = {Ann. of Math. (2)},
  FJOURNAL = {Annals of Mathematics. Second Series},
    VOLUME = {183},
      YEAR = {2016},
    NUMBER = {1},
     PAGES = {229--295},
      ISSN = {0003-486X,1939-8980},
   MRCLASS = {14G25 (14C15 14M22)},
  MRNUMBER = {3432584},
MRREVIEWER = {Amanda\ Knecht},
       DOI = {10.4007/annals.2016.183.1.5},
       URL = {https://doi.org/10.4007/annals.2016.183.1.5},
}

@article {FH:HOFA-mult-func,
    AUTHOR = {Frantzikinakis, Nikos and Host, Bernard},
     TITLE = {Higher order {F}ourier analysis of multiplicative functions
              and applications},
   JOURNAL = {J. Amer. Math. Soc.},
  FJOURNAL = {Journal of the American Mathematical Society},
    VOLUME = {30},
      YEAR = {2017},
    NUMBER = {1},
     PAGES = {67--157},
      ISSN = {0894-0347,1088-6834},
   MRCLASS = {11N37 (05D10 11B30 11N60 37A45)},
  MRNUMBER = {3556289},
MRREVIEWER = {Vilius\ Stakenas},
       DOI = {10.1090/jams/857},
       URL = {https://doi.org/10.1090/jams/857},
}

@article {Kuca,
    AUTHOR = {Kuca, Borys},
     TITLE = {Multidimensional polynomial patterns over finite fields:
              bounds, counting estimates and {G}owers norm control},
   JOURNAL = {Adv. Math.},
  FJOURNAL = {Advances in Mathematics},
    VOLUME = {448},
      YEAR = {2024},
     PAGES = {Paper No. 109700, 61},
      ISSN = {0001-8708,1090-2082},
   MRCLASS = {11B30},
  MRNUMBER = {4741361},
MRREVIEWER = {Ben\ Joseph\ Green},
       DOI = {10.1016/j.aim.2024.109700},
       URL = {https://doi.org/10.1016/j.aim.2024.109700},
}

@article {PP:nonlinear-Roth,
    AUTHOR = {Peluse, Sarah and Prendiville, Sean},
     TITLE = {Quantitative bounds in the nonlinear {R}oth theorem},
   JOURNAL = {Invent. Math.},
  FJOURNAL = {Inventiones Mathematicae},
    VOLUME = {238},
      YEAR = {2024},
    NUMBER = {3},
     PAGES = {865--903},
      ISSN = {0020-9910,1432-1297},
   MRCLASS = {11B30 (37A30)},
  MRNUMBER = {4824730},
MRREVIEWER = {Wolfgang\ A.\ Schmid},
       DOI = {10.1007/s00222-024-01293-x},
       URL = {https://doi.org/10.1007/s00222-024-01293-x},
}

@article {HWW:fibration,
    AUTHOR = {Harpaz, Yonatan and Wei, Dasheng and Wittenberg, Olivier},
     TITLE = {Rational points on fibrations with few non-split fibres},
   JOURNAL = {J. Reine Angew. Math.},
  FJOURNAL = {Journal f\"{u}r die Reine und Angewandte Mathematik. [Crelle's
              Journal]},
    VOLUME = {791},
      YEAR = {2022},
     PAGES = {89--133},
      ISSN = {0075-4102,1435-5345},
   MRCLASS = {14G12 (11G35)},
  MRNUMBER = {4489626},
MRREVIEWER = {Mikhail\ Borovoi},
       DOI = {10.1515/crelle-2022-0042},
       URL = {https://doi.org/10.1515/crelle-2022-0042},
}

@article {Matthiesen:norm,
    AUTHOR = {Matthiesen, Lilian},
     TITLE = {On the square-free representation function of a norm form and
              nilsequences},
   JOURNAL = {J. Inst. Math. Jussieu},
  FJOURNAL = {Journal of the Institute of Mathematics of Jussieu. JIMJ.
              Journal de l'Institut de Math\'{e}matiques de Jussieu},
    VOLUME = {17},
      YEAR = {2018},
    NUMBER = {1},
     PAGES = {107--135},
      ISSN = {1474-7480,1475-3030},
   MRCLASS = {11B30 (11D57 11N37 14G05)},
  MRNUMBER = {3742557},
MRREVIEWER = {Sean\ M.\ Prendiville},
       DOI = {10.1017/S1474748015000389},
       URL = {https://doi.org/10.1017/S1474748015000389},
}

@incollection {SD:Applications-of-Schinzel,
    AUTHOR = {Swinnerton-Dyer, Peter},
     TITLE = {Some applications of {S}chinzel's hypothesis to {D}iophantine
              equations},
 BOOKTITLE = {Number theory in progress, {V}ol. 1
              ({Z}akopane-{K}o\'scielisko, 1997)},
     PAGES = {503--530},
 PUBLISHER = {de Gruyter, Berlin},
      YEAR = {1999},
      ISBN = {3-11-015715-2},
   MRCLASS = {11G35 (14D10 14J28)},
  MRNUMBER = {1689527},
MRREVIEWER = {Tam\'as\ Szamuely},
}

@article {BM:norm-forms,
    AUTHOR = {Browning, Tim D. and Matthiesen, Lilian},
     TITLE = {Norm forms for arbitrary number fields as products of linear
              polynomials},
   JOURNAL = {Ann. Sci. \'{E}c. Norm. Sup\'{e}r. (4)},
  FJOURNAL = {Annales Scientifiques de l'\'{E}cole Normale Sup\'{e}rieure.
              Quatri\`eme S\'{e}rie},
    VOLUME = {50},
      YEAR = {2017},
    NUMBER = {6},
     PAGES = {1383--1446},
      ISSN = {0012-9593,1873-2151},
   MRCLASS = {11D57 (14G05)},
  MRNUMBER = {3742196},
MRREVIEWER = {P.\ Bundschuh},
       DOI = {10.24033/asens.2348},
       URL = {https://doi.org/10.24033/asens.2348},
}

@article {NS:saturation,
    AUTHOR = {Nevo, Amos and Sarnak, Peter},
     TITLE = {Prime and almost prime integral points on principal
              homogeneous spaces},
   JOURNAL = {Acta Math.},
  FJOURNAL = {Acta Mathematica},
    VOLUME = {205},
      YEAR = {2010},
    NUMBER = {2},
     PAGES = {361--402},
      ISSN = {0001-5962,1871-2509},
   MRCLASS = {22F30 (11N36 11N45 20G30 22E40)},
  MRNUMBER = {2746350},
MRREVIEWER = {B.\ Sury},
       DOI = {10.1007/s11511-010-0057-4},
       URL = {https://doi.org/10.1007/s11511-010-0057-4},
}

@incollection {Sarnak:eqi-and-primes,
    AUTHOR = {Sarnak, Peter},
     TITLE = {Equidistribution and primes},
      NOTE = {G\'eom\'etrie diff\'erentielle, physique math\'ematique,
              math\'ematiques et soci\'et\'e. II},
   JOURNAL = {Ast\'erisque},
  FJOURNAL = {Ast\'erisque},
    NUMBER = {322},
      YEAR = {2008},
     PAGES = {225--240},
      ISSN = {0303-1179,2492-5926},
      ISBN = {978-285629-259-4},
   MRCLASS = {11N05 (11N36 11N75 52C15)},
  MRNUMBER = {2521658},
}

@article {Marasingha:saturation,
    AUTHOR = {Marasingha, Gihan},
     TITLE = {Almost primes represented by binary forms},
   JOURNAL = {J. Lond. Math. Soc. (2)},
  FJOURNAL = {Journal of the London Mathematical Society. Second Series},
    VOLUME = {82},
      YEAR = {2010},
    NUMBER = {2},
     PAGES = {295--316},
      ISSN = {0024-6107,1469-7750},
   MRCLASS = {11N36 (11E76 11N05 11N32 11R42)},
  MRNUMBER = {2725041},
MRREVIEWER = {S.\ W.\ Graham},
       DOI = {10.1112/jlms/jdq026},
       URL = {https://doi.org/10.1112/jlms/jdq026},
}

@inbook{DH, place={Cambridge}, series={London Mathematical Society Lecture Note Series}, title={Some Applications of Sieves of Dimension exceeding 1}, booktitle={Sieve Methods, Exponential Sums, and their Applications in Number Theory}, publisher={Cambridge University Press}, author={Diamond, H. and Halberstam, H. and Greaves, G. R. H. and Harman, G. and Huxley, M. N.}, year={1997}, pages={101–108}, collection={London Mathematical Society Lecture Note Series}}

\end{document}